\documentclass[a4paper, leqno, 11pt]{amsart}  
\usepackage{amsfonts}
\usepackage{amssymb}
\usepackage{graphicx}
\usepackage[applemac]{inputenc}
\usepackage{color}
\usepackage{hyperref}
\usepackage{mathdots}
\usepackage{float}
\newtheorem*{TheoremA}{Theorem A}
\newtheorem*{CorollaryA}{Corollary A}
\newtheorem{lemma}{Lemma}[section]

\newtheorem{definition}[lemma]{Definition}
\newtheorem{proposition}[lemma]{Proposition}
\newtheorem{theorem}[lemma]{Theorem}
\newtheorem{remark}[lemma]{Remark}

\newtheorem{corollary}[lemma]{Corollary}

\def\CCC{\mathbb{C}}

\def\ZZZ{\mathbb{Z}}

\def\RRR{\mathbb{R}}
\def\TTT{\mathbb{T}}

\def\NNN{\mathbb{N}}
\def\PPP{\mathbb{P}}

\def\aaa{\mathtt{a}}
\def\eee{\mathtt{e}}

\def\hhh{\mathbf{h}}

\def\sss{\mathtt{s}}
\def\uuu{\mathtt{u}}
\def\vvv{\mathtt{v}}
\def\www{\mathtt{w}}
\def\xxx{\mathtt{x}}

\def\gotA{\mathfrak{A}}

\def\gotL{\mathfrak{L}}

\def\nucirc{\raise 0pt\hbox{$\mathop{\#\mathcal{I}}\limits^{_\circ\!}$}}
\def\etacirc{\raise -1.5pt\hbox{$\mathop{\eta}\limits^{_\circ\!}$}}

\def\emptyword{{\circ\!\!\!/}}

\def\aa{{\rm a}}

\def\ee{{\bf e}}

\def\hh{{\rm h}}

\def\kk{{\bf k}}
\def\ww{{\bf m}}

\def\ss{{\bf s}}

\def\ww{{\bf w}}

\def\yy{{\bf y}}

\def\YY{{\bf Y}}
\def\KK{{\bf K}}
\def\XX{\mathbf{X}}

\def\II{\mathbf{I}}
\def\JJ{{\bf J}}
\def\ZZ{{\bf Z}}

\def\VV{{\bf V}}

\def\bow#1{\mathop{\approx}\limits_n}

\def\1{\mathtt{1}}
\def\0{\mathtt{0}}
\def\2{\mathtt{2}}
\def\3{\mathtt{3}}
\def\4{\mathtt{4}}
\def\5{\mathtt{5}}

\def\[{[\![}
\def\]{]\!]}

\title[Projective convergence of columns for inhomogeneous matrix products]{Projective convergence of columns for inhomogeneous products of matrices with nonnegative entries}
\author[E. Olivier]{\'Eric Olivier}
\address{\'Eric Olivier,
GDAC-I2M UMR 7373 CNRS  Universit\'e d'Aix-Marseille, 3, place Victor Hugo,\hfil\break
13331 Marseille Cedex 03, France}
\email{eric.olivier@univ-amu.fr}
\author[A. Thomas]{Alain Thomas}
\address{Alain Thomas,
448 all\'ee des Cantons, 83640 Plan d'Aups Sainte Baume, France}
\email{alain-yves.thomas@laposte.net}
\keywords{Infinite products of nonnegative matrices, Multifractal analysis, Bernoulli convolutions}
\subjclass{15B48, 28A12}
\date{}
  
\begin{document}

\baselineskip=15pt

\begin{abstract}

Let $P_n$ be the $n$-step right product $A_1\cdots A_n$, where $A_1,A_2,\dots$ is a given infinite sequence of  $d\times d$ matrices with nonnegative entries. In a wide range of situations, the normalized matrix product $P_n/{\Vert P_n\Vert}$ does not converge and we shall be rather interested in the asymptotic behavior of the normalized columns $P_nU_i/\Vert P_nU_i\Vert$, where $U_1,\dots,U_d$ are the canonical $d\times 1$ vectors. 
Our main result in Theorem~A gives a sufficient condition  ${\bf (C)}$ over the sequence $A_1,A_2,\dots$ ensuring the existence  of {\it dominant columns} of $P_n$, having the same projective limit $V$: more precisely, for any rank $n$, there exists a partition of  $\{1,\dots,d\}$ made of two subsets $J_n\ne\emptyset$ and $J_n^c$ such that each one of the sequences of normalized columns, say  $P_nU_{j_n}/\Vert P_nU_{j_n}\Vert$ with $j_n\in J_n$ tends to $V$ as $n$ tends to $+\infty$ and are {\it dominant} in the sense that the ratio $\Vert P_nU_{j_n'}/\Vert P_nU_{j_n}\Vert$ tends to $0$, as soon as  $j_n'\in J_n^c$.  The existence of sequences of such {\it dominant columns} implies that for any probability vector $X$ with positive entries, the probability vector $P_nX/\Vert P_nX\Vert$, converges as $n$ tends to $+\infty$. 
Our main application of Theorem~A (and our initial motivation) is related to an {\it Erd\H os problem} concerned with a family of probability measures $\mu_\beta$ (for $1<\beta<2$ a real parameter) fully supported by a subinterval of the real line, known as {\it Bernoulli convolutions}. For some parameters $\beta$ (actually the so-called PV-numbers) such  measures are known to be {\it linearly representable}: the $\mu_\beta$-measure of a suitable family of nested generating intervals may be computed by means of matrix products of the form $P_nX$, where $A_n$ takes only finitely many values, say $A(\0),\dots,A(\aaa)$,  and $X$ is a probability vector with positive entries. Because, $A_n=A(\xi_n)$, where $\xi=\xi_1\xi_2\cdots$ is a sequence (one-sided infinite word) with  $\xi_n\in\{\0,\dots,\aaa\}$,  we shall write $P_n=P_n(\xi)$ the dependence of the $n$-step product with~$\xi$: when the convergence of ${P_n(\xi)X}/{\Vert P_n(\xi)X\Vert}$ is uniform w.r.t. $\xi$, a sharp  analysis of the measure $\mu_\beta$ (Gibbs structures and multifractal decomposition) becomes possible. However, most of the matrices involved in the decomposition of these Bernoulli convolutions are large, sparse and it is usually not easy to prove the condition ${\bf (C)}$ of Theorem~A. To illustrate the technics, we  consider one parameter $\beta$  for which the matrices are neither too small nor too large and we verify condition ${\bf (C)}$: this leads to the Gibbs properties of $\mu_\beta$.

\end{abstract}

\maketitle

\vskip20pt

\section{\bf Introduction}

\smallskip

\baselineskip=15pt

\subsection{Generalities}Given $\mathcal{A}=(A_1,A_2,\dots)$ an infinite sequence of $d\times d$ matrices we consider the right products
$
P_n=A_1\dots A_n.
$
Among the numerous ways for studying such products, the problem of the (entrywise) convergence of the sequence $P_1,P_2,\dots$ itself (notion of Right Convergent Product) is solved by Daubechies and Lagarias \cite{DL92,DL01}. The probabilistic approach is exposed in the book of Bougerol and Lacroix \cite{BL85} with a large range of results about the products of random matrices. Given a probability measure $\mu$ on the set of complex-valued $d\times d$ matrices, \cite[Part A III Theorem 4.3]{BL85} gives two sufficient conditions for the convergence in probability of the normalized columns of~$P_n$ to a random vector, as well as for the almost sure convergence to $0$ of the angle between any couple of rows of~$P_n$. The first condition (strong irreducibility in \cite[Part A III Definition 2.1]{BL85}) means the non existence of a reunion of proper subspaces of $\mathbb C^d$ being stable by left multiplication by each element in the support of $\mu$. The second condition (contraction condition in
\cite[Part A III Definition 1.3]{BL85}) is equivalent to the existence for any $n$, of a matrix~product of matrices in the support of $\mu$ and whose normalization  converges to a rank $1$ matrix as $n$ tends to $+\infty$. This is completed in \cite[Part A VI Theorem 3.1]{BL85}: with a supplementary hypothesis, the normalized columns of $P_n$ converge almost surely to a random vector. (\cite{BL85} is a synthetic version of many results contained for instance in \cite{FK60}\cite{GR85}\cite{LeP84}\cite{Rau77}.)

Throughout the paper, we shall be mostly concerned by matrices  with  nonnegative entries and we shall focus in a non probabilistic projective convergence within the cone of the nonzero vectors with non negative entries. We start with a first remark: if one avoids the case when each one of the matrices $A$ such that $A_i=A$ infinitely many times, have a common left eigenvector (see Section~\ref{limpoints}), then the normalized matrix ${P_n}/\Vert P_n\Vert$ diverges as $n$ tends to $+\infty$. This is why we shall focus our attention with the asymptotic behavior of the (column) probability vectors ${P_nX}/{\Vert P_nX\Vert}$ for $X$ ranging over the whole  simplex $\mathcal{S}_d$ of the  probability vectors in $\RRR^d$. Here and throughout the paper, $\Vert\cdot\Vert$ stands for the norm applying to the column vectors (or more generally to matrices) and whose value is the sum of (the modulus of) the entries: hence, $\mathcal{S}_d$ is made of the nonnegative column vectors $X$ such that $\Vert X\Vert=1$. Actually the main result established in Theorem~A, gives a sufficient condition called ${\bf (C)}$, which allows to reduce (in some sense) the question of the convergence of a general probability vector ${P_nX}/{\Vert P_nX\Vert}$, for $X\in\mathcal{S}_d$, to the case of the {\it normalized columns of $P_n$}, i.e.  ${P_nU_j}/{\Vert P_nU_j\Vert}$ ($1\le j\le d$): here $U_1,\dots,U_d$ are the $d$ canonical $d\times 1$ vectors (and the extremal points of $\mathcal{S}_d$). More precisely, condition ${\bf (C)}$ in  Theorem~A, ensures the existence of nonempty disjoint sets $J_{h}(n)\subset\{1,\dots,d\}$, for each $n\ge 1$ and $h$ between $1$ and a constant $H\in\{1,\dots,d\}$, such that: (1)~:~each sequence of normalized columns of $P_n$, say $n\mapsto P_nU_{j_n}/\Vert P_nU_{j_n}\Vert$ with $j_n\in J_h(n)$, tends to a common  limit probability vector $V_h$; moreover (for $n$ large enough)  $P_nU_{j_n}$ and $V_h$ have the same row indices corresponding to  nonzero entries; 
(2)~:~the ratio $\Vert P_nU_{j_n'}\Vert/\Vert P_nU_{j_n}\Vert$ tends to $0$ as $n$ tends to $+\infty$ whenever $(j_n,j_n')\in J_h(n)\times J_{h'}(n)$ for $h<h'$;
(3)~:~$P_nU_j=0$ if and only if $j\not\in J_{1}(n)\cup\cdots\cup J_{H}(n)$. Finally, under condition ${\bf (C)}$, the limit of ${P_nX}/{\Vert P_nX\Vert}$ as $n$ tends to $+\infty$ exists for any probability vector $X$ with positive entries: this is  a consequence of a more accurate result (part (iii) of Theorem~A) concerned with any arbitrary probability vector $X$ s.t. $P_nX\ne0$, for  any $n\ge 1$.

 Let's mention that for $A_n=A(\xi_n)$, where $A(\0),\dots,A(\aaa)$ are nonnegative $2\times2$ matrices  the question is solved: one finds in \cite{OT13b} (resp. \cite{OT13a}) a necessary and sufficient condition for the pointwise  convergence (resp. uniform convergence w.r.t. the sequence $\xi_1\xi_2\cdots$) of ${P_nX}/{\Vert P_nX\Vert}$, assuming that $X$ has positive entries.  However the technics developed in these papers use specificities of the $2\times 2$ matrices and seems to us difficult to generalize. Our motivation comes from several works  concerned with  the geometry of fractal/multifractal sets and measures  (see \cite{McM84}\cite{Bed84}\cite{KP96a,KP96b}\cite{Yay09}\cite{Oli10a,Oli10b}\cite{Fen11}), with a special importance for the family  of the {\it Bernoulli convolutions} $\mu_\beta$ ($1<\beta<2$): the mass distributions $\mu_\beta$, whose support is a subinterval of the real line, arise in an {\it Erd\H os problem} \cite{Erd39} and are related to measure theoretic aspects of Gibbs structures for numeration with redundant digits
(see \cite{SV98}\cite{DST99}\cite{FO03}\cite{Fen05}\-\cite{OST05}\cite{BO07,BO10}\cite{Oli12}\cite{FLT13}). 
For some parameters $\beta$ (actually the so-called Pisot-Vijayaraghavan (PV) numbers) such a measure is {\it linearly representable}: the $\mu_\beta$-measure of a suitable family of nested generating intervals may be computed by means of matrix products of the form $P_nX$, where $A_n$ takes only finitely many values, say $A(\0),\dots,A(\aaa)$, and $X$ is a probability vector with positive entries. Because, $A_n=A(\xi_n)$, where $\xi=\xi_1\xi_2\cdots$ is a sequence with  $\xi_n\in\{\0,\dots,\aaa\}$,  we shall write $P_n=P_n(\xi)$ the dependence of the $n$-step product with~$\xi$: when the convergence of ${P_n(\xi)X}/{\Vert P_n(\xi)X\Vert}$ is uniform w.r.t. $\xi$, a sharp  analysis of the Bernoulli convolution $\mu_\beta$ (Gibbs structures and multifractal decomposition) becomes possible. 

The paper is organized as follows. Section~\ref{background} is devoted to the presentation of the  background necessary to establish Theorem~A, the proof of the theorem itself being detailed in Section~\ref{ProofThA}. We illustrate in Section~\ref{heuristic} many aspects of Theorem~A through several elementary examples. In Section~\ref{GibbsSection}, we show to what extent Theorem~A may be used to analyse the Gibbs properties of the linearly representable measures : we give two examples, the first one (Section~\ref{KAMAESEC}) being the {\it Kamae measure}, related with a self-affine graph studied by Kamae \cite{Kam86}, the second one (Section~\ref{ICBM}) being concerned with the Bernoulli convolutions $\mu_\beta$ already mentioned. In the latter case, most of the matrices involved in the linear decomposition of $\mu_\beta$ are large, sparse and it is usually not easy to verify condition ${\bf (C)}$ of Theorem~A. To illustrate the technics involved, Section~\ref{ICBM} is devoted to the PV-number $\beta\approx1.755$   s.t. $\beta^3=2\beta^2-\beta+1$, leading to 
\begin{equation}\label{matrix7X7}
\begin{matrix}
A(\0):=\begin{pmatrix}1&0&0&0&0&0&0\\0&0&1&0&0&0&0\\0&0&0&1&1&0&0\\0&0&0&0&0&0&0\\1&0&0&0&0&0&1\\0&0&0&0&1&0&0\\0&1&0&0&0&0&0\end{pmatrix}&
A(\1):=\begin{pmatrix}0&0&1&1&0&0&0\\0&0&0&0&0&1&0\\0&0&0&1&1&0&0\\1&0&0&0&0&0&0\\0&0&1&0&0&0&0\\0&0&0&0&0&0&0\\0&0&0&0&0&0&0\end{pmatrix}\\\\

A(\2):=\begin{pmatrix}1&0&0&0&1&0&1\\0&0&0&0&0&0&0\\1&0&0&0&0&0&1\\0&0&0&1&1&0&0\\0&0&0&0&1&0&0\\0&0&0&0&0&0&0\\0&0&0&0&0&0&0\end{pmatrix}&
\end{matrix}
\end{equation}
In Section~\ref{7X7}, we verify (in details) that condition ${\bf (C)}$ holds for (most of) the products $P_n(\xi)=A(\xi_1)\cdots A(\xi_n)$, where $A(\0),A(\1)$ and $A(\2)$ are the matrices in (\ref{matrix7X7}) related to $\mu_\beta$, so that we are able to establish the Gibbs properties of this measure.   
\bigskip

A complete bibliography about infinite products of matrices can be found in the paper of Victor Kozyakin \cite{Koz13}.

{\bf Acknowledgement. --} {\it Both authors are grateful to Ludwig Elsner and collaborators for their comments on a preliminary version of the present work; in particular this has incitated us to clarify the relations between Theorem~A and the rank $1$ asymptotic approximation of the normalized matrix products under condition ${\bf (C)}$ (See Section~\ref{limpoints})}.

\subsection{Statement of Theorem~A}
Given a $r\times 1$ vector $V$, we note
\begin{equation}\label{jfdhgdk}
\mathcal{I}(V):=\Big\{1\le i\le r\;;\;V(i)\ne0\Big\}
\quad\text{and}\quad
\Delta(V):=\sum\nolimits_{i\in\mathcal{I}(V)}U_i\;;
\end{equation}
the inclusion  
$\mathcal{I}(V)\subset\mathcal{I}(V')$ is thus equivalent to the inequality $\Delta(V)\le\Delta(V')$.
Major ingredients for Theorem~A are the sets ${\mathcal H}_1$, ${\mathcal H}_2(\Lambda)$ and ${\mathcal H}_3(\lambda)$ (for $\lambda\ge0$ and $\Lambda\ge1$), whose elements are $d_1\times d_2$ matrices and respectively defined by setting
\vskip0pt 
\begin{eqnarray}\label{H1 of the eqnarray}
A\in {\mathcal H}_1&\iff&\Big[\forall j_0,j_1,\Delta(AU_{j_0})\ge\Delta(AU_{j_1})\quad\text{or}\quad\Delta(AU_{j_0})\le \Delta(AU_{j_1})\Big]\;;\\\nonumber\\
\label{H2 of the eqnarray}A\in{\mathcal H}_2(\Lambda)&\iff&\Big[A(i_0,j_0)\ne 0\;\Longrightarrow\;\Vert AU_{j_0}\Vert\le \Lambda A(i_0,j_0)\Big]\;;\\\nonumber\\
\label{H3 of the eqnarray}A\in {\mathcal H}_3(\lambda)&\iff&\Big[A(i_0,j_0)\ne0,A(i_0,j_1)=0\;\Longrightarrow\;\Vert AU_{j_1}\Vert\le\lambda A(i_0,j_0)\Big].
\end{eqnarray}
\vskip15pt
\noindent For any nonnegative matrix $A\ne 0$, the following two constants are well defined:
\begin{equation}\label{lambda}
\Lambda_A:=\min\Big\{\Lambda\ge1\;;\;A\in\mathcal{H}_2(\Lambda)\Big\}
\quad\text{and}\quad
\lambda_A:=\min\Big\{\lambda\ge0\;;\;A\in\mathcal{H}_3(\lambda)\Big\}.
\end{equation}

\begin{definition}\label{C}Let $\mathcal{A}=(A_1,A_2,\dots)$ be a sequence of $d\times d$ matrices and $0=s_0=s_1<s_2<\dots$ a sequence of integers; given $n\ge0$ there exists $k=\kk(n)\ge0$ s.t. $s_{k+1}\le n<s_{k+2}$ and we note 
$$
P_{n}:=A_1\cdots A_{n}\quad\text{and}\quad Q_{n}:=A_{s_{k}+1}\cdots A_n,
$$
where by  convention $P_0=Q_0$ is the $d\times d$ identity matrix; 
hence, for any $n\ge 0$, 
\begin{equation}\label{decomposition}
P_n=P_{s_k}Q_n\quad
\quad\text{and}\quad
k\ge1\Longrightarrow P_{s_k}=Q_{s_1}\cdots Q_{s_k}.
\end{equation} 
The sequence $\mathcal{A}$ satisfies condition ${\bf (C)}$ w.r.t. $0\le\lambda<1\le\Lambda<+\infty$ and  
$(s_0,s_1,\dots)$ if
\begin{equation}\label{hypot}
n\ge s_2\Longrightarrow Q_n\in\mathcal{H}_1\cap\mathcal{H}_2(\Lambda)\cap \mathcal{H}_3(\lambda).
\end{equation} 
\end{definition}

\begin{figure}[H]
\begin{center}
\includegraphics[scale=0.45]{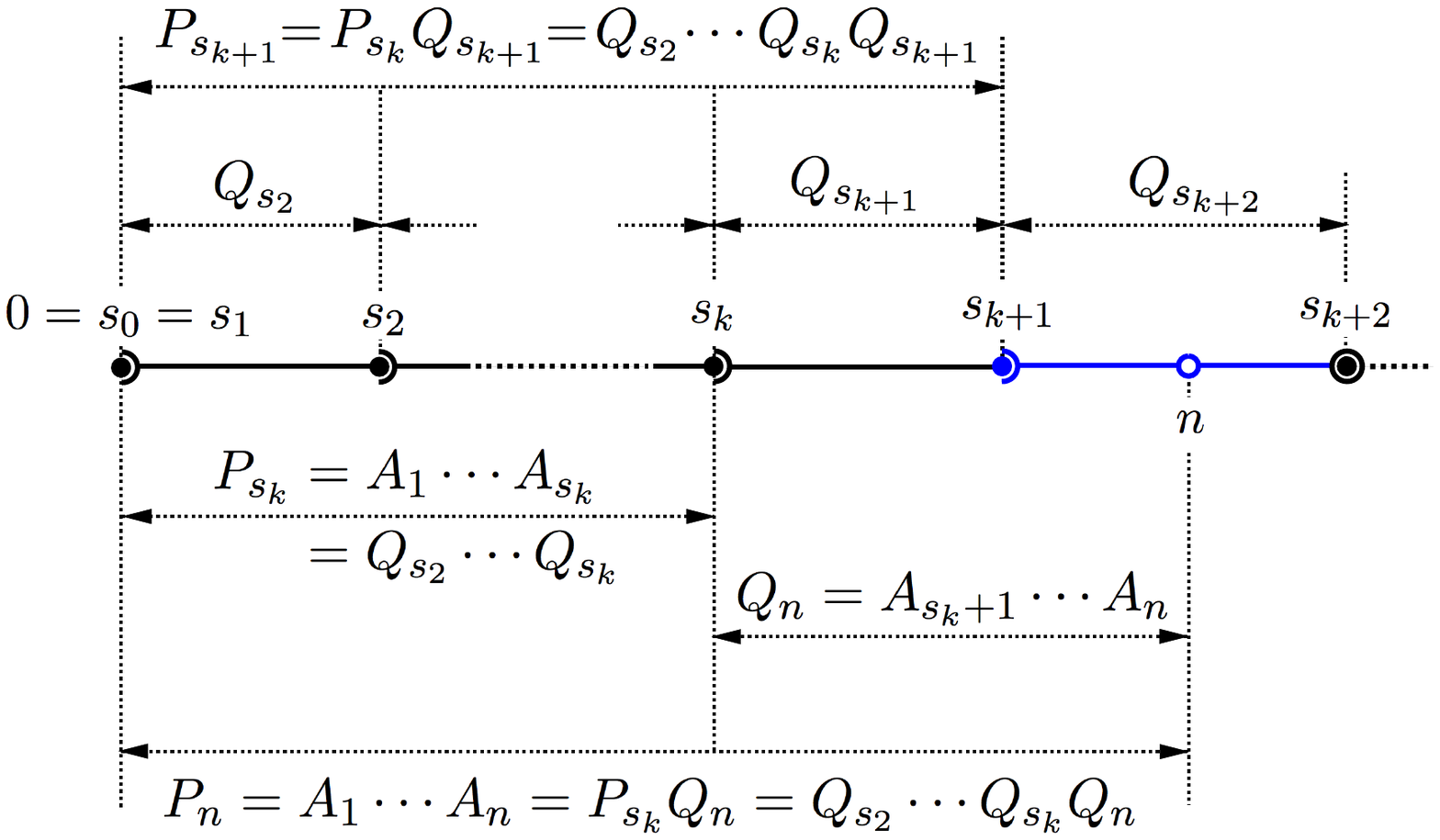}
\caption{\label{scalesk}\footnotesize\it The sequence $0=s_0=s_1<s_2<\cdots$ determines the definition of the $d\times d$ matrices $Q_n$ s.t. $P_{s_k}Q_n=P_n$ for any $n\ge 0$ and $s_{k+1}\le n<s_{k+2}$ (with $k=\kk(n)\ge0$); notice that $P_{n}=Q_{n}$ for $0\le n<s_3$ (i.e. $\kk(n)=0$ or $1$), with   
$P_0=Q_0$ being the $d\times d$ identity matrix.} 
\end{center}
\end{figure}




\begin{TheoremA}\label{rankone}
Let $\mathcal{A}=(A_1,A_2,\dots)$ be a sequence of $d\times d$ matrices satisfying condition ${\bf (C)}$; then, there exist $H$  probability vectors $V_1,\dots,V_H$ ($1\le H\le d$) with 
\begin{equation}\label{hhhhhhhhh}
\Delta(V_1)\ge \cdots \ge\Delta(V_H)\quad\text{while}\quad V_{h-1}\ne V_h\quad\text{(for any $1< h\le H$)}
\end{equation}
and there exists $\emptyset\ne J_h(n)\subset\{1,\dots,d\}$
($1\le h\le H$ and $n\in\mathbb N$ large enough), giving the partition
\begin{equation}\label{partition}
\Big\{(i,j)\;;\;P_n(i,j)\ne0\Big\}=\bigsqcup\nolimits_{h=1}^H I_h\times  J_{h}(n),
\quad\text{(where $I_h:=\mathcal{I}(V_h)$)}
\end{equation}
for which the following assertions hold:

(i) : for $1\le h\le H$ and $j_1,j_2,\dots$ a sequence in $\{1,\dots,d\}$ s.t. $j_n\in  J_h(n)$, then
$$
\lim_{n\to+\infty}{P_nU_{j_n}/\Vert P_nU_{j_n}\Vert}=V_h\;;
$$

(ii) : for  $1< h\le H$ and $j_1,j_2,\dots$ and $j_1',j_2',\dots$ two sequences in $\{1,\dots,d\}$,
$$
(j_n,j_n')\in  J_{h-1}(n)\times  J_{h}(n)\;\Longrightarrow\;\lim_{n\to+\infty}{\Vert P_nU_{j_n'}\Vert/\Vert P_nU_{j_n}\Vert}=0\;;
$$

(iii) : there exists  real numbers $\varepsilon_1,\varepsilon_2,\cdots$  with $\varepsilon_n\to0$ as $n\to+\infty$ such that, for any $X\in\mathcal{S}_d$ for which each  $P_nX\ne0$, one has for any $n\ge1$: 
\begin{equation}\label{Vhn}
\left\Vert {P_nX\over \Vert P_nX\Vert}-V_{\hhh_X(n)}\right\Vert\le \Lambda_X\cdot \varepsilon_n
\quad\text{where}\quad \hhh_X(n)=\min\Big\{h\;;\; \mathcal{I}(X)\cap J_{h}(n) \ne\emptyset\Big\}.
\end{equation}
\end{TheoremA}

We define the {\it upper/lower top Lyapunov exponent} of $\mathcal{A}=(A_1,A_2,\dots)$ to be the quantities 
\begin{equation}\label{gkhjfdskjfh}
\underline{\chi}_{\rm top}(\mathcal{A}):=\liminf_{n\to+\infty}{1\over n}\log\Vert A_1\cdots A_n\Vert
\quad\text{and}\quad
\overline{\chi}_{\rm top}(\mathcal{A}):=\limsup_{n\to+\infty}{1\over n}\log\Vert A_1\cdots A_n\Vert
\end{equation}
and we speak of the {\it top Lyapunov exponent} denoted $\chi_{\rm top}(\mathcal{A})$ whenever both lower and upper exponents $\underline{\chi}_{\rm top}(\mathcal{A})$ and $\overline{\chi}_{\rm top}(\mathcal{A})$ coincide. A straightforward consequence of part (iii) of Theorem~A, is the following corollary.

\begin{CorollaryA}\label{CorollaryA}Let  $\mathcal{A}=(A_1,A_2,\dots)$ be a sequence of $d\times d$ matrices satisfying condition ${\bf (C)}$; then,  here exists a unique $V\in\mathcal{S}_d$ (the top Lyapunov direction) such that for any vector $X$ with positive entries, 
$$
\lim_{n\to+\infty}{P_nX/\left\Vert P_nX\right\Vert}=V\;.
$$ 
\end{CorollaryA}

We notice that, if $e^{n\chi_1(n)}\ge\dots\ge e^{n\chi_d(n)}$ is the ordered list of the singular values of $P_n$ ($n \ge  1$), $e^{n\chi_1(n)}$ is also the euclidean norm of $P_n$, hence from (\ref{gkhjfdskjfh})
$$
\liminf_{n\to+\infty}\chi_1(n)=\underline{\chi}_{top}(\mathcal A)\quad\hbox{and}\quad\limsup_{n\to+\infty}\chi_1(n)=\overline{\chi}_{top}(\mathcal A).
$$

The exponents of the form $\chi_{\rm top}(\mathcal{A})$ are found in a probabilistic framework which roots in the seminal work by Furstenberg \& Kesten~\cite{FK60}. To fix ideas, let $T:\Omega\to\Omega$ be a continuous transformation, where (for simplicity) $\Omega$ is a compact metric space and suppose that $\mu$ is a $T$-ergodic borelian probability measure. 
We also consider that $A:\Omega\to\mathcal{M}_d(\CCC)$ is a given (borelian) map from $\Omega$ to the space $\mathcal{M}_d(\CCC)$ of the complex $d\times d$ matrices. We note $P_n(\omega)=A(\omega)A(T(\omega))\cdots A(T^{n-1}(\omega))$, with the convention that $P_0(\omega)$ is the $d\times d$ identity matrix;  
under these conditions 
$(n,\omega)\mapsto\Vert P_n(\omega)\Vert$ is a submultiplicative process since 
$\Vert P_{n+m}(\omega)\Vert\le \Vert P_n(\omega)\Vert\cdot \Vert P_m(T^n(\omega))\Vert$
and  Kingman's Subadditive Ergodic Theorem ensures the existence of a constant $\chi_1(\mu)\ge-\infty$ s.t. 
\begin{equation}\label{toplyap}
\lim_{n\to+\infty}{1\over n} \log\Vert P_n(\omega)\Vert=\inf\left\{{1\over n} \log\Vert P_n(\omega)\Vert\;;\;n\ge1\right\}=\chi_1(\mu)\quad\text{$\mu$-a.s.}
\end{equation}
The quantity $\chi_1(\mu)$ is the top-Lyapunov exponent of the random process $(n,\omega)\mapsto\Vert P_n(\omega)\Vert$, for $\Omega$ weighted by $\mu$. In view of the definition in (\ref{gkhjfdskjfh}), notice that for $\omega$ being a {\it $\mu$-generic point} in $\Omega$ (hence for $\mu$-a.e. $\omega$)
one has 
$\chi_{\rm top}\big(\mathcal{A}(\omega)\big)=\chi_1(\mu)$, where $\mathcal A(\omega) := (A(\omega), A(T(\omega)), A(T^2(\omega)), \dots )$.  
A more general framework for characteristic exponents associated with matrix products is given by Oseledets Theorem \cite{Ose68}:~indeed,~if 
$$
e^{n\chi_1(n,\omega)}\ge\cdots \ge e^{n\chi_d(n,\omega)}
$$
form the ordered sequence of the singular values of $P_n(\omega)$,
then (for $\mu$ being $T$-ergodic) each exponent $\chi_k(n,\omega)$ tends $\mu$-a.s. toward a limit $\chi_k(\mu)$ which is the $k$-th Lyapunov exponent of $\mu$. (The larger Lyapunov exponent $\chi_1(\mu)$ coincides with the {\it top Lyapunov exponent} as defined in (\ref{toplyap}), hence his name.) The theory of Lyapunov exponents related to matrix products is a wide domain of research: we mention relationships with Hausdorff dimension of stationary probabilities and multifractal analysis of positive measures (see for instance \cite{Led84}\cite{FL02}\cite{Fen03,Fen04,Fen09}\cite{FH10}\cite{DHX11}).

\section{\bf Notations and background}\label{background}

\subsection{\bf Projective distance and contraction coefficient} 
The definitions and the properties of the Hilbert metric $\delta_\mathcal{H}(\cdot,\cdot)$ and of the Birkhoff (or ergodic) contraction coefficient $\tau_\mathcal{B}(\cdot)$ may be found in the very beginning  of Subsection 3.4 of Seneta's book \cite{Sen81};
we shall need (in particular in the proof of Theorem A) an adaptation/generalization of some concepts. Indeed, the Birkhoff contraction coefficient $\tau_\mathcal{B}(A)$ of a square matrix $A$ (with nonnegative entries) belongs to the unit interval $[0\,;1]$ with the crucial property that $\tau_\mathcal{B}(A)<1$ if and only if $A$ has positive entries: however, the usual framework of Theorem~A is concerned with sparse matrices. Hence, we shall consider a contraction coefficient map $A\mapsto\tau(A)$ whose domain is made of the (non necessarily square) matrices having nonnegative entries, and such that $\tau(A)<1$ if and only if the  positive entries are positioned on a rectangular submatrix.
\begin{definition}(i) : The $d_1\times d_2$ nonnegative matrix $A=(A(i,j))$, distinct from the null matrix, is said to satisfy hypothesis ${\bf (H)}$ if there exist two nonempty sets $I_A\subset\{1,\dots,d_1\}$ and $J_A\subset\{1,\dots,d_2\}$ such that
$$
A(i,j)\ne0\iff (i,j)\in I_A\times J_A\;;
$$
(ii) : the $\delta$-coefficient of $A$ is either $\delta(A):=+\infty$ if $A$ does not satisfy ${\bf (H)}$, or otherwise, 
\begin{equation}\label{Defdelta}
\delta(A):=\max\left\{\log\left({A(i,j)\cdot A(k,\ell)\over A(k,j)\cdot A(i,\ell)}\right)\;;\;(i,k)\in I_A\times I_A,\;(j,\ell)\in J_A\times J_A\right\}\;;
\end{equation}
(iii) : the generalized contraction coefficient of $A$ is:
\begin{equation}\label{Deftau}
\tau(A):=\tanh\left(\frac{\delta(A)}4\right).
\end{equation}
\end{definition}
\begin{proposition}
$\tau(A)<1$ if $A$ satisfies ${\bf (H)}$ and $\tau(A)=1$ otherwise.
\end{proposition}
If $A=(X\ Y)$ is a $d\times 2$ nonnegative matrix (i.e. $X=AU_1$ and $Y=AU_2$), then  we shall make an abuse of notations writing $\delta(X,Y)$ instead of $\delta(A)$: we call $\delta(X,Y)$ the projective distance between $X$ and $Y$
(likewise we shall note $\delta(X^\star,Y^\star):=\delta(A^\star)=\delta(A)$ for row vectors). So the restriction of $\delta$ to the simplex $S_d$ is an extended metric, and defines a topology on $S_d$.
Let $\emptyset\ne I\subset\{1,\dots,d\}$; if one assumes that $\Delta(X)=\Delta(Y)$ with $\mathcal{I}(X)=\mathcal{I}(Y)=I$ then, the projective distance $\delta(X,Y)$ ($=\delta(X^\star,Y^\star)$) is finite and given by the two following equivalent expressions:
\begin{equation}\label{equivHILBERT}
\delta(X,Y)=
\max_{i,j\in I}\left\{\log\left(\frac{X(i)Y(j)}{Y(i)X(j)}\right)\right\}
=\max_{i\in I}\left\{\log\left(\frac{X(i)}{Y(i)}\right)\right\}+\max_{i\in I}\left\{\log\left(\frac{Y(i)}{X(i)}\right)\right\}.
\end{equation}
Moreover, $\delta(X,Y)$ coincides with $\delta\big(X/\Vert X\Vert,Y/\Vert Y\Vert\big)$: hence $\delta(\cdot,\cdot)$ is entirely determined by its values on the simplex $\mathcal{S}_d$. The map $\delta(\cdot,\cdot):\mathcal{S}_d\times\mathcal{S}_d\to[0\,;+\infty]$ is closed to a metric but recall that $\delta(X, Y )=+\infty$ means $\Delta(X)\ne\Delta(Y)$. 
\begin{definition}\label{DefSI}  For $\emptyset\ne I\subset\{1,\dots,d\}$, we denote by $\mathcal{S}_I$ the set of vectors $X\in\mathcal{S}_d$ s.t. 
$\mathcal{I}(X)=I$.
\end{definition}


\begin{proposition}\label{face}Let  $\emptyset\ne I\subset\{1,\dots,d\}$: then, (i) : for  $X,Y\in\mathcal{S}_d$, one has $\delta(X,Y)<+\infty\iff\Delta(X)=\Delta(Y)$; (ii)~:~the restriction $\delta(\cdot,\cdot):\mathcal{S}_I\times\mathcal{S}_I\to[0\,;+\infty[$ defines a metric on $\mathcal{S}_I$ (when $I=\{1,\dots,d\}$, the metric $\delta(\cdot,\cdot)$ over $\mathcal{S}_I$ (or $\mathcal{S}_I^\star$) coincides with Hilbert projective metric $\delta_\mathcal{H}(\cdot,\cdot)$ already mentioned);  
(iii)~:~if $X,Y\in\mathcal{S}_I$,  then
\begin{equation}\label{doublebound}
{1\over d}\cdot\Vert X-Y\Vert\le \delta(X,Y)\le{\Vert X-Y\Vert\over \min_{i\in I}\{X(i),Y(i)\}}.
\end{equation}
\end{proposition}

\begin{proof}[{\bf Proof}]Part (i) is straightforward from the definition of $\delta(\cdot,\cdot)$ while part (ii) is a key property of Hilbert projective metric (see \cite{Sen81}\cite{Har02}). To prove the double inequality in (\ref{doublebound}) of part (iii), we shall use the following inequalities,
\begin{equation}\label{kqds}
a-b\le \log\left({a\over b}\right)\le {a-b\over b}
\end{equation}
valid as soon as  $1\ge a\ge b>0$. Fix  $\emptyset\ne I\subset\{1,\dots,d\}$ with $\#I\ge2$ (otherwise the result is trivial): we then consider $X,Y\in\mathcal{S}_I$. 
First, up to a permutation of $X$ and $Y$, we assume  $\Vert X-Y\Vert_\infty=X(i_0)-Y(i_0)$ for $i_0\in I$, so that $1\ge X(i_0)\ge Y(i_0)>0$; however, since $X,Y\in\mathcal{S}_d$, it is necessary that $\sum_i(X(i)-Y(i))=0$ and thus, there exits $i_1$ s.t. $X(i_1)-Y(i_1)\le 0$; in particular $\max\nolimits_{i\in I}\big\{\log (Y(i)/X(i))\big\}\ge0$: then, we use the lower bound in (\ref{kqds}) together with the second expression of  $\delta(X,Y)$ in (\ref{equivHILBERT}) to write
\begin{align*}
\delta(X,Y)&\ge \log\big(X(i_0)/Y(i_0)\big)
= X(i_0)-Y(i_0)\ge\Vert X-Y\Vert_\infty\ge{1\over d}\cdot\Vert X-Y\Vert,
\end{align*}
proving the lower bound in (\ref{doublebound}). For the corresponding upper bound, let $i_0'$ (resp. $i_1'$) s.t. $X(i'_0)/Y(i'_0) = \max_i(X(i)/Y(i))$ (resp. $Y(i'_1)/X(i'_1) = \max_i(Y(i)/X(i))$ ). If $i'_0 = i'_1$ then $X =Y$, so that $\delta(X,Y)=\Vert X-Y\Vert=0$ and the desired upper bound holds. Now suppose $i_0'\ne i_1'$: because $\sum_i(X(i)-Y(i))=0$, we know that  $1\ge X(i_0')\ge Y(i_0')\ge\varepsilon$ and $1\ge Y(i_1')\ge X(i_1')\ge\varepsilon$, where $\varepsilon:=\min_{i\in I}\{X(i),Y(i)\}$: using the upper bound in (\ref{kqds}) together with the second expression of  $\delta(X,Y)$ in (\ref{equivHILBERT}) gives
\begin{align*}
\delta(X,Y)&= \log\big(X(i_0')/Y(i_0')\big)+\log\big(Y(i_1')/X(i_1')\big)\\
&\le {1\over \varepsilon}\cdot\big(X(i_0')-Y(i_0')\big)+{1\over \varepsilon}\cdot\big(Y(i_1')-X(i_1')\big)\le{1\over \varepsilon}\cdot \Vert X-Y\Vert.
\end{align*}

\end{proof}

\begin{remark}
Let $X_*,X_1,X_2,\dots\in\mathcal{S}_d$;  if $\delta(X_n,X_*)\to0$ as $n\to+\infty$ then $\Delta(X_n)=\Delta(X_*)$ as soon as $\delta(X_n,X_*)<+\infty$. So one has the equivalence 
\begin{equation}\label{EI}
\lim_{n\to+\infty}\delta(X_n,X_*)=0\iff
\Big(\lim_{n\to+\infty}\left\Vert X_n-X_*\right\Vert=0\hbox{ and }\Delta(X_n) = \Delta(X_*)\hbox{ for }n\hbox{ large enough}\Big).
\end{equation}
We make the simple remark that if each one of the nonzero entries $X_n(i)$ are bounded from bellow by $\varepsilon>0$, then the normed convergence $X_n\to X_*$ is equivalent to $\delta(X_n,X_*)\to0$, with more precisely a convergence $X_n\to X_*$ w.r.t. the metric $\delta(\cdot,\cdot)$ over $\mathcal{S}_{\mathcal{I}(X_*)}$. This is the key point of Lemma~\ref{trapping}  in the proof of Theorem~A in Section~\ref{ProofThA}.
\end{remark}

Given $X$ and $Y$ two arbitrary $d\times 1$ vectors with positive entries and $B$ a square $d\times d$ allowable (without null row nor null column) nonnegative matrix, 
one recovers that (see \cite[\S~3.4]{Sen81})
$$
\delta(X^\star,Y^\star)=\delta_{\mathcal H}(X^\star,Y^\star)\quad\text{and}\quad\tau(B)=\tau_{\mathcal B}(B):=\sup\left\{{\delta_{\mathcal H}(X^\star B,Y^\star B)\over \delta_{\mathcal H}(X^\star,Y^\star)}\;;\;X,Y\hbox{ non colinear}\right\},
$$
which gives the  contraction property of $\delta_\mathcal{H}(\cdot,\cdot)$ and  $\tau_\mathcal{B}(\cdot)$ that is 
\begin{equation}\label{Seneta-equation}
\delta_{\mathcal H}(X^\star B,Y^\star B)\le\delta_{\mathcal H}(X^\star,Y^\star)\tau_{\mathcal B}(B).
\end{equation}

\begin{lemma}\label{deltatau}
For any $d_1\times d_2$ nonnegative matrix $A$ satisfying ${\bf(H)}$ and any $d_2\times d_3$ nonnegative matrix $B$, if $AB\ne0$ one has:
\begin{equation}\label{variante}
\delta(AB)\le\delta(A)\tau(B).
\end{equation}

\end{lemma}

\begin{proof}[{\bf Proof}]To begin with, consider that $A$ is a $d_1\times d_2$ matrix with positive entries together with $B$ a nonnegative allowable  $d_2\times d_3$ matrix with $d_3\le d_2$. We claim the contraction formula
\begin{equation}\label{Seneta-equation-Bis}
\delta(AB)\le \delta(A)\tau(B)
\end{equation}
to be valid in this case. The case $d_3=d_2$ is actually equivalent to (\ref{Seneta-equation}). On the other hand, (\ref{Seneta-equation-Bis}) is still valid if $d_3<d_2$. Indeed, consider the square allowable matrix $B'=(BU_1\cdots BU_{d_3}\cdots BU_{d_3})$: by the first case,  $\delta(AB')\le \delta(A)\tau(B')$ and (\ref{Seneta-equation-Bis}) holds, because $\delta(AB')=\delta(AB)$ and $\tau(B')=\tau(B)$.

Now, suppose that $A$ is a $d_1\times d_2$ nonnegative matrix satisfying ${\bf(H)}$ and $B$ a $d_2\times d_3$ nonnegative matrix. Clearly $AB$, if it differs from the null matrix, satisfies ${\bf(H)}$, and either $AB$ does not have any $2\times 2$ submatrix with positive entries (in this case $\delta(AB)=0\le\delta(A)\tau(B)$), or there exists $i\ne j$ both in $\{1,\dots,d_1\}$ and $k\ne\ell$ both in $\{1,\dots, d_3\}$ so that $\delta(AB)=\delta(A'B')$, where $A'B'$ has positive entries and
$$
A':=\begin{pmatrix}
A(i,1)&\cdots&A(i,d_2)\\
A(j,1)&\cdots&A(j,d_2)\\
\end{pmatrix}
\quad
B':=\begin{pmatrix}
B(k,1)&B(\ell,1)\\
\vdots&\vdots\\
B(k,d_2)&B(\ell ,d_2)
\end{pmatrix}.
$$
Since $A$ is supposed to satisfy ${\bf (H)}$, the set $J_A\subset\{1,\dots,d_2\}$ is non empty and  because $A'B'$ has positive entries, it is licit to consider the maximal (necessarily non empty) set $I\subset J_A$  for which the $I\times\{1,2\}$ submatrix $B''$ of $B'$ is allowable. Now, if $A''$ is the $\{1,2\}\times I$ submatrix of $A'$, one has 
$\delta(AB)=\delta(A'B')=\delta(A''B'')$: because $A''$ has positive entries and $B''$ is allowable, one deduces by (\ref{Seneta-equation-Bis}) that $\delta(A''B'')\le \delta(A'')\tau(B'')$ (the case where $B''$ has only one row is trivial: $A''B''$ has rank $1$ and $\delta(A''B'')=0$). However, for $A''$ (resp. $B''$) being a submatrix of $A$ (resp. $B$), one has $\delta(A'')\le \delta(A)$ (resp. $\tau(B'')\le \tau(B)$).

\end{proof}

\subsection{\bf Some properties of $\mathcal{H}_1,\mathcal{H}_2(\cdot )$ and $\mathcal{H}_3(\cdot)$}
To begin with, suppose that $A$ is a nonnegative $d_1\times d_2$ matrix  in ${\mathcal H}_2(\Lambda)\cap {\mathcal H}_3(\lambda)$ with 
$A(i_0,j_0)\ne0$ while  $A(i_0,j_1)=0$: if $A(i,j_0)\ne0$, then 
\begin{equation}\label{j0j1}
\Vert AU_{j_1}\Vert\le\lambda\cdot A(i_0,j_0)\le (\lambda\Lambda)\cdot A(i,j_0).
\end{equation}

\begin{lemma}\label{majoration}
Given $\lambda\ge0$, $\Lambda\ge1$ and $A\in\mathcal{H}_2(\Lambda)\cap\mathcal{H}_3(\lambda)$, the following proposition holds:  

if $\Delta(AU_{j_0})>\Delta(AU_{j_1})$ then
$\max\nolimits_i\big\{A(i,j_1)\big\}\le \Vert AU_{j_1}\Vert\le(\lambda\Lambda)\min_i\big\{A(i,j_0)\ne0\big\}$.

\end{lemma}
\vskip-20pt
\begin{figure}[H]
\begin{center}
\includegraphics[trim=0 170 0 90,clip,scale=0.5]{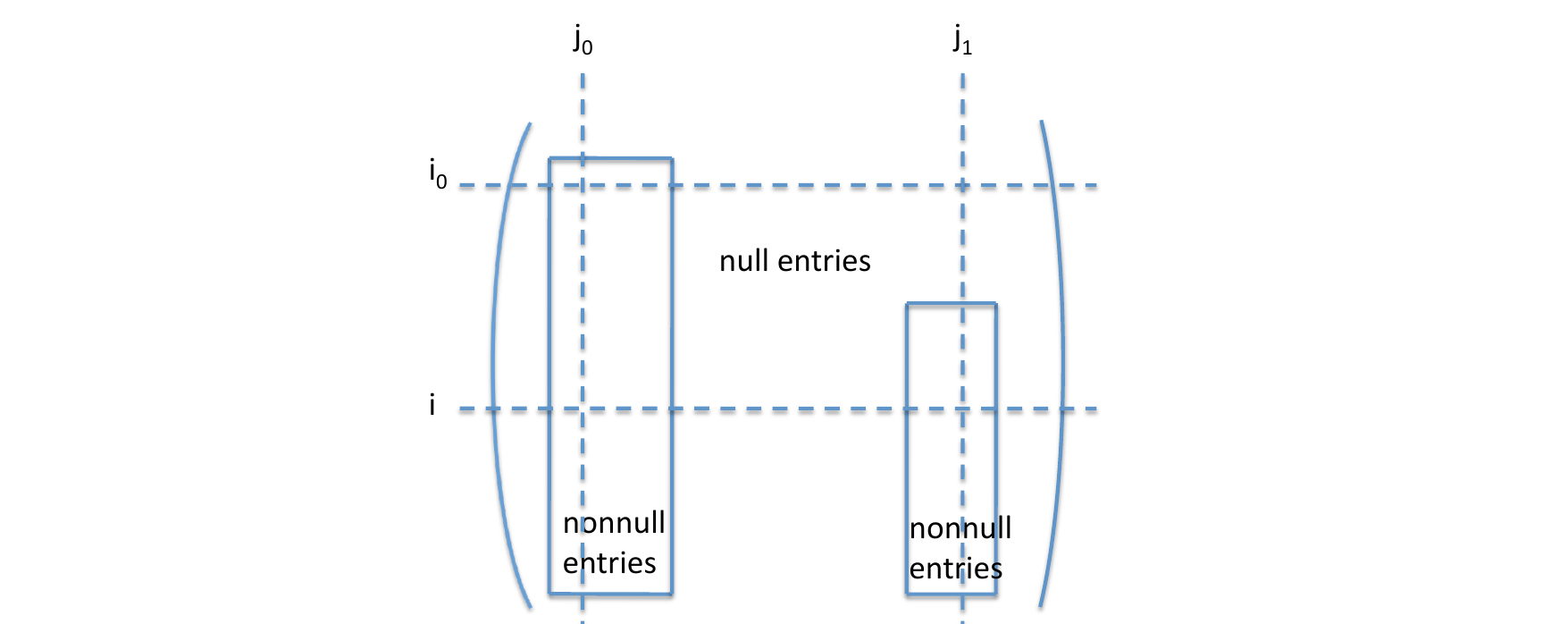}
\caption{\label{Un}\footnotesize\it Illustration of inequality (\ref{j0j1}).} 
\end{center}
\end{figure}

If $A\in\mathcal{H}_2(\Lambda)$ satisfies condition ${\bf (H)}$ (i.e. its nonzero entries are positioned on a rectangular submatrix), then the  $\delta$-coefficient of $A$ --~which is  finite~-- may actually be bounded by means of $\Lambda$. Indeed, if $A(i_0,j_0)A(i_1,j_1)A(i_1,j_0)A(i_0,j_1)\ne0$ then, the condition $A \in \mathcal H_2(\Lambda)$ implies 
$A(i_0,j_0)\le \Vert AU_{j_0}\Vert\le \Lambda\cdot A(i_1,j_0)$ and $A(i_1,j_1)\le \Vert AU_{j_1}\Vert\le \Lambda\cdot A(i_0,j_1)$, so that 
$$
\log\left({A(i_0,j_0)\over A(i_1,j_0)}\cdot {A(i_1,j_1)\over A(i_0,j_1)}\right)\le
\log\left({\Lambda\cdot A(i_1,j_0)\over A(i_1,j_0)}\cdot {\Lambda\cdot A(i_0,j_1)\over A(i_0,j_1)}\right)\le \log(\Lambda^2).
$$

\begin{lemma}\label{logLambda}
If a $d_1\times d_2$ matrix $0\ne A\in\mathcal{H}_2(\Lambda)$ satisfies condition ${\bf (H)}$ then, $\delta(A)\le\log(\Lambda^2)$.
\end{lemma}

We shall also need some stability properties of $\mathcal{H}_2(\cdot)$ and $\mathcal{H}_3(\cdot)$ w.r.t.  matrix multiplication.

\begin{lemma}\label{hypotheses}Let $A$ and $B$ be two matrices with non negative entries, of size $d_1\times d_2$ and $d_2\times d_3$ respectively; then,

(i) :  if $A\in\mathcal{H}_3(\lambda_a)$ and $B\in\mathcal{H}_3(\lambda_b)$ then $AB\in\mathcal{H}_3(\lambda_a\lambda_b)$;

(ii) :  if $A\in\mathcal{H}_2(\Lambda_a)\cap\mathcal{H}_3(\lambda_a)$ and $B\in\mathcal{H}_2(\Lambda_b)$ then $AB\in\mathcal{H}_2\Big(\Lambda_a+\lambda_a\Lambda_b\Big)$.

\end{lemma}

\begin{proof}[{\bf Proof}](i) : 
Suppose that  $AB(i_0,j_0)\ne0$ and $AB(i_0,j_1)=0$. Notice that $AB(i_0,j_0)\ne0$ means the existence of $j_*$ s.t. $A(i_0,j_*)B(j_*,j_0)\ne0$, while  $AB(i_0,j_1)=0$ implies $A(i_0,j)B(j,j_1)=0$ for any $j$; we also emphasize that $A(i_0,j_*)\ne0$ and $A(i_0,j_*)B(j_*,j_1)=0$ implies $B(j_*,j_1)=0$. To compare $\Vert ABU_{j_1}\Vert=\sum_iAB(i,j_1)$ with $AB(i_0,j_0)$ we write successively
\begin{eqnarray*}
\Vert ABU_{j_1}\Vert &=&	\sum_{A(i_0,j)\ne0}\Vert AU_j\Vert B(j,j_1)+\sum_{A(i_0,j)=0} \Vert AU_j\Vert B(j,j_1)\\
(B(j,j_1)=0\hbox{ when }A(i_0,j)\ne0)&\le& \lambda_aA(i_0,j_*)\sum_jB(j,j_1)\\
&=& \lambda_aA(i_0,j_*)\Vert BU_{j_1}\Vert\\
(B(j_*, j_0) \ne	0\hbox{ and }B(j_*, j_1) = 0)&\le& \lambda_aA(i_0,j_*)\lambda_bB(j_*,j_0)\\
&=& \lambda_a\lambda_bAB(i_0,j_0)\frac{A(i_0,j_*)B(j_*,j_0)}{AB(i_0,j_0)}
\end{eqnarray*}
and $\Vert ABU_{j_1}\Vert\le \lambda_a\lambda_b AB(i_0,j_0)$: this proves that $AB\in\mathcal{H}_3(\lambda_a\lambda_b)$.

(ii) : Suppose that $AB(i_0,j_0)\ne 0$ with $A(i_0,j_*)B(j_*,j_0)\ne0$; then, one has:
\begin{eqnarray*}
\Vert ABU_{j_0}\Vert
&=&\sum\nolimits_{A(i_0,j)\ne0}\Vert AU_j\Vert B(j,j_0)+
\sum\nolimits_{A(i_0,j)=0}\Vert AU_j\Vert B(j,j_0)\\
\text{($A(i_0,j_*)\ne0$)}&\le&\Lambda_a\sum\nolimits_{j}A(i_0,j)B(j,j_0)+
\lambda_a A(i_0,j_*)\sum\nolimits_{j}B(j,j_0)\\
&=&\Lambda_aAB(i_0,j_0)+\lambda_a A(i_0,j_*)\Vert BU_{j_0}\Vert\\
\text{($B(j_*,j_0)\ne0$)}&\le&\Lambda_aAB(i_0,j_0)+\lambda_a A(i_0,j_*)\Lambda_bB(j_*,j_0)\\
&=&AB(i_0,j_0)\left(\Lambda_a+\lambda_a\Lambda_b{A(i_0,j_*)B(j_*,j_0)\over AB(i_0,j_0)}\right),
\end{eqnarray*}
that is $\Vert ABU_{j_0}\Vert\le AB(i_0,j_0)\big(\Lambda_a+\lambda_a\Lambda_b\big)$: this proves that $AB\in\mathcal{H}_2\Big(\Lambda_a+\lambda_a\Lambda_b\Big)$.

\end{proof}


Given  $M$ a $r\times s$ matrix, we define 
\begin{equation}\label{jfdhg}
\#{\bf Col}(M):=\#\Big\{\Delta(MU_1),\dots,\Delta(MU_s)\Big\}\setminus\Big\{0\Big\}.
\end{equation}
\begin{lemma}\label{hypothesesbis}
Let $A,B$ and $C$  be three $d\times d$ matrices 
with non negative entries; then 

(a) : $\Delta(AU_{j_0})\le\Delta(AU_{j_1}) \Longrightarrow \Delta(BAU_{j_0})\le\Delta(BAU_{j_1})$;

(b) : $\Delta(AU_{j_0})=\Delta(AU_{j_1}) \Longrightarrow \Delta(BAU_{j_0})=\Delta(BAU_{j_1})$;

(c) : $\Delta(AU_{j_0})=0 \Longrightarrow \Delta(BAU_{j_0})=0$;

\noindent moreover if  $A\in{\mathcal H}_1$ then 

(d) : $BACU_j\ne 0\Longrightarrow \Delta(BACU_j)\in\big\{\Delta(BAU_1),\dots,\Delta(BAU_d)\big\}$; 

(e) : $\#{\bf Col}(BAC)\le\#{\bf Col}(A)$; 

(f) : $BAC\in{\mathcal H}_1$.
\end{lemma}

\begin{proof}[{\bf Proof}]To begin with, any column of $BA$ is a linear combination of the columns of $B$: more precisely, for any $j=1,\dots,d$, one has $AU_j=\sum_iA(i,j)U_i=\sum_{i\in\mathcal{I}(AU_j)}A(i,j)U_i$ and thus
\begin{equation}\label{fgj}
BAU_j=\sum\nolimits_{i\in\mathcal{I}(AU_j)}A(i,j)BU_i,
\end{equation}
which gives (a-b-c). If $A\in\mathcal{H}_1$, then (\ref{fgj}) with (a-b-c)  prove that 
$BA\in{\mathcal H}_1$ and $\#{\bf Col}(BA)\le\#{\bf Col}(A)$.
Like in (\ref{fgj}), each column of $BAC$ is a linear combination of the columns of $BA$, with 
\begin{equation}\label{fgjhdk}
BACU_j=\sum\nolimits_{i\in\mathcal{I}(CU_j)}C(i,j)BAU_i.
\end{equation}
Suppose that $BACU_j\ne0$ and let $i_1,\dots,i_s$ be the indices in $\mathcal{I}(CU_j)$ s.t. $\Delta(BAU_{i_k})\ne0$; because $BA\in\mathcal{H}_1$, one may suppose that 
$\Delta(BAU_{i_1})\ge\cdots\ge\Delta(BAU_{i_s})$
and by (\ref{fgjhdk}), one deduces $\Delta(BACU_j)=\Delta(BAU_{i_1})$: this is the content of part (d). Part (d) together with $BA\in\mathcal{H}_1$ shows that  $BAC\in{\mathcal H}_1$, proving (f). Part (d) also ensures $\#{\bf Col}(BAC)\le\#{\bf Col}(BA)$ and  (e) holds, for we already know that $\#{\bf Col}(BA)\le\#{\bf Col}(A)$.

\end{proof}

\section{\bf Proof of Theorem~A}\label{ProofThA}

\subsection{Preparatory lemmas}\label{convergence}
We shall consider that $\mathcal{A}=(A_1,A_2,\dots)$ is a given sequence of $d\times d$ matrices with nonnegative entries and  satisfying condition ${\bf (C)}$. A key point of the  argument leading to Theorem~A, is to {\it reenforce}  condition ${\bf (C)}$ (see condition in (\ref{hypotplus}) and (\ref{asterixbis}) below) without further assumptions about the sequence $\mathcal{A}$.
\begin{lemma}\label{increasing}
Suppose that  $\mathcal{A}=(A_1,A_2,\dots)$ satisfies condition ${\bf (C)}$ w.r.t. $0\le\lambda<1\le\Lambda<+\infty$ and the sequence of integers  $0=s_0 =s_1 <s_2<\dots$;  then for any $s_{k+1}\le n<s_{k+2}$ with $k\ge 1$,   
\begin{equation}\label{asterix}
1\le p\le k\Longrightarrow
Q_n\hbox{ and }Q_{s_p}\cdots Q_{s_{k}}Q_n\in\mathcal H_2(\Lambda')
\quad\text{where}\quad
\Lambda'={\Lambda\over 1-\lambda}\;;
\end{equation}
(by convention $Q_{s_1}=Q_0$ is the identity matrix); moreover, up to a change of the sequence $(s_0,s_1,\dots)$, it is licit to consider that
\begin{equation}\label{obelix}
n\ge s_2\Longrightarrow Q_n\hbox{ and }Q_{s_p}\cdots Q_{s_{k}}Q_n\in{\mathcal H}_1\cap {\mathcal H}_2(\Lambda')\cap{\mathcal H}_3(\lambda^k).
\end{equation}
\end{lemma}
\begin{proof}[{\bf Proof}]Let $\mathcal{A}=(A_1,A_2,\dots)$ satisfy ${\bf (C)}$ w.r.t. $0\le\lambda<1\le\Lambda$ and $(s_0,s_1,\dots)$. A direct application of part (ii) in Lemma~\ref{hypotheses} implies (\ref{asterix}).
To prove (\ref{obelix}), define the sequence $1=\gamma(0),\gamma(1),\gamma(2),\dots$ s.t. $\gamma(k+1)=\gamma(k)+k$; setting $S_k=s_{\gamma(k)}$, one gets 
 $0 = S_0 = S_1 < S_2 < \dots $.
Given $S_{k+1}\le n<S_{k+2}$, let $Q_n':=A_{S_k+1}\cdots A_n=A_{S_k+1}\cdots A_{S_{k+1}}\cdots A_n$, so that 
\begin{align*}
Q_{n}'&=
A_{s_{\gamma(k)}+1}\cdots A_{s_{\gamma(k+1)}}\cdots A_n
=
Q_{s_{\gamma(k)+1}}\cdots Q_{s_{\gamma(k)+k}}\cdots Q_n.
\end{align*}
Now, if  $S_2\le S_{k+1}\le n< S_{k+2}$, then $k\ge1$ and $s_{\gamma(k)+1}\ge s_2$; by condition ${\bf (C)}$, we know that
$Q_{s_{\gamma(k)+1}},\dots,Q_n\in\mathcal{H}_1\cap\mathcal{H}_2(\Lambda)\cap\mathcal{H}_3(\lambda)$:
applying (\ref{asterix}) together with Lemma \ref{hypotheses} (i) it is necessary that  $Q_n'\in{\mathcal H}_2(\Lambda')\cap{\mathcal H}_3(\lambda^k)$.  
Finally, part (f) of Lemma~\ref{hypothesesbis} also ensures that $Q_n'\in\mathcal{H}_1$: replacing the sequence $(s_0,s_1,\dots)$ by the subsequence  $(S_0,S_1,\dots)$ (if necessary), one may assume that (\ref{obelix}) holds.

\end{proof}

From now on and throughout, we assume (Lemma~\ref{increasing}) the existence of  $0\le\lambda< 1\le \Lambda$ and of a sequence  $0=s_0=s_1<s_2<\cdots$ of integers s.t. for any $n\ge0$ and $k={\bf k}(n)\ge 0$ s.t. $s_{k+1} \le n<s_{k+2}$, 
\begin{equation}\label{hypotplus} 
n\ge s_2\Longrightarrow Q_n\in{\mathcal H}_1\cap{\mathcal H}_2(\Lambda)\cap{\mathcal H}_3(\lambda^k),
\end{equation}
and moreover 
\begin{equation}\label{asterixbis}
s_2\le s_p< s_{k+1}\le n<s_{k+2}\Longrightarrow Q_{s_p}\cdots Q_{s_{k}}Q_n\in\mathcal{H}_2(\Lambda).
\end{equation}

By definition, any $d\times d$ matrix $M\in{\mathcal H}_1$ is associated in an unique way with disjoint nonempty subsets of $\{1,\dots,d\}\times\{1,\dots,d\}$ that are of the form $ \II_h(M)\times  \JJ_h(M)$, for $1\le h\le\kappa(M):=\#{\bf Col}(M)$ and such that
\begin{equation}\label{natural}
\Big\{(i,j)\;;\;M(i,j)\ne0\Big\}=\bigsqcup_{h=1}^{\kappa(M)}  \II_h(M)\times  \JJ_h(M)\quad\hbox{and}\quad  \II_1(M)\supsetneq  \II_2(M)\supsetneq\dots\supsetneq  \II_{\kappa(M)}(M).
\end{equation}
In particular, given any $1\le h,h'\le\kappa(M)$ and $(j,j')\in \JJ_h(M)\times \JJ_{h'}(M)$,
\begin{equation}\label{naturalbis}
h=h'\Longrightarrow \Delta(MU_j)=\Delta(MU_{j'})\ne0
\quad\text{while}\quad h>h'\Longrightarrow \Delta(MU_j)>\Delta(MU_{j'})\ne	0.
\end{equation}

\begin{figure}[H]
\begin{center}
\includegraphics[scale=0.5]{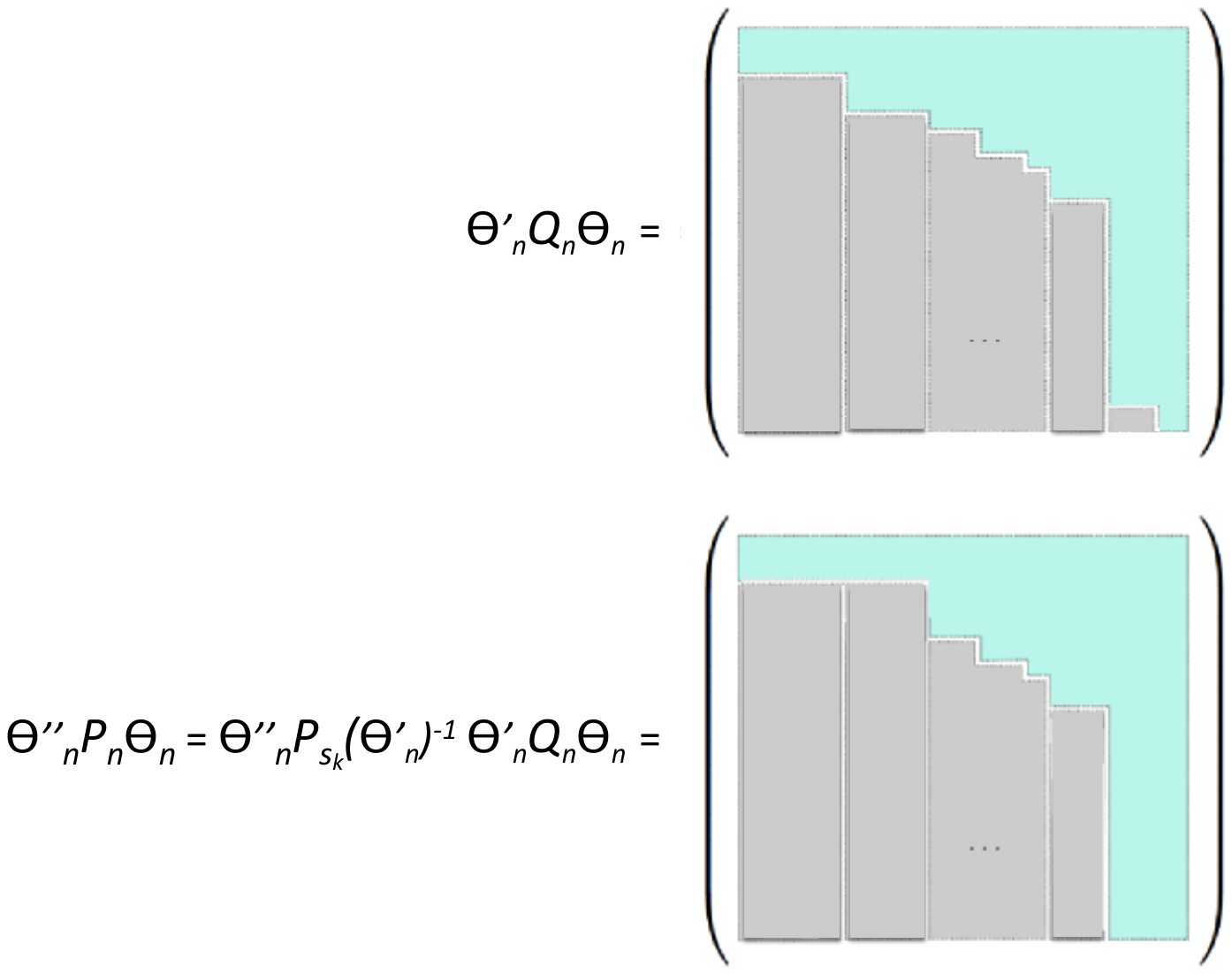}
\caption{\label{RUSSE}\footnotesize\it How the index sets $ I_h'(n)$ and $ J_h'(n)$ are obtained for each $1\le h\le H_n\le\kappa(Q_n)$: here, $\Theta_n$ is a permutation matrix which rearranges the columns in $Q_n$ and $P_n$ so that both $(\mathcal I(Q_n\Theta_nU_j))_{j=1}^d$ and $(\mathcal I(P_n\Theta_nU_j))_{j=1}^d$,  form a non-increasing sequence of row-index sets. Similarly, $\Theta'_n$ (resp. $\Theta''_n$) is a permutation matrix which rearranges the rows in $Q_n$ (resp. $P_n$) so that $(\mathcal I((\Theta'_nQ_n)^*U_j))_{j=1}^d$ (resp. $(\mathcal I((\Theta''_nP_n)^*U_j))_{j=1}^d$),  form a non-decreasing sequence of row-index sets.}
\end{center}
\end{figure}

\noindent The definitions of $1\le H\le d$, of the index domains $ I_h$ and $J_{h}(n)$ ($n\ge s_2$ and $1\le h\le H$) and of the vectors $V_h$ --~as introduced in Theorem~A~-- will be given as a consequence of Theorem~\ref{constant} below. Before that, we  first consider the intermediate index domains $ I_h'(n)$ and $J_h'(n)$, for $1\le h\le H_n$ as follows; 
by part b of Lemma~\ref{hypothesesbis}, given $1\le h\le\kappa(Q_n)$, the row-index set $\mathcal I(P_nU_j)$ is the same for any $j\in\JJ_h(Q_n)$; it may be empty if $h$ exceed some value, say $H_n$: so for any $1\le h\le H_n$, we define
\begin{equation}\label{defIJprime}
 J_h'(n):= \JJ_h(Q_n),
\end{equation}
and $I_h'(n)$ is by definition the common value of $\mathcal I(P_nU_j)$ for any $j\in\JJ_h(Q_n)$. By part a of Lemma~\ref{hypothesesbis}, 
$ \II_1(Q_n)\supsetneq 
\dots\supsetneq  \II_{H_n}(Q_n)$ implies $ I_{1}'(n)\supset 
\dots \supset  I_{H_n}'(n)$ (see Figure~\ref{RUSSE}).

\begin{definition}Given any $1\le h\le H_n$,
\begin{equation}\label{hstarPn}
h\diamond P_n:=\sum_{(i,j)\in I_h'(n)\times J_h'(n)}P_n(i,j)U_iU_j^\star
\end{equation}
(i.e. $h\diamond P_n$ is  obtained from $P_n$ by replacing by zero each entry $P_n(i,j)$, for $(i,j)\not\in I_h'(n)\times J_h'(n)$); hence,  each $h\diamond P_n$ satisfies condition ${\bf (H)}$ together with the identity
$$
P_n=\sum\nolimits_{h=1}^{H_n}h\diamond P_n.
$$
\end{definition}
\begin{lemma}\label{infinitelemma}
Given $s_1<s_{k+1}\le n<s_{k+2}$ (i.e. $\kk(n)=k\ge1$) and $X\in\mathcal{S}_d$ with  $P_nX\ne0$, define 
\begin{equation}\label{BBB}
\hhh_X'(n)=\min\Big\{\;1\le h\le H_n\;;\; (h\diamond P_n)X\ne0\Big\}=\min\Big\{\;1\le h\le H_n\;;\; \mathcal{I}(X)\cap J_h'(n)\ne\emptyset\Big\}\;;
\end{equation}
then, $\Delta(P_nX)=\Delta\big((\hhh_X'(n)\diamond P_n)X\big)$ with
\begin{equation}\label{INFTY}
\mathcal{I}(P_nX)=\mathcal{I}\Big(\big(\hhh_X'(n)\diamond P_n\big)X\Big)=I_{\hhh_X'(n)}'(n).
\end{equation}
\end{lemma}

The matrix $P_n$ is not likely to satisfy condition ${\bf (H)}$ and one can reasonably expect that $\delta(P_n)=+\infty$ for most of ranks $n$. However, it is relevant to look at 
$$
\delta(h\diamond P_n)=\max\Big\{\delta(P_nU_{j_0},P_nU_{j_1})\;;\;j_0,j_1\in J_h'(n)\Big\}
\quad\text{and}\quad\max\Big\{\delta(h\diamond P_n)\;;\;1\le h\le H_n\Big\}\;;
$$ 
we shall prove (see Lemma~\ref{rexpk} below) that the latter maximum tends to $0$ with an exponential rate depending on $k=\kk(n)$. Before proving this convergence, we begin with crucial inequalities.

\begin{lemma}\label{rexpkbis}Let $s_{k+1}\le n<s_{k+2}$ (with $k=\kk(n)\ge1$) and  $X,Y\in\mathcal{S}_d$ such that both $P_nX,P_nY\ne0$; then, 
with $M_{X,Y}:=\max\{\Lambda_X,\Lambda_Y\}$, one has  
\begin{eqnarray*}
(i)&:&\hhh_X'(n)=\hhh_Y'(n)\Longrightarrow \delta\big(P_nX,P_nY\big)\le\lambda^k\big(2\Lambda M_{X,Y}\big)+\delta\big(\hhh_X'(n)\diamond P_n\big)\begin{matrix}\,\\\,\end{matrix}\;;\\
(ii)&:&
\Delta(X)=\Delta(Y)\Longrightarrow \delta\big(P_nX,P_nY\big)\le\lambda^k\big(2\Lambda M_{X,Y}\big)+\delta\big(\hhh_X'(n)\diamond P_n\big)\left(\frac{M_{X,Y}-1}{M_{X,Y}+1}\right).
\end{eqnarray*}

\end{lemma}

\begin{proof}[{\bf Proof}]For $s_{k+1}\le n<s_{k+2}$ with $k\ge1$ and $X\in\mathcal{S}_d$ such that $P_nX\ne0$ we note $h:=\hhh_X'(n)$. We first prove in (\ref{flhglds}) below that each $P_nX(i)$ is closed to $(h\diamond P_n)X(i)$. By definition of $h$, there exists at least one index $j_0\in  J_h'(n)$ such that $X(j_0)\ne0$. On the one hand, because $X(j_0)\ne0$, the definition of $\Lambda_X$ gives $\Vert X\Vert\le \Lambda_XX(j_0)$; on the other hand,  recall that condition (\ref{hypotplus}) means $Q_n\in\mathcal{H}_2(\Lambda)\cap\mathcal{H}_3(\lambda^k)$: hence, for any $h<r\le H_n$ and  $j_1\in J_{r}'(n)$, one has $\Delta(Q_nU_{j_0})>\Delta(Q_nU_{j_1})\ne 0$ and thus (Lemma~\ref{majoration}) for any $1\le i\le d$ (and $n\ge s_2$):
\begin{equation*}
P_n(i,j_1)=
\sum\nolimits_jP_{s_k}(i,j)Q_n(j,j_1)\le 
\sum\nolimits_jP_{s_k}(i,j)\big(\lambda^k\Lambda\big)Q_n(j,j_0)=
\lambda^k\Lambda P_{n}(i,j_0)\;;
\end{equation*}
therefore,
$$
\begin{array}{rcl}\displaystyle\sum\nolimits_{r=h+1}^{H_n}(r\diamond P_n)X(i)=\displaystyle\sum\nolimits_{r=h+1}^{H_n}\displaystyle\sum\nolimits_{j\in J_{r}'(n)}P_n(i,j)X(j)&\le&
\lambda^k\Lambda P_n(i,j_0)\Vert X\Vert\\&\le&\lambda^k\big(\Lambda\Lambda_X\big)P_n(i,j_0)X(j_0),
\end{array}
$$
which allows to write:
\begin{eqnarray*}
(h\diamond P_n)X(i)\le P_nX(i)
&=&\sum\nolimits_{r=h}^{H_n}(r\diamond P_n)X(i)\\
&\le&(h\diamond P_n)X(i)+\sum\nolimits_{r=h+1}^{H_n}(r\diamond P_n)X(i)\\
&\le&(h\diamond P_n)X(i)+\lambda^k\big(\Lambda\Lambda_X\big)P_n(i,j_0)X(j_0).
\end{eqnarray*}
Because $P_n(i,j_0)X(j_0)\le(h\diamond P_n)X(i)$ one finally gets 
\begin{equation}\label{flhglds}
1\le {P_nX(i)\over (h\diamond P_n)X(i)}\le {(h\diamond P_n)X(i)+\lambda^k\big(\Lambda\Lambda_X\big)P_n(i,j_0)X(j_0)\over (h\diamond P_n)X(i)}\le 1+\lambda^k\big(\Lambda\Lambda_X\big).
\end{equation}
Since $h:=\hhh_X'(n)$, Lemma~\ref{infinitelemma} and the  definition of the $\delta$-projective distance give 
$$
\delta\big(P_nX,(h\diamond P_n)X\big)=\max\left\{\log\left(\frac{P_nX(i)\cdot(h\diamond P_n)X(i')}{P_nX(i')\cdot(h\diamond P_n)X(i)}\right)\;;\;i,i'\in  I_h'(n)\right\}
$$ 
and since by a double application of (\ref{flhglds}), 
$$
\frac{P_nX(i)\cdot(h\diamond P_n)X(i')}{P_nX(i')\cdot (h\diamond P_n)X(i)}\le\frac{P_nX(i)}{(h\diamond P_n)X(i)}\le1+\lambda^k\big(\Lambda\Lambda_X\big),
$$
it follows that
\begin{equation}\label{YY*}
\delta\big(P_nX,(h\diamond P_n)X\big)\le\lambda^k\big(\Lambda\Lambda_X\big).
\end{equation}
Now, let  $X,Y\in\mathcal{S}_d$ with $P_nX,P_nY\ne 0$ and satisfying the additional condition ${\bf h}_n'(X)={\bf h}_n'(Y)=:h$. 
Using triangular inequality for $\delta(\cdot,\cdot)$ together with (\ref{YY*})~gives
\begin{eqnarray*}
\delta\big(P_nX,P_nY\big)
&\le&\delta\big(P_nX,(h\diamond P_n)X\big)+\delta\big((h\diamond P_n)X,(h\diamond P_n)Y\big)+\delta\big((h\diamond P_n)Y,P_nY\big)\\
&\le&\lambda^k\big(2\Lambda  M_{X,Y}\big)+\delta\big((h\diamond P_n)B\big),
\end{eqnarray*}
where we have introduced the $d\times 2$ matrix $B=(X\ Y)$.
On the one hand, $h\diamond P_n$ satisfies (by definition) the condition ${\bf (H)}$ and thus, by Lemma~\ref{deltatau}, 
\begin{equation}\label{XY}
\delta\big(P_nX,P_nY\big)\le \lambda^k(2\Lambda M_{X,Y})+\delta(h\diamond P_n)\tau(B).
\end{equation}
This proves part (i) of the lemma, since $\tau(\cdot)$ is bounded by $1$. On the other hand, if one assumes that $\Delta(X)=\Delta(Y)$ then (Lemma~\ref{logLambda})  $\delta(B)\le\log(M_{X,Y}^2)$: because $\tau(B)=\tanh(\delta(B)/4)$, one deduces part (ii) from (\ref{XY}) and the lemma is proved. 

\end{proof}

We are now in position to prove the following key lemma which is the first step for proving assertions (i) and (ii) of Theorem~A.
\begin{lemma}\label{rexpk}There exist two constants $C>0$ and $0<r<1$ 
for which the following properties hold for any $s_{k+1}\le n<s_{k+2}$ (with $k=\kk(n)\ge1$): 
for any $1\le h\le H_n$ 
\begin{equation}\label{AAA}
\max\Big\{\delta(P_nU_{j_0},P_nU_{j_1})\;;\;j_0,j_1\in J_h'(n)\Big\}=\delta(h\diamond P_n)\le C r^{k}\;;
\end{equation}
moreover, given $X\in\mathcal{S}_d$ with  $P_nX\ne0$ 
\begin{equation}\label{Xandc}
j\in  J_{\hhh_X'(n)}'(n)\Longrightarrow \delta(P_nU_j,P_nX)\le \big(C\Lambda_X\big)r^{k}.
\end{equation}

\end{lemma}

\begin{proof}[{\bf Proof}] 
To prove (\ref{AAA}) we proceed by induction over $k\ge1$. Consider  $r$ and~$C$~such that
\begin{equation}\label{EEE}
r:=\max\left\{\lambda,\frac{\Lambda}{\Lambda+1}\right\}
\quad\text{and}\quad
C:=4\Lambda^3
\end{equation}
(in particular $\frac12\le r<1$ and  $C>0$). The case $k=1$ means that $0=s_1<s_2\le n<s_3$, so that $P_n=P_{s_1}Q_n=Q_n$ and condition (\ref{hypotplus}) --~deduced from condition ${\bf (C)}$~-- gives $P_n\in{\mathcal H}_2(\Lambda)$: for an arbitrary $1\le h\le H_n$ (Lemma~\ref{logLambda}) $\delta(h\diamond P_n)\le\log(\Lambda^2)\le \Lambda^2=\frac C{4\Lambda}\le \frac C4\le Cr^1$ and the induction is initialized for  rank $k=1$. 
Suppose (\ref{AAA}) satisfied for rank $k\ge1$, that is $\delta(h\diamond P_n)\le Cr^k$, for each $1\le h\le H_n$ and any $s_{k+1}\le n<s_{k+2}$. Let $s_{k+2}\le n<s_{k+3}$ so that $P_n=P_{s_{k+1}}Q_n$ and take two columns of $Q_n$, say  $X=Q_nU_{j_0}$ and $Y=Q_nU_{j_1}$ with $j_0,j_1\in J_h'(n)$ for an arbitrary $1\le h\le H_n$. Notice that $\Delta(X)=\Delta(Y)$ with $P_{s_{k+1}}X,P_{s_k+1}Y\ne0$ so that 
$$
\hhh_X'\big(s_{k+1}\big)=\hhh_Y'\big(s_{k+1}\big)=:\ell.
$$
Since $\kk(s_{k+1})=k$ the induction hypothesis  gives  $\delta(\ell\diamond P_{s_{k+1}})\le Cr^k$; moreover, condition (\ref{hypotplus}) implies that $\Lambda_X$ and $\Lambda_Y$ are bounded by $\Lambda$ and thus $M_{X,Y}\le\Lambda$:  using part (ii) of Lemma~\ref{rexpkbis} 
\begin{align*}
\delta(P_nU_{j_0},P_nU_{j_1})
&=\delta(P_{s_{k+1}}X,P_{s_{k+1}}Y)\\
&\le \lambda^{k}(2\Lambda M_{X,Y})+\delta(\ell\diamond P_{s_{k+1}})\left(\frac{M_{X,Y}-1}{M_{X,Y}+1}\right)\\
&\le\lambda^{k}(2\Lambda^2)+Cr^k\left(\frac{\Lambda-1}{\Lambda+1}\right)=\lambda^{k}(2\Lambda^2)+Cr^{k+1}-Cr^k\left(\frac1{\Lambda+1}\right);
\end{align*}
but (\ref{EEE}) gives $\lambda^{k}(2\Lambda^2)\le r^{k}(2\Lambda^2)=r^{k}\left(\frac C{2\Lambda}\right)\le r^{k}\left(\frac C{\Lambda+1}\right)=Cr^k\left(\frac1{\Lambda+1}\right)$ so that $\delta(P_nU_{j_0},P_nU_{j_1})\le Cr^{k+1}$:  the induction holds and (\ref{AAA}) is established.

To prove (\ref{Xandc})  consider $s_{k+1}\le n<s_{k+2}$ for $k\ge1$ and $X\in\mathcal{S}_d$ s.t. $P_nX\ne0$. Given $Y=U_j$, for $j\in J_{\hhh_X'(n)}'(n)$ it is necessary (and sufficient) that $\hhh_{Y}'(n)=\hhh_X'(n)=:h$; moreover, $\Lambda_{Y}=1$ and $M_{X,Y}=\Lambda_X$. From (\ref{AAA}) there exists a constant $C>0$ s.t. $\delta(h\diamond P_n)\le Cr^k$ and thus using  part (i) of Lemma~\ref{rexpkbis} gives
\begin{eqnarray*}\delta(P_nX,P_nU_j)
\le\lambda^k(2\Lambda\Lambda_X)+\delta(h\diamond P_n)\le\left(2\Lambda\Lambda_X+C\right)r^k\le\big(C'\Lambda_X\big)r^k,
\end{eqnarray*}
where $C'=2\Lambda+C$ : this proves (\ref{Xandc}).

\end{proof}

\subsection{Proof of Theorem~A} Recall that $\mathcal{A}=(A_1,A_2,\dots)$ (where each $A_i$ is a $d\times d$ matrix with non negative entries) is supposed to satisfy condition (\ref{hypotplus}) w.r.t. $0\le \lambda<1\le \Lambda$ and the sequence of integers $0=s_0=s_1<s_2<\cdots$. The (compact) simplex $\mathcal{S}_d$ of the nonnegative $d\times 1$ vectors $V$ s.t. $\Vert V\Vert=1$ is endowed with the topology induced from the normed topology on~$\RRR^d$. The set $\mathcal{S}_d(\mathcal{A})$ of the {\it projective limit vectors} of $\mathcal{A}$ is made of the $W\in\mathcal{S}_d$ for which there exists a sequence  $1\le m_1<m_2<\cdots$ of integers and a sequence $(j_1,j_2,\dots)$  in $\{1,\dots,d\}$ such that $X_t=P_{m_t}U_{j_t}/\Vert P_{m_t}U_{j_t}\Vert\to W$ as ${t\to+\infty}$; in other words,
\begin{equation}\label{cvg}
\lim_{t\to+\infty}\Vert X_t-W\Vert=0.
\end{equation}
The following lemma ensures the convergence in (\ref{cvg}) to imply that $X_t$ is {\it trapped} in the face of the projective limit vector $W$: more precisely --~for $t$ large enough
\begin{equation}\label{t0}
X_t\in\mathcal{S}_{\mathcal{I}(W)}=\{X\in\mathcal{S}_d\;;\; \mathcal{I}(X)=\mathcal{I}(W)\}.
\end{equation}

\begin{lemma}[Trapping lemma]\label{trapping}
(i) : Let $X_{n,j}:=P_nU_j/\Vert P_nU_j\Vert$, for $n\ge s_2$ ($1\le j\le d$):~then, 
\begin{equation}\label{projectivecvgtris}
X_{n,j}(i)\ne0\Longrightarrow X_{n,j}(i)\ge 1/\Lambda\;;
\end{equation}
(ii) : given any $W\in\mathcal{S}_d(\mathcal{A})$, one has 
\begin{equation}\label{projectivecvgbis}
W(i)\ne 0\Longrightarrow W(i)\ge 1/\Lambda\;;
\end{equation}
(iii) : for any sequence $1\le m_1<m_2<\cdots$ (resp. $j_1,j_2,\dots$) made of integers (resp. indices in $\{1,\dots,d\}$) and $W_1,W_2,\dots\in \mathcal{S}_d(\mathcal{A})$, one has the equivalence
\begin{equation}\label{projectivecvg}
\lim_{t\to+\infty}\left\Vert X_{m_t,j_t}-W_t\right\Vert=0\iff \lim_{t\to+\infty}\delta\left(X_{m_t,j_t},W_t\right)=0.
\end{equation}
\end{lemma}
\begin{proof}[{\bf Proof}] (i) : For $s_{k+1}\le n< s_{k+2}$ ($k\ge1$): because  
$P_{n}=Q_{s_2}\cdots Q_{s_k}Q_{n}$, we know from condition  (\ref{hypotplus})-(\ref{asterixbis}) that  $X_{n,j}(i)\ne0$ implies $1=\Vert X_{n,j}\Vert\le\Lambda\cdot X_{n,j}(i)$, so that $X_{n,j}(i)\ge 1/\Lambda$.

(ii) : Let $W\in\mathcal{S}_d(\mathcal{A})$ and $1\le m_1<m_2<\cdots$  s.t. $\Vert X_{m_t,j_t}-W\Vert\to 0$ as $t\to+\infty$; from part~(i) one deduces that  $W(i)\ne 0$ implies $W(i)\ge 1/\Lambda$.

(iii) : If $\Vert X_{m_t,j_t}-W_t\Vert\to 0$ as $t\to+\infty$ then, an immediate consequence of (i)  and (ii) is that $X_{m_t,j_t}\in\mathcal{S}_{\mathcal{I}(W_t)}$ ($t$ large enough): hence, it is licit to apply  the upper bound in (\ref{doublebound}) of Proposition~\ref{face} so that $\delta(X_{m_t,j_t},W_t)\le \Lambda\cdot \Vert X_{m_t,j_t}-W_t\Vert$ and  $\delta(X_{m_t,j_t},W_t)\to0$. Conversely, if $\delta(X_{m_t,j_t},W_t)\to0$, then $\delta(X_{m_t,j_t},W_t)$ is finite ($t$ large enough) and $X_{m_t,j_t}\in\mathcal{S}_{\mathcal{I}(W_t)}$: the lower bound in (\ref{doublebound}) of Proposition~\ref{face} gives $ \Vert X_{m_t,j_t}-W_t\Vert\le d\cdot\delta(X_{m_t,j_t},W_t)$ and  $\Vert X_{m_t,j_t}-W_t\Vert\to0$.

\end{proof}
\begin{lemma}
The set $\mathcal{S}_d(\mathcal{A})$ of the projective limit vectors is finite with $1\le \#\mathcal{S}_d(\mathcal{A})\le d$.
\end{lemma}

\begin{proof}[{\bf Proof}]Let $C>0$ and $0<r<1$ be the two constants given by Lemma~\ref{rexpk} and 
fix $\varepsilon>0$, $q\ge1$ such that $C\Lambda r^q\le\varepsilon$; moreover, assume $X_t:=P_{m_t}U_{j_t}/\Vert P_{m_t}U_{j_t}\Vert\to V$, as $t\to+\infty$: according to (\ref{t0}) and by continuity of $\delta(\cdot,\cdot)$, it is licit to consider that $\lim_{t\to+\infty}\delta(X_t,V)=0$. For  $m_t\ge s_{q+1}$, write $X_t={P_{s_{q}}Y_t/\Vert P_{s_{q}}Y_t\Vert}$ where 
$$
Y_t=
\begin{cases}
Q_{m_t}U_{j_t}&\text{if $m_t<s_{q+2}$ (i.e. $\kk(m_t)=q$)}\;;\\
Q_{s_{q+1}}\cdots Q_{s_{\kk(m_t)}}Q_{m_t}U_{j_t}&\text{if $m_t\ge s_{q+2}$ (i.e. $\kk(m_t)\ge q+1$)}.
\end{cases}
$$
By definition $X_t\ne0$ and  (\ref{Xandc}) in Lemma~\ref{rexpk} ensures the existence of a column index $j_t'$ (actually $j_t'\in J_{\hhh_{s_q}'}(Y_t)$) for which
\begin{equation}\label{dhjkflqs}
\delta(P_{s_{q}}U_{j_t'},X_t)=\delta(P_{s_{q}}U_{j_t'},P_{s_{q}}Y_t)\le C\Lambda_{Y_t} r^q.
\end{equation}
Since (conditions  (\ref{hypotplus})-(\ref{asterixbis})) both $Q_{m_t}$ and $Q_{s_{q+1}}\cdots Q_{s_{\kk(m_t)}}Q_{m_t}$ belong to $\mathcal{H}_2(\Lambda)$, one gets  that $\Lambda_{Y_t}\le \Lambda$ and  (\ref{dhjkflqs}) implies 
$\delta(P_{s_{q}}U_{j_t'},X_t)\le\varepsilon$. This last inequality being valid for any $t$ such that  $m_t\ge s_{q+1}$, there exists a column index $j_0$, for which $j_t'=j_0$ for infinitely many $t$, meaning that $\delta(P_{s_{q}}U_{j_0},V)\le \varepsilon$. Suppose --~for a contradiction~-- the existence of at least  $d+1$ projective limit vectors, say  $V_1,\dots,V_{d+1}$, with $\varepsilon$ sufficiently small so that $\delta\left(V_{\ell},V_{\ell'}\right)\ge3\varepsilon$ as soon as $\ell\ne \ell'$: the sets
$\{X\in\mathcal{S}_d\;;\;\delta(X,V_\ell)\le\varepsilon\}$
are disjoint for $\ell=1,\dots,d+1$ and each ones must contain a probability vector of the form $P_{s_{q}}U_j/\Vert P_{s_{q}}U_j\Vert$: this is a contradiction.

\end{proof}

\begin{lemma}\label{distance}
Let $0<\delta_0\le+\infty$ be the smallest $\delta$-distance between two distinct projective limit vectors (and $\delta_0=+\infty$ if $\#\mathcal{S}_d(\mathcal{A})=1$); for $0<\varepsilon\le \delta_0/3$ there exists $N(\varepsilon)\ge1$~s.t.    
\begin{equation}\label{Vhnbis}
\forall n\ge N(\varepsilon),\;1\le \forall h\le H_n,\;\exists! W_h(n)\in\mathcal{S}_d(\mathcal{A}),\;j\in J_h'(n)\Longrightarrow\delta\big(P_nU_j,W_h(n)\big)< \varepsilon\;;
\end{equation}
moreover, $\delta\big(P_nU_j,W_h(n)\big)< \varepsilon$ ($\le+\infty$) implies $\Delta(P_nU_j)=\Delta(W_h(n))$  and from the definition of $J_h'(n)$,
\begin{equation}\label{Vhntris}
\Delta(W_1(n))\ge\cdots\ge\Delta(W_{H_n}(n)).
\end{equation}
\end{lemma}

\begin{proof}[{\bf Proof}] 
Let $0<\varepsilon\le +\infty$ be arbitrary given  and define
$\mathcal{S}_d^\varepsilon(\mathcal{A})$ the set of the vectors $X\in\mathcal{S}_d$ for which there exists $V\in\mathcal{S}_d(\mathcal{A})$  for which $\delta(X,V)< \varepsilon$. We claim the existence of a rank $N(\varepsilon)\ge1$ s.t. each $P_nU_j/\Vert P_nU_j\Vert$, for $1\le j\le d$ and $n\ge N(\varepsilon)$, belongs to $\mathcal{S}_d^\varepsilon(\mathcal{A})$. Suppose --~for a contradiction~-- the existence of an increasing sequence $m_1<m_2<\cdots$ of ranks together with a sequence $(j_1,j_2,\dots)$ of indices  with  $X_t:=P_{m_t}U_{j_t}/\Vert P_{m_t}U_{j_t}\Vert\in\mathcal{S}_d\setminus\mathcal{S}_d^\varepsilon(\mathcal{A})$, for any $t\ge1$: taking a subsequence if necessary, this means the existence of at least one projective limit vector $X_*$ for which the normed convergence $X_t\to X_*$ holds as $t\to+\infty$. However, by the {\it trapping lemma}  (see  (\ref{projectivecvg}) in Lemma~\ref{trapping}) we know that $\Vert X_t-X_*\Vert\to 0$ is equivalent to the projective convergence $\delta(X_t,X_*)\to0$ ensuring the existence of  $t_0\ge1$ s.t. $\delta(X_t,X_*)< \varepsilon$ for any $t\ge t_0$: this is a contradiction, because $X_*\in\mathcal{S}_d(\mathcal{A})$, while $\delta(X_t,V)\ge \varepsilon$, for any $t\ge1$ and any $V\in \mathcal{S}_d(\mathcal{A})$.

Now, assume in addition that $\varepsilon\le \delta_0/3$, let $q$ satisfy $Cr^q\le\varepsilon$ (with $C,r$ given by Lemma~\ref{rexpk}) and choose $N(\varepsilon)\ge s_{q+1}$: by the first part of the argument, for any $n\ge N(\varepsilon)$ and any $1\le j\le d$, there exists a unique projective limit vector $W_j'(n)\in\mathcal{S}_d(\mathcal{A})$ for which $\delta\big(P_nU_j,W_j'(n)\big)< \varepsilon$. If  $j_0$ and $j_1$ are arbitrary taken in $J_h'(n)$, it follows from  (\ref{AAA}) in Lemma~\ref{rexpk} that $\delta(P_nU_{j_0},P_nU_{j_1})\le\varepsilon$, and $W_{j_0}'(n)=W_{j_1}'(n)$: by definition  $W_h(n)$ stands for this common projective limit vector.

\end{proof}

Let  $\gotL$ be the set of the   $X=(x_1,\dots,x_n)$, with $x_i$ unspecified and where by abuse of notation $\#X:=n\ge1$ stands for the length of $X$.  We shall consider the map $\Xi:\gotL\to\gotL$ defined by setting
$$
\Xi(x_1,\dots,x_n)=(x_{\varphi(1)},\dots,x_{\varphi(k_*)}),
$$
where 
$\varphi(1)=1$ while 
$\varphi(k+1)$ is the minimum of the $\varphi(k)<i\le n$ s.t. $x_i\ne x_{\varphi(k)}$, provided that $(x_{\varphi(k)},\dots,x_n)\ne(x_{\varphi(k)},\dots,x_{\varphi(k)})$; here $k_*$ is the minimum of the $k\ge1$ s.t. $(x_{\varphi(k)},\dots,x_n)=(x_{\varphi(k)},\dots,x_{\varphi(k)})$: for instance $\Xi(1,1,1,2,3,3,3,1,1)=(1,2,3,1)$ with $k_*=4$. According to Lemma~\ref{distance} it is licit to define 
\begin{equation}\label{defHN*}
H:=\min\Big\{\#\Xi\big(W_1(s_k),\dots,W_{H_{s_k}}(s_k)\big)\;;\;s_k\ge N(\delta_0/3)\Big\}.
\end{equation}

\begin{remark}\label{thelast}
The definition of the index sets $J_h(n)$ and of the vectors $V_h$ (for $1\le h\le H$) --~that are the main ingredients of Theorem~A~-- are given in Theorem~\ref{constant} below. Let's mention now that the argument for Theorem~\ref{constant} stands  on the simple remark that for $X=(x_1,\dots,x_n)$ any ordered lists and $\varphi:\{1,\dots,n\}\to\{1,\dots,n\}$ any nondecreasing map, one has 
\begin{equation}\label{propXi}
\#\Xi(x_{\varphi(1)},\dots,x_{\varphi(n)})\le \#\Xi(x_1,\dots,x_n).
\end{equation}
This is a motivation for the introduction of the maps $\varphi_n(\cdot):\{1,\dots,H_n\}\to\{1,\dots,H_{s_k}\}$ and $\psi_k(\cdot):\{1,\dots,H_{s_{k+2}}\}\to\{1,\dots,H_{s_k}\}$ in (\ref{defhnell})
and (\ref{defhnellbis}) below and the associated Lemma~\ref{bidule}.
\end{remark}

Given  $s_{k+1}\le n<s_{k+2}$ and $1\le h\le H_n$, one has $J_h'(n)=\JJ_h(Q_n)$ and:   
$$
j_0,j_1\in J_h'(n)\Longrightarrow \Delta(Q_nU_{j_0})=\Delta(Q_nU_{j_1})=:X_{n,h}.
$$
We suppose that $k\ge2$ because $J_h'(n)$ is defined for $n\ge s_2$ and we need the map $X\mapsto\hhh_X'(s_k)$ --~as introduced in (\ref{BBB})~-- to make sense. The map $\varphi_n(\cdot):\{1,\dots,H_n\}\to\{1,\dots,H_{s_k}\}$ is defined by setting
\begin{equation}\label{defhnell}
\varphi_n(h)=\hhh_{X_{n,h}}'(s_k)=\min\Big\{1\le \ell\le H_{s_k}\;;\;\mathcal{I}(X_{n,h})\cap J_\ell'(s_k)\ne\emptyset\Big\}.
\end{equation}
One defines $\psi_k(\cdot):\{1,\dots,H_{s_{k+2}}\}\to\{1,\dots,H_{s_k}\}$ analogously to $\varphi_n(\cdot)$. 
Recall that $P_n=P_{s_k}Q_n$, while 
$P_{s_{k+2}}=P_{s_k}Q_{s_{k+1}}Q_{s_{k+2}}$: clearly, for any  $1\le h\le H_{s_{k+2}}$, 
$$
j_0,j_1\in J_h'(s_{k+2})\Longrightarrow \Delta(Q_{s_{k+1}}Q_{s_{k+2}}U_{j_0})=\Delta(Q_{s_{k+1}}Q_{s_{k+2}}U_{j_1})=:Y_{k,h}.
$$
By definition
\begin{equation}\label{defhnellbis}
\psi_k(h)=\hhh_{Y_{k,h}}'(s_k)=\min\Big\{1\le \ell\le H_{s_k}\;;\;\mathcal{I}(Y_{k,h})\cap J_\ell(s_k)\ne\emptyset\Big\}.
\end{equation}

\begin{lemma}\label{bidule} For any $s_{k+1}\le n<s_{k+2}$ (with $k\ge2$), one has the following propositions: 
\begin{eqnarray}
\label{increasingh_n}
\varphi_n(1)\le\dots\le \varphi_n(H_n)\quad&\text{and}&\quad \psi_k(1)\le\dots\le \psi_k(H_{s_{k+2}})\;;\\
\label{Xmn}
(j_0,j_1)\in J_{\varphi_n(h)}'(s_k)\times J_h'(n)&\Longrightarrow& \delta\big(P_{s_k}U_{j_0},P_nU_{j_1}\big)\le (C\Lambda) r^{k-1}\;;\\
\label{Xmnbis}
(j_0,j_1)\in J_{\psi_k(h)}'(s_k)\times J_{h}'(s_{k+2})&\Longrightarrow& \delta\big(P_{s_k}U_{j_0},P_{s_{k+2}}U_{j_1}\big)\le (C\Lambda) r^{k-1}\;;
\end{eqnarray}
moreover, if both inequalities $s_{k}\ge N(\delta_0/3)$ and $(C\Lambda) r^{{k}-1}\le \delta_0/3$ are satisfied, then
\begin{eqnarray}
\label{ellhell}
\big(W_{1}(n),\dots,W_{H_n}(n)\big)&=&\big(W_{\varphi_n(1)}(s_k),\dots,W_{\varphi_n(H_n)}(s_k)\big)\;;\\
\label{ellhellbis}
\big(W_{1}(s_{k+2}),\dots,W_{H_{s_{k+2}}}(s_{k+2})\big)&=&\big(W_{\psi_k(1)}(s_k),\dots,W_{\psi_k(H_{s_{k+2}})}(s_k)\big).
\end{eqnarray}\end{lemma}

\begin{proof}[{\bf Proof}]If  $1\le h_0<h_1\le H_n$ and $(j_0,j_1)\in J'_{h_0}(n)\times J'_{h_1}(n)$ then, by definition of the sets $J'_h(n)$, $\Delta(Q_nU_{j_0})>\Delta(Q_nU_{j_1})$ and thus by minimality, it follows that $\varphi_n(h_0)\le \varphi_n(h_1)$. Suppose now that $1\le h_0<h_1\le H_{s_{k+2}}$ and $(j_0,j_1)\in J'_{h_0}(s_{k+2})\times J'_{h_1}(s_{k+2})$, then $\Delta(Q_{s_{k+2}}U_{j_0})>\Delta(Q_{s_{k+2}}U_{j_1})$, thus $\Delta(Q_{s_{k+1}}Q_{s_{k+2}}U_{j_0})\ge\Delta(Q_{s_{k+1}}Q_{s_{k+2}}U_{j_1})$ and by minimality, it follows that $\psi_k(h_0)\le \psi_k(h_1)$.

For part (\ref{Xmn}) --~and (\ref{Xmnbis}) similarly~-- let $j_1\in J_{h}'(n)$, where $1\le h\le H_n$. It is licit to apply (\ref{Xandc}) in Lemma~\ref{rexpk}: indeed,  $P_{s_k}Q_nU_{j_1}\ne0$ (by definition of $J_{h}'(n)$) and thus
$$
\delta\big(P_{s_k}U_{j_0},P_nU_{j_1}\big)
=\delta\big(P_{s_k}U_{j_0},P_{s_k}Q_nU_{j_1}\big)
\le \big(C\Lambda\big) r^{{k}-1}
$$ 
(here we used the fact that $\kk(s_k)=k-1$ and condition (\ref{hypotplus} ) ensuring that $\Lambda_{Q_nU_{j_1}}\le\Lambda$). 

To prove (\ref{ellhell}) --~and (\ref{ellhellbis}) similarly~-- assume $s_{k}\ge N\left(\delta_0/3\right)$ and $(C\Lambda) r^{{k}-1}\le \delta_0/3$. Consider 
$j_0\in J'_{\varphi_n(h)}(s_k)$ and $j_1\in J_{h}'(n)$, for $1\le h\le H_n$; by definition of $\varphi_n(h)$  and according to (\ref{Xmn}) we know that
$\delta(P_{s_k}U_{j_0},P_{n}U_{j_1})\le  (C\Lambda)r^{k-1}\le \delta_0/3$.
However, since  $s_k\ge N(\delta_0/3)$,  Lemma \ref{distance} ensures that both  $\delta(P_{s_k}U_{j_0},W_{\varphi_n(h)}(s_k))$ and $\delta(P_{n}U_{j_1},W_{h}(n))$ are strictly upper bounded by $\delta_0/3$: for  $\delta_0$ being the minimal $\delta$-distance between two different projective limit vectors, it is necessary that $W_{\varphi_n(h)}(s_k)=W_{h}(n)$ and (\ref{ellhell}) is proved.

\end{proof}

\begin{theorem}\label{constant}
Let $1\le H\le d$ be defined in (\ref{defHN*}); then, there exists a (surjective but not necessarily injective) finite sequence of probability vectors $V_1,\dots,V_H$ in $\mathcal{S}_d(\mathcal{A})$ satisfying 
\begin{equation}\label{hhhhh}
\Delta(V_1)\ge \cdots \ge\Delta(V_H)\quad\text{while}\quad V_{h-1}\ne V_h\quad\text{(for any $1< h\le H$)}
\end{equation}
such that for any $n$ (large enough) there exists a partition
$$
\{1,\dots,H_n\}=E_1(n)\sqcup\cdots\sqcup E_{H}(n)
$$ 
into $H$ nonempty sets $E_{1}(n),\dots,E_H(n)$ for which 
\begin{eqnarray}\label{Vii}
(\ell,\ell')\in E_{h-1}(n)\times E_{h}(n)\Longrightarrow \ell<\ell'\quad\text{and}\quad
\ell\in E_{h}(n)\Longrightarrow W_{\ell}(n)=V_h\;;
\end{eqnarray}
in other words $\Xi(W_{1}(n),\dots,W_{H_n}(n))=(V_1,\dots,V_H)$ and (with a rough representation)
\begin{equation}\label{Vi}
(W_{1}(n),\dots,W_{H_n}(n))=(\underbrace{V_1,\dots,V_1}_{E_1(n)},\underbrace{V_2,\dots,V_2}_{E_2(n)},\dots,\underbrace{V_H,\dots,V_H}_{E_H(n)})\;;
\end{equation}
\end{theorem}

\begin{proof}[{\bf Proof}]Let  $k_0\ge1$ be the minimal integer satisfying the three constraints
\begin{equation}\label{minimality}
(a)\;:\;s_{k_0}\ge N(\delta_0/3)\;;\;\; (b)\;\;:\;(C\Lambda)r^{{k_0}-1}\le\delta_0/3\;;\;\;
(c)\;:\;\#\Xi\big(W_1(s_{k_0}),\dots,W_{H_{s_{k_0}}}(s_{k_0})\big)=H.
\end{equation}
The condition (\ref{minimality}-a) is needed to ensure (Lemma~\ref{distance}) the existence  of the projective limit vectors $W_1(n),\dots W_{H_n}(n)$,  for any $n\ge s_{k_0}$ : hence, by the condition (\ref{minimality}-c) over $k_0$ together with the definition of $H$ in (\ref{defHN*}), it is licit to fix $V_1,\dots,V_H\in\mathcal{S}_d(\mathcal{A})$ s.t.
\begin{equation}\label{DefVh}
\Xi\big(W_1(s_{k_0}),\dots,W_{H_{s_{k_0}}}(s_{k_0})\big)=(V_1,\dots,V_H).
\end{equation}
By (\ref{Vhntris}) in Lemma \ref{distance} we know  that $\Delta(V_1)\ge\cdots\ge \Delta(V_H)$, while (definition of $\Xi$) $V_{h-1}\ne V_h$ for any $1<h\le H$.
The theorem holds as soon as  for any $k\ge k_0$ and any  $s_{k+1}\le n<s_{k+2}$,
\begin{equation}\label{Kn}
\Xi\big(W_1(n),\dots,W_{H_n}(n)\big)=(V_1,\dots,V_H)
\end{equation}
(then, the partition $\{1,\dots,H_n\}=E_1(n)\sqcup\cdots\sqcup E_{H}(n)$ is completely determined by (\ref{Kn}), according to the specification in (\ref{Vii})).

$\bullet$ To prove (\ref{Kn}), we begin with an induction showing that  for any $k\ge k_0$
\begin{equation}\label{Ksksss}
\#\Xi\big(W_{1}(s_k),\dots,W_{H_{s_k}}(s_k)\big)=H.
\end{equation}
The initialization being satisfied  for  $k_0$ --~see (\ref{minimality}-c)~-- let $k\ge k_0$ for which (\ref{Ksksss}) is satisfied as well: because (\ref{minimality}-b) ensures $(C\Lambda)r^{k_0-1}\le\delta_0/3$, it is licit to use 
 (\ref{ellhell}) in  Lemma~\ref{bidule}, so that
\begin{equation}\label{iiiiisss}
\Xi\big(W_{1}(s_{k+1}),\dots,W_{H_{s_{k+1}}}(s_{k+1})\big)=\Xi\big(W_{\varphi_{s_{k+1}}(1)}(s_k),\dots,W_{\varphi_{s_{k+1}}(H_{s_{k+1}})}(s_k)\big).
\end{equation}
Using the fact (see (\ref{increasingh_n}) in Lemma~\ref{bidule}) 
that $\varphi_{s_{k+1}}(\cdot)$ is nondecreasing together with the definition of $H$ in (\ref{defHN*}), one deduces from (\ref{iiiiisss})  --~and the induction hypothesis over $k$~-- that 
\begin{align*}
H\le \#\Xi\big(W_{1}(s_{k+1}),\dots,W_{H_{s_{k+1}}}(s_{k+1})\big)
&\le \#\Xi\big(W_{1}(s_k),\dots,W_{H_{s_k}}(s_k)\big)=  H
\end{align*}
and  (\ref{Ksksss}) is inductive.

$\bullet$ We now prove by induction that for any $k\ge k_0$, 
\begin{equation}\label{Kskss}
\Xi\big(W_{1}(s_k),\dots,W_{H_{s_k}}(s_k)\big)=(V_1,\dots,V_H).
\end{equation}
The initialization holds for $k_0$ due to definition of the projective limit vectors $V_1,\dots,V_H$ in (\ref{DefVh}). Let $k\ge k_0$ s.t. (\ref{Kskss}) holds: then, by (\ref{Ksksss}) and (\ref{iiiiisss}), 
\begin{equation}\label{THEENDsss}
H=\#\Xi\big(W_{1}(s_{k+1}),\dots,W_{H_{s_{k+1}}}(s_{k+1})\big)
=\#\Xi\Big(W_{\varphi_{s_{k+1}}(1)}(s_k),\dots,W_{\varphi_{s_{k+1}}(H_{s_{k+1}})}(s_k)\Big).
\end{equation}
However, by the induction hypothesis $\Xi(W_1(s_k),\dots,W_{H_{s_k}})=(V_1,\dots,V_H)$ and for  $\varphi_{s_{k+1}}(\cdot)$ being  nondecreasing,  (\ref{THEENDsss}) implies that  
$$
\Xi\Big(W_{\varphi_{s_{k+1}}(1)}(s_k),\dots,W_{\varphi_{s_{k+1}}(H_{s_{k+1}})}(s_k)\Big)=(V_1,\dots,V_H)\;;
$$
finally, using again (\ref{iiiiisss}) one concludes
$$
\Xi\big(W_{1}(s_{k+1}),\dots,W_{H_{s_{k+1}}}(s_{k+1})\big)=(V_1,\dots,V_H),
$$
which means that (\ref{Kskss}) is inductive.

$\bullet$ We shall now prove that (\ref{Kn}) holds,  for any $s_{k+1}\le n<s_{k+2}$ ($k\ge k_0$).
On the one hand, recall that  (\ref{minimality}-b) ensures $(C\Lambda)r^{k_0-1}\le\delta_0/3$: hence, according to (\ref{ellhell}) in  Lemma~\ref{bidule},
\begin{equation}\label{iiiii}
\big(W_{1}(n),\dots,W_{H_{n}}(n)\big)=\big(W_{\varphi_{n}(1)}(s_k),\dots,W_{\varphi_{n}(H_{n})}(s_k)\big).
\end{equation}
Because  (see (\ref{increasingh_n}) in Lemma~\ref{bidule})
$\varphi_{n}(\cdot)$ is nondecreasing, it follows from (\ref{Ksksss}) that 
\begin{equation}\label{increasingh_nbis}
\#\Xi\big(W_{1}(n),\dots,W_{H_{n}}(n)\big)
\le \#\Xi\big(W_{1}(s_k),\dots,W_{H_{s_k}}(s_k)\big)=  H.
\end{equation}
On the other hand (\ref{Kskss}) together with (\ref{ellhellbis}) in Lemma~\ref{bidule} give
\begin{equation*}
(V_1,\dots,V_H)=\Xi\big(W_{1}(s_{k+2}),\dots,W_{H_{s_{k+2}}}(s_{k+2})\big)=\Xi\big(W_{\psi_k(1)}(s_k),\dots,W_{\psi_k(H_{s_{k+2}})}(s_k)\big)\;;
\end{equation*}
in particular, this means the existence of a nondecreasing $\alpha:\{1,\dots,H\}\to\{1,\dots,H_{s_{k+2}}\}$~s.t. 
\begin{equation}\label{jjj}
\big(W_{\psi_k\circ\alpha(1)}(s_k),\dots,W_{\psi_k\circ\alpha(H)}(s_k)\big)
=(V_1,\dots,V_H).
\end{equation}
We claim that (\ref{Kn}) will be established (and the theorem as well) as soon has we have proved the existence of a nondecreasing  $\beta:\{1,\dots,H\}\to\{1,\dots,H_{n}\}$ for which  
\begin{equation}\label{jjjjj}
\big(W_{\varphi_n\circ\beta(1)}(s_k),\dots,W_{\varphi_n\circ\beta(H)}(s_k)\big)=
\big(W_{\psi_k\circ\alpha(1)}(s_k),\dots,W_{\psi_k\circ\alpha(H)}(s_k)\big).
\end{equation}
Indeed, (\ref{jjjjj}) together with (\ref{jjj}) gives $\big(W_{\varphi_n\circ\beta(1)}(s_k),\dots,W_{\varphi_n\circ\beta(H)}(s_k)\big)=(V_1,\dots,V_H)$ and thus (\ref{Kskss}) implies $\Xi\big(W_{\varphi_{n}(1)}(s_k),\dots,W_{\varphi_{n}(H_{n})}(s_k)\big)=(V_1,\dots,V_H)$: finally, with  (\ref{iiiii}) one concludes
$$
\Xi\big(W_1(n),\dots,W_{H_n}(n)\big)=(V_1,\dots,V_H).
$$

$\bullet$ To prove  the existence of  the nondecreasing map $\beta(\cdot)$  satisfying (\ref{jjjjj}), let $1\le h_0\le h_1\le H_{s_{k+2}}$ and take $j_i\in J_{\alpha(h_i)}'(s_{k+2})$; because $\alpha:\{1,\dots,H\}\to\{1,\dots,H_{s_{k+2}}\}$ is a nondecreasing map, and from the definition of the sets $J_h'(s_{k+2})$, we know that 
$\Delta(Q_{s_{k+2}}U_{j_0})\ge\Delta(Q_{s_{k+2}}U_{j_1})$. 
Now, recall that 
$$
Q_{s_{k+1}}Q_{s_{k+2}}=Q_nA_{n+1}\cdots A_{s_{k+2}}\;;
$$
here it is crucial that $Q_n\in{\mathcal H}_1$: indeed, according to part (d) of Lemma~\ref{hypothesesbis}, we know that 
$$
\Delta\big(Q_{s_{k+1}}Q_{s_{k+2}}U_{j_i}\big)=\Delta\big(Q_nA_{n+1}\cdots A_{s_{k+2}}U_{j_i}\big)\in\Big\{\Delta(Q_nU_{1}),\dots,\Delta(Q_nU_{d})\Big\}
$$
ensuring the existence of  $1\le j_0',j_1'\le d$ s.t. 
\begin{equation}\label{kkkkk}
\Delta(Q_nU_{j_0'})=\Delta(Q_{s_{k+1}}Q_{s_{k+2}}U_{j_0})\ge \Delta(Q_{s_{k+1}}Q_{s_{k+2}}U_{j_1})=\Delta(Q_nU_{j_1'}).
\end{equation}
Given  $1\le \beta(h_i)\le H_N$ s.t. $j_i'\in J'_{\beta(h_i)}(n)= \JJ_{\beta(h_i)}(Q_n)$, since $\II_h(Q_n)\supsetneq \II_{h'}(Q_n)$ for any $1<h<h'\le H_n$, the inequality in (\ref{kkkkk}) implies $\beta(h_0)\le \beta(h_1)$, meaning that $\beta:\{1,\dots, H\}\to\{1,\dots,H_N\}$ is nondecreasing. Furthermore,  by definition of $\varphi_n$ in (\ref{defhnell})
\begin{equation}\label{orph1}
\varphi_n\circ\beta(h_i)=\hhh_{\Delta(Q_nU_{j_i'})}'(s_k)=\min\Big\{1\le \ell\le H_{s_k}\;;\;\mathcal{I}\big(\Delta(Q_nU_{j_i'})\big)\cap J_\ell'(s_k)\ne\emptyset\Big\}
\end{equation}
while by definition of $\psi_k$ in (\ref{defhnellbis}) 
\begin{equation}\label{orph2}
\psi_k\circ\alpha(h_i)=\hhh_{\Delta(Q_{s_{k+1}}Q_{s_{k+2}}U_{j_i})}'(s_k)=\min\Big\{1\le \ell\le H_{s_k}\;;\;\mathcal{I}\big(\Delta(Q_{s_{k+1}}Q_{s_{k+2}}U_{j_i}))\big)\cap J_{\ell}'(s_{k})\ne\emptyset\Big\}\;;
\end{equation}
however, the equality $\Delta(Q_nU_{j_i'})=\Delta(Q_{s_{k+1}}Q_{s_{k+2}}U_{j_i})$ in (\ref{kkkkk}) allows to deduce from (\ref{orph1}) and (\ref{orph2}) that $\varphi_n\circ\beta(h_i)=\psi_k\circ\alpha(h_i)$: the theorem is proved.

\end{proof}

\begin{proof}[{\bf Proof of Theorem~A}]Let $1\le H\le d$, $\{1,\dots,H_n\}=E_1(n)\sqcup\dots\sqcup E_{H}(n)$ the partition and $V_1,\dots,V_H$ the limit projective vectors as  defined in Theorem~\ref{constant}.  By definition
\begin{equation}\label{DefJhn}
J_{h}(n):=\bigcup\nolimits_{\ell\in E_{h}(n)}J'_{\ell}(n).
\end{equation}

(i) : Fix $1\le h\le H$ and let $(j_1,j_2,\dots)$ be a sequence  in $\{1,\dots,d\}$ s.t. $j_n\in J_h(n)$, i.e. $j_n\in J_{\ell_n}'(n)$ for $\ell_n\in E_h(n)$, for any $n$ large enough; in particular $W_{\ell_n}(n)=V_h$. Hence, Lemma~\ref{distance} implies 
\begin{equation}\label{trap}
\lim_{n\to+\infty}\delta(P_nU_{j_n},V_h)=0
\end{equation}
and one concludes with part (iii) in Lemma~\ref{trapping} that 
$$
\lim_{n\to+\infty}\left\Vert {P_nU_{j_n}\over\Vert P_nU_{j_n}\Vert} -V_h\right\Vert=0.
$$
Part (i) of Theorem~A is proved.

(ii) : Let $s_{k+1}\le n< s_{k+2}$ (with $k\ge1$) and fix  $(j_0,j_1)\in  J_{h-1}'(n)\times J_h'(n)$, so that  $\Delta(Q_nU_{j_0})>\Delta(Q_nU_{j_1})\ne0$. Recall that condition (\ref{hypotplus}) means $Q_n\in\mathcal{H}_2(\Lambda)\cap\mathcal{H}_3(\lambda^k)$: hence, Lemma~\ref{majoration} allows to write for any $1\le i\le d$ (and $n\ge s_2$):
\begin{equation}\label{ii}
P_n(i,j_1)=
\sum\nolimits_jP_{s_k}(i,j)Q_n(j,j_1)\le 
\sum\nolimits_jP_{s_k}(i,j)\big(\lambda^k\Lambda\big)Q_n(j,j_0)=
\big(\lambda^k\Lambda\big)P_{n}(i,j_0).
\end{equation}
Given $(\ell_1,\ell_2,\dots)$ s.t. $1< \ell_n\le H_n$, for any $n\ge1$ and $(j_1,j_2,\dots)$, $(j_1',j_2',\dots)$ two index sequences, it follows from (\ref{ii})~that
\begin{equation*}\label{proofii}
\Big(\forall n\ge s_2,\;(j_n,j_n')\in  J_{\ell_n-1}'(n)\times  J_{\ell_n}'(n)\Big)\;\Longrightarrow\;\lim_{n\to+\infty}{\Vert P_nU_{j_n'}\Vert/\Vert P_nU_{j_n}\Vert}=0
\end{equation*}
and thus, for $1< h\le H$, the definition of $J_h(n)$ in (\ref{DefJhn}) and the  property (\ref{Vii}) of $E_n(h)$ gives 
\begin{equation*}
\Big(\forall n\ge s_2,\;(j_n,j_n')\in  J_{h-1}(n)\times  J_{h}(n)\Big)\;\Longrightarrow\;\lim_{n\to+\infty}{\Vert P_nU_{j_n'}\Vert/\Vert P_nU_{j_n}\Vert}=0.
\end{equation*}
Part (ii) of Theorem~A is proved.

(iii) : Let $X\in\mathcal{S}_d$ such that $P_n(X)\ne0$ for any $n\in\mathbb N$, and let $j_n\in J'_{\hhh_X'(n)}(n)$, where   $\hhh_X'(n):=\min\{1\le \ell\le H_n\;;\;\mathcal{I}(X)\cap J'_{\ell}(n)\ne\emptyset\}$: by  inequality (\ref{Xandc}) in Lemma~\ref{rexpk} one has  
\begin{equation}\label{THAEQ1}
\delta(P_nX,P_nU_{j_n})\le (C\Lambda_X)r^k
\end{equation}
(for $k=\kk(n)$ s.t. $s_{k+1}\le n< s_{k+2}$).  According to (\ref{DefJhn}) and the property (\ref{Vii}) of $E_n(h)$,
$
j_n\in J'_{\hhh_X'(n)}(n)\subset J_{\hhh_X(n)}(n)$ and 
$W_{\hhh_X'(n)}(n)=V_{\hhh_X(n)},
$
where $\hhh_X(n):=\min\{1\le h\le H\;;\;\mathcal{I}(X)\cap J_{h}(n)\ne\emptyset\}$ (as defined in Theorem~A).
However, we know (Lemma \ref{distance})  that 
$$
\varepsilon_n:=\sup\bigcup_{h=1}^{H_n}\Big\{\delta\big(P_nU_j,W_h(n)\big)\;;\;j\in J_h'(n)\Big\}\to0
$$
as $n\to+\infty$ and because
$$
\delta\big(P_nU_{j_n},V_{\hhh_X(n)}\big)=\delta\big(P_nU_{j_n},W_{\hhh_X'(n)}(n)\big)\le\varepsilon_n,
$$
it follows from  (\ref{THAEQ1}) that
\begin{align*}
\delta\left({P_nX\over \Vert P_nX\Vert},V_{\hhh_X(n)}\right)
&=\delta\left(P_nX,V_{\hhh_X(n)}\right)\\
&\le \delta\left(P_nX,P_nU_{j_n}\right)+
\delta\left(P_nU_{j_n},V_{\hhh_X(n)}\right)\le (C\Lambda_X)r^{\kk(n)}+\varepsilon_n.
\end{align*}
For $n$ large enough, $\delta\big({P_nX/\Vert P_nX\Vert},V_{\hhh_X(n)}\big)<+\infty$, meaning that $\Delta(P_nX/\Vert P_nX\Vert)=\Delta(V_{\hhh_X(n)})$, and with (\ref{doublebound}) in Proposition~\ref{face} we conclude
$$
\left\Vert{P_nX\over \Vert P_nX\Vert}-V_{\hhh_X(n)}\right\Vert\le d\delta\left({P_nX\over \Vert P_nX\Vert},V_{\hhh_X(n)}\right)\le\ d\Big(\varepsilon_n+C\cdot r^{\kk(n)}\Big)\Lambda_X. 
$$
Part (iii) of Theorem~A is proved.

\end{proof}

\section{\bf Heuristics for Theorem~A -- basic examples and applications}\label{heuristic}

\subsection{Example 1: Product of block-triangular matrices}This is a natural and relatively simple example of application of Theorem~A: one suppose that each $A_n$ is lower block-triangular (resp. upper block-triangular) and that each block has only positive entries and size independent of $n$. 

\begin{corollary}\label{tbbcase}(i) Suppose that
$$
A_n=\begin{pmatrix}B_n(1,1)&0&\cdots&0\\B_n(2,1)&B_n(2,2)&\cdots&0\\\vdots&\vdots&\ddots&\vdots\\B_n(\delta,1)&B_n(\delta,2)&\cdots&B_n(\delta,\delta)\end{pmatrix}
$$
where the entries of each submatrix $B_n(i,j)$ are positive. We suppose that the size $d_i\times d_j$ of the matrix $B_n(i,j)$ is independent of $n$, as well as the real $\Lambda\ge1$ such that $A_n\in\mathcal H_2(\Lambda)$. Denoting by $U_1^{(k)},\dots,U_k^{(k)}$ the canonical $k$-dimensional column vectors we suppose that
\begin{equation}\label{thelowercase}
\forall i,j,\quad\sum_{n=1}^{+\infty}\frac{\Vert B_1(i+1,i+1)\cdots B_n(i+1,i+1)\Vert}{\big\Vert B_1(i,i)\cdots B_n(i,i)U_j^{(d_i)}\big\Vert}<+\infty.
\end{equation}
Then the sequence $(A_1,A_2,\dots)$ satisfies condition ${\bf (C)}$ and Theorem~A applies.

(ii) This conclusion remains true if we replace the hypothesis "$A_n$ lower block-triangular" by "$A_n$ upper block-triangular" and (\ref{thelowercase}) by
$$
\forall i,j,\quad\sum_{n=1}^{+\infty}\frac{\Vert B_1(i-1,i-1)\cdots B_n(i-1,i-1)\Vert}{\big\Vert B_1(i,i)\cdots B_n(i,i)U_j^{(d_i)}\big\Vert}<+\infty.
$$
\end{corollary}

\begin{proof}[{\bf Proof}](i) The norm we use in this proof is $\Vert M\Vert=\max_jL^{(k)}MU^{(\ell)}_j$ for any $k\times\ell$ matrix $M$, where $L^{(k)}=\begin{pmatrix}1&\cdots&1\end{pmatrix}$ is the $k$-dimensional row vector with entries $1$. Let
\begin{equation}\label{summableseries}
\varepsilon_n=\max_{1\le i<\delta\atop1\le j\le d_i}\frac{\Vert B_1(i+1,i+1)\cdots B_n(i+1,i+1)\Vert}{\big\Vert B_1(i,i)\cdots B_n(i,i)U_j^{(d_i)}\big\Vert}\quad(n\ge N),
\end{equation}
and $\displaystyle S=1+\Lambda\sum_{n=0}^{+\infty}\varepsilon_n<+\infty$ (where $\varepsilon_0=1$). We put $B_n(i)=B_n(i,i)$ and we define $A_n(i),C_n(i)$ by
\begin{equation}\label{bigmatrix}
A_n(i)=\begin{pmatrix}B_n(i,i)&\cdots&0\\\vdots&\ddots&\vdots\\B_n(\delta,i)&\cdots&B_n(\delta,\delta)\end{pmatrix}=\begin{pmatrix}B_n(i)&0\\C_n(i)&A_n(i+1)\end{pmatrix}.
\end{equation}
Let us prove by descending induction on $1<i\le\delta$ that -- for any $n\in\mathbb N$ and $1\le j\le d_{i-1}$
\begin{equation}\label{toprove}
A_1(i)\cdots A_n(i)\in\mathcal H_2(S^{\delta-i}\Lambda)\quad\hbox{and}\atop\Vert A_1(i)\cdots A_n(i)\Vert\le\varepsilon_nS^{\delta-i}\big\Vert B_1(i-1,i-1)\cdots B_n(i-1,i-1)U_j^{(d_{i-1})}\big\Vert.\end{equation}

Let $i=\delta$, one has $A_1(\delta)\cdots A_n(\delta)=B_1(\delta,\delta)\cdots B_n(\delta,\delta)$ and this  product of positive matrices belongs to $\mathcal H_2(\Lambda)$ by Lemma \ref{hypotheses} (ii). Using the definition of $\varepsilon_n$ in (\ref{summableseries}), the induction hypotheses are satisfied at the rank $i=\delta$.

We suppose now that the induction hypotheses are satisfied at some rank $i+1\le\delta$ and we prove it at the rank $i$. From the second equality in~(\ref{bigmatrix}),
\begin{equation}\label{bigtriangle}
\displaystyle\prod_{k=1}^nA_k(i)=\begin{pmatrix}\displaystyle\prod_{k=1}^nB_k(i)&0\\\displaystyle\sum_{k=1}^n\left(\prod_{\ell=1}^{k-1}A_\ell(i+1)C_k(i)\prod_{\ell=k+1}^nB_\ell(i)\right)&\displaystyle\prod_{k=1}^nA_k(i+1)\end{pmatrix}
\end{equation}
where the notation $\displaystyle\prod_{k=a}^bM_k$ means the identity matrix if $b<a$. Since we have supposed that the induction hypotheses are true at the rank $i+1$, the columns of $\displaystyle\prod_{k=1}^nA_k(i)$ whose index is larger than $d_i$ belong to $\mathcal H_2(S^{\delta-(i+1)})$ and their norms are bounded by $\varepsilon_nS^{\delta-(i+1)}\big\Vert B_1(i,i)\cdots B_n(i,i)U_j^{(d_i)}\big\Vert$ ($1\le j\le d_i$). Consequently they belong to $\mathcal H_2(S^{\delta-i})$. Moreover their norms are bounded by $\varepsilon_nS^{\delta-i}\big\Vert B_1(i-1,i-1)\cdots B_n(i-1,i-1)U_{j'}^{(d_{i-1})}\big\Vert$ because $\frac{\big\Vert B_1(i,i)\cdots B_n(i,i)U_j^{(d_i)}\big\Vert}{\big\Vert B_1(i-1,i-1)\cdots B_n(i-1,i-1)U_{j'}^{(d_{i-1})}\big\Vert}\le\varepsilon_n\le S$ ($1\le j'\le d_{i-1}$). So it remains to look at the columns of $A_1(i)\cdots A_n(i)$ with index $j\le d_i$.

Let $D_i=d_i+\dots+d_\delta$. Using the inequality $L^{(D_{i+1})}C_k(i)\le\Lambda L^{(d_i)}B_k(i)$ (consequence of $A_k\in\mathcal H_2(\Lambda)$) we obtain for any $j,j'\in\{1,\dots,d_i\}$
\begin{equation}\label{with}
\begin{array}{rcl}\displaystyle\Big\Vert\prod_{k=1}^nA_k(i)U_j^{(D_i)}\Big\Vert&=&\displaystyle\Big\Vert\prod_{k=1}^nB_k(i)U_j^{(d_i)}\Big\Vert+\sum_{k=1}^n\Big\Vert\prod_{\ell=1}^{k-1}A_\ell(i+1)C_k(i)\prod_{\ell=k+1}^nB_\ell(i)U_j^{(d_i)}\Big\Vert\\
&\le&\displaystyle\Big\Vert\prod_{k=1}^nB_k(i)U_j^{(d_i)}\Big\Vert+\sum_{k=1}^n\Big\Vert\prod_{\ell=1}^{k-1}A_\ell(i+1)\Big\Vert L^{(D_{i+1})}C_k(i)\prod_{\ell=k+1}^nB_\ell(i)U_j^{(d_i)}\\&\le&\displaystyle\Big\Vert\prod_{k=1}^nB_k(i)U_j^{(d_i)}\Big\Vert\\&&\displaystyle+\sum_{k=1}^n\varepsilon_{k-1}S^{\delta-(i+1)}\Big\Vert\prod_{\ell=1}^{k-1}B_\ell(i)U_{j'}^{(d_i)}\Big\Vert\Lambda L^{(d_i)}\prod_{\ell=k}^nB_\ell(i)U_j^{(d_i)}.\end{array}
\end{equation}
Clearly, for any $d_i$-dimensional row vector $L$ and any $d_i$-dimensional column vector $R$,
$$
\min_{1\le j'\le d_i}\big(LU_{j'}^{(d_{i})}L^{(d_i)}R\big)\le LR.
$$
Applying this to $\displaystyle L=L^{(d_i)}\prod_{\ell=1}^{k-1}B_\ell(i)$ and $\displaystyle R=\prod_{\ell=k}^nB_\ell(i)U_j^{(d_i)}$ we obtain
\begin{equation}\label{lastbutonestep}
\begin{array}{rcl}\displaystyle\Big\Vert\prod_{k=1}^nA_k(i)U_j^{(D_i)}\Big\Vert&\le&\displaystyle\Big\Vert\prod_{k=1}^nB_k(i)U_j^{(d_i)}\Big\Vert+\sum_{k=1}^n\Lambda\varepsilon_{k-1}S^{\delta-(i+1)}L^{(d_i)}\prod_{\ell=1}^nB_\ell(i)U_j^{(d_i)}\\&\le&\displaystyle S^{\delta-i}\Big\Vert\prod_{\ell=1}^nB_\ell(i)U_j^{(d_i)}\Big\Vert.\end{array}
\end{equation}
Using the the definition of $\varepsilon_n$ in (\ref{summableseries}), we deduce that the second condition of (\ref{toprove}) is fulfilled.

In order to check the first condition, we consider a entry in the $j^{\rm th}$ column of $A_1(i)\cdots A_n(i)$. From (\ref{bigtriangle}), either it is a entry of the $j^{\rm th}$ column of $B_1(i)\cdots B_n(i)$, or it is at least equal to a entry of the $j^{\rm th}$ column of $C_1(i)B_2(i)\cdots B_n(i)$. From the second equality in (\ref{bigmatrix}) one has $\begin{pmatrix}B_1(i)\\C_1(i)\end{pmatrix}\in\mathcal H_2(\Lambda)$ and, from Lemma \ref{hypotheses} (ii), the product matrix $\begin{pmatrix}B_1(i)\\C_1(i)\end{pmatrix}B_2(i)\cdots B_n(i)$ also belongs to $\mathcal H_2(\Lambda)$. So, the entry we have considered is at least equal to $\big\Vert B_1(i)\cdots B_n(i)U_j^{(d_i)}\big\Vert/\Lambda$. Using (\ref{lastbutonestep}) it is at least equal to $\big\Vert A_1(i)\cdots A_n(i)U_j^{(D_i)}\big\Vert/(S^{\delta-i}\Lambda)$, proving that the first condition of (\ref{toprove}) is fulfilled.

So( \ref{toprove}) is true for any $n\in\mathbb N$, and we deduce that $A_1\cdots A_n\in\mathcal H_2(S^{\delta-1}\Lambda)\cap\mathcal H_3(\varepsilon_nS^{\delta-2}\Lambda)$. For any $n$ large enough, 
$$
A_1\cdots A_n\in\mathcal H_2(S^{\delta-1}\Lambda)\cap\mathcal H_3(1/2).
$$
Notice that for any nonnegative integer $N$ the hypotheses of Corollary \ref{tbbcase}, and in particular the condition (\ref{thelowercase}), are satisfied by $(A_{N+n})_{n\in\mathbb N}$ because, denoting by $m$ (resp. $M$) the minimum (resp. the maximum) of the norms of $B_1(i,i)\cdots B_N(i,i)U_j^{(d_i)}$ for any $i,j$, we have
$$
m\big\Vert B_{N+1}(i,i)\cdots B_n(i,i)U_j^{(d_i)}\big\Vert\le\big\Vert B_1(i,i)\cdots B_n(i,i)U_j^{(d_i)}\big\Vert\le M\big\Vert B_{N+1}(i,i)\cdots B_nN(i,i)U_j^{(d_i)}\big\Vert.
$$
So we can define by induction the sequence $(s_0,s_1,\dots)$ such that $(A_n)_{n\in\mathbb N}$ satisfies condition ${\bf (C)}$ w.r.t. $\lambda=1/2$ and $\Lambda=S^{\delta-1}\Lambda$: one put $s_0=s_1=0$ and for any $k$, $s_{k+1}$ is the smallest integer such that
$$
\forall n\ge s_{k+1},\ A_{s_k+1}\cdots A_n\in\mathcal H_3(1/2).
$$

(ii) Let $\Delta=\begin{pmatrix}0&\cdots&0&1\\0&\cdots&1&0\\\vdots&\ddots&\vdots&\vdots\\1&\cdots&0&0\end{pmatrix}$. If the $A_n$ are upper block-triangular, the $\Delta A_n\Delta$ are lower block-triangular and the order of the diagonal blocs is reversed.

\end{proof}

\subsection{Example 2}Let us look at the case of the $2\times2$ upper triangular matrices:
$$
P_n:=
\begin{pmatrix}a_1&b_1\\0&d_1\end{pmatrix}\cdots
\begin{pmatrix}a_n&b_n\\0&d_n\end{pmatrix}
$$
with $a_nd_n>0$ and $b_n\ge0$. The sequence of normalized columns of $P_n$ converges in $\mathcal{S}_2$, with more precisely
$$
\lim_{n\to+\infty}{P_nU_j\over \Vert P_nU_j\Vert}=
\begin{cases}
\begin{pmatrix}\frac s{s+1}\\\frac1{s+1}\end{pmatrix}&\text{with $s=\displaystyle \sum_{n=1}^{+\infty}\frac{a_1\dots a_{n-1}b_n}{d_1\dots d_{n-1}d_n}$, if $j=2$ and $s<+\infty$}\\\\
\begin{pmatrix}1\\0\end{pmatrix}&\text{if either $j=1$ or $s=+\infty$.}
\end{cases}
$$
If $s<+\infty$ and $\displaystyle\lim_{n\to+\infty}\frac1n\log(d_1\dots d_n)$ exists (resp.  if $s<+\infty$ and $\displaystyle\lim_{n\to+\infty}\frac1n\log(a_1\dots a_n)$ exists), it is the first (resp. the second) Lyapunov exponent of the sequence of product matrices $(P_1,P_2,\dots)$.

\subsection{Example 3}Corollary A does not hold for the following products of triangular matrices. Let $n_k=1+\cdots+k$ and define $P_n=A_1\cdots A_n$, where

$$
A_{n_k}=\begin{pmatrix}1&1&1\\0&0&0\\0&0&2\end{pmatrix}
\quad\text{while}\quad n\ne n_k\Longrightarrow A_n=\begin{pmatrix}1&3&1\\0&4&0\\0&0&1\end{pmatrix}.
$$
Then, one gets 
$$
P_{n_k}=
\begin{pmatrix}1&1&3\cdot2^k-2k-3\\0&0&0\\0&0&2^k\end{pmatrix}
\quad \hbox{and}\quad
P_{n_k+k}=\begin{pmatrix}1&2\cdot4^k-1&3\cdot2^k-k-3\\0&0&0\\0&0&2^k\end{pmatrix}
$$
and thus,  provided that $\Delta(V)\ge(0,1,1)$, 
$$
\lim_{k\to+\infty}
\frac{P_{n_k}V}{\Vert P_{n_k}V\Vert}=\begin{pmatrix}\frac34\\0\\\frac14\end{pmatrix}
\quad\text{while}\quad
\lim_{k\to+\infty}\frac{P_{n_k+k}V}{\Vert P_{n_k+k}V\Vert}=\begin{pmatrix}1\\0\\0\end{pmatrix}.
$$

\subsection{Example 4: positive matrices}Suppose that $A_1,A_2,\dots$ is an infinite sequence of $d\times d$ matrices each ones having positive entries. We note $I_d:=\{1,\dots,d\}$, so that (see Definition~\ref{DefSI}) $\mathcal{S}_{I_d}$ is the subset of $\mathcal{S}_d$  whose element are the probability vectors with positive entries. Fix $X\in\mathcal{S}_{d}$ arbitrarily (i.e with possibly zero entries) and suppose that $\gamma:=\sup_i\{\tau(A_i)\}<1$; then, as an immediate consequence of Lemma \ref{deltatau} the probability vectors
$P_1X/\Vert P_1X\Vert$, $P_2X/\Vert P_2X\Vert,\dots$ form a 
Cauchy sequence in the (complete, non compact) metric space $(\mathcal{S}_{I_d},\delta(\cdot,\cdot))$ and thus has a limit $X_*\in\mathcal{S}_{I_d}$. More precisely, one has the following proposition.
\begin{proposition}[Folklore]\label{Folk}Let $\mathcal{A}=(A_1,A_2,\dots)$ be a sequence of $d\times d$ matrices having positive entries such that $\gamma:=\sup_i\{\tau(A_i)\}<1$: then there exist  $X_*\in\mathcal{S}_{I_d}$ and $C>0$ such that 
$$
\forall X\in\mathcal{S}_d,\quad \left\Vert {P_nX\over \Vert P_nX\Vert}-X_*\right\Vert\le C\gamma^n.
$$
\end{proposition}
In the framework of Proposition~\ref{Folk}, one deduces the matrix product $P_n$ to be asymptotically closed to a rank one matrix, in the sense that, there exists $d$ sequences $n\mapsto\alpha_i(n)$ ($1\le i\le d$) each ones being made of positive real numbers and such that 
$$
P_n\approx\Big(\alpha_1(n)X_*\ \cdots\ \alpha_d(n)X_*\Big).
$$
Here, the condition ${\bf (C)}$ is trivially fulfilled and Theorem~A is applicable: according to Proposition~\ref{Folk}  --~and with the notations of Theorem~A~-- one has $H=1$ with $V_1=X_*$ and $J_{1}(n)=\{1,\dots,d\}$, for any $n\ge1$. 
The sequences $n\mapsto\frac{\alpha_i(n)}{\sum_j\alpha_j(n)}$ are not necessarily convergent. Consider for instance $A_i=\beta_iB_i$ where
$$
B_{2i-1}=\begin{pmatrix}1/2&1/2\\1/2&1/2\end{pmatrix}
\quad\text{and}\quad
B_{2i}=\begin{pmatrix}1/3&2/3\\1/3&2/3\end{pmatrix}
$$ 
so that $P_n=(\beta_1\cdots\beta_n)B_n$.
Suppose now that $A_i$ takes values in a finite set, say $\{A(\0),\dots,A(\aaa)\}$: the condition that $\gamma:=\sup_n\{\tau(A_i)\}<1$ is automatically satisfied and Proposition~\ref{Folk} holds. 

\subsection{Example 5}\label{toplyapmap} 
\begin{figure}[H]
\begin{center}
       \includegraphics[scale=0.7]{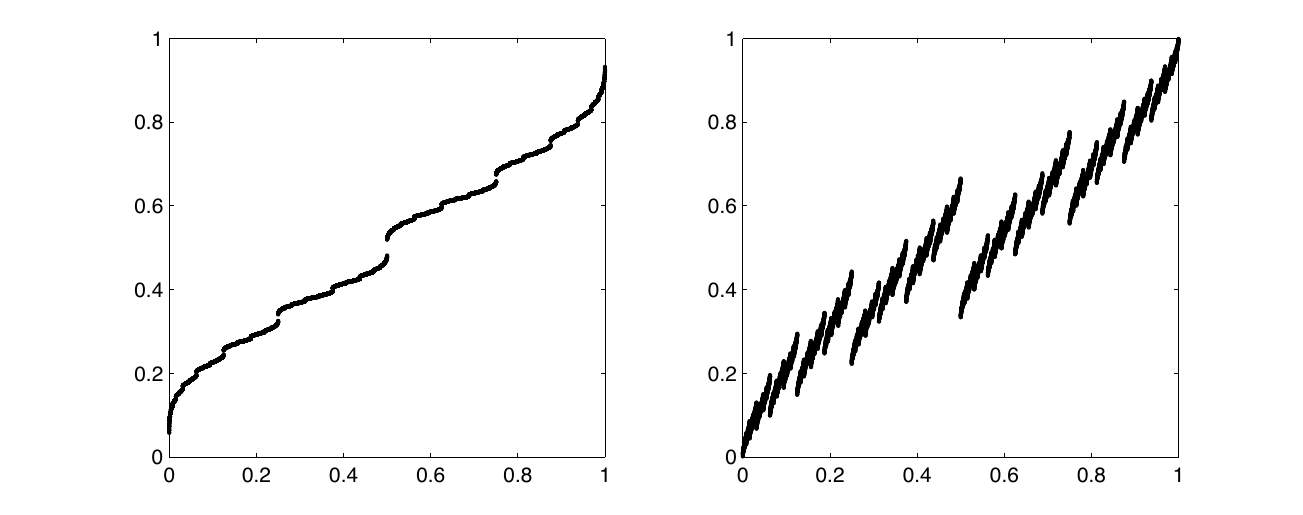}
\caption{\label{complicated}\leftskip=5mm \rightskip=5mm\small  \it Representation of the map $[0\,;1[\ni x\mapsto\mathtt{p}(x_1x_2\cdots)$, where $0\cdot x_1x_2\cdots$ is the binary expansion of $x$ and $\mathtt{p}:\{\0,\1\}^\NNN\to[0\,;1]$ is associated with $A(\0)$ and $A(\1)$ given in (\ref{contfrac}) and (\ref{UIBCM}) respectively; in the left box plot, one recognizes the inverse of the Minkowski question mark function (also known as the Conway's function), while in the right box plot one finds a limit Rademacher function.} 
\end{center}
\end{figure}

We shall illustrate the possible {\it underlying fractal nature} of the {\it top Lyapunov direction} as introduced in Corollary~A. Let $A(\0),\dots, A(\aaa)$ be fixed $d\times d$ matrices with nonnegative entries and for $\omega=\omega_1\omega_2\cdots\in\Omega:=\{\0,\dots,\aaa\}^\NNN$, consider the sequence $\mathcal{A}(\omega)=(A(\omega_1),A(\omega_2),\dots)$. Provided it makes sense, the top Lyapunov direction map is  $\VV:\Omega\to\mathcal{S}_d$~s.t.
$$
\VV(\omega):=\lim_{n\to+\infty}{P_n(\omega)U\over \Vert P_n(\omega)U\Vert}
$$
(here we recall that $P_n(\omega)=A(\omega_1)\cdots A(\omega_n)=:A(\omega_1\cdots\omega_n)$ and $U=U_1+\cdots+U_d$). First consider that each $A(i)$ is a $2\times 2$ matrix of the form $\begin{pmatrix}a(i)&b(i)\\1-a(i)&1-b(i)\end{pmatrix}$ with the additional condition that $\sup_i\vert \det(A(i))\vert<1$. Then (see \cite{MNR99,MNR99a,MN02}), one gets
\begin{equation}\label{mukerpi}
\VV(\omega)=\begin{pmatrix}\mathtt{p}(\omega)\\1-\mathtt{p}(\omega)\end{pmatrix}\quad\text{where}\quad\mathtt{p}(\omega)=\sum_{n=1}^{+\infty} b(\omega_n)\det A(\omega_1\cdots\omega_{n-1}).
\end{equation}
To fix ideas let $1\le \beta\le 2$ and consider that $\aaa=1$ with the two matrices
\begin{equation}\label{UIBCM}
A(\0)=\begin{pmatrix}
1/\beta&0\\
1-1/\beta&1\\
\end{pmatrix}
\quad\text{and}\quad
A(\1)=\begin{pmatrix}
1&1-1/\beta\\
0&1/\beta\\
\end{pmatrix},
\end{equation}
so that (\ref{mukerpi}) gives
$$
\mathtt{p}(\omega)=(\beta-1)\sum_{n=1}^{+\infty}{\omega_n\over \beta^{n}}
$$ 
(see Figure~\ref{complicated}-right).

An other very classical case is related to the expansion of numbers by means of continued fraction: to see this, take  $\aaa=1$ and the two matrices: 
\begin{equation}\label{contfrac}
A_\0=
\begin{pmatrix}
1&0\\1&1
\end{pmatrix}
\quad\text{and}\quad
A_\1=\begin{pmatrix}
1&1\\0&1
\end{pmatrix}.
\end{equation}
Also in this case, the top Lyapunov direction map  about a sequence $\omega\in\Omega=\{\0,\1\}^\NNN$, may be sketched by  direct  computation. Consider  $\omega=\0^{a_0}\1^{a_1}\0^{a_2}\cdots$, where $a_0,a_1,\cdots$ are integers (with $a_0\ge0$ and $a_i>0$, for $i\ge1$); for any $n\ge0$ and $\varepsilon\in\{\0,\1\}$, a simple induction gives  
\begin{equation}\label{stroumpf}
A(\1^{a_0}\0^{a_1}\1^{a_2}\cdots\varepsilon^{a_{n}})=
\begin{pmatrix}
a_0&1\\1&0
\end{pmatrix}
\begin{pmatrix}
a_{1}&1\\1&0
\end{pmatrix}
\cdots
\begin{pmatrix}
a_{n}&1\\1&0
\end{pmatrix}
\begin{pmatrix}
0&1\\1&0
\end{pmatrix}^\varepsilon=
\begin{pmatrix}
q_{n}&q_{n-1}\\p_{n}&p_{n-1}
\end{pmatrix}
\begin{pmatrix}
0&1\\1&0
\end{pmatrix}^\varepsilon
\end{equation}
where by convention $(q_{-1},p_{-1})=(1,0)$, while for $n\ge0$ 
\begin{equation}\label{defcontfrac} 
{p_n\over q_n}={1\over a_0+\displaystyle{{1\over \ddots\raise-8pt\hbox{$+\displaystyle{1\over a_n}$}}}}=:\[a_0,\dots,a_n\]
\end{equation}
converges toward an irrational real number $x=\[a_0,a_1,\dots\]\in[0\,;1]$; we use (\ref{stroumpf}),  and put $\theta_n(x):=q_{n-1}/q_n$, so that by the approximation $p_n\approx xq_n$, one gets 
$$
{A(\1^{a_0}\0^{a_1}\1^{a_2}\cdots\varepsilon^{a_{n}})U\over \Vert A(\1^{a_0}\0^{a_1}\1^{a_2}\cdots\varepsilon^{a_{n}})U\Vert}
\approx
{1\over (1+x)(1+\theta_n(x))}\begin{pmatrix}
1&\theta_{n}(x)\\x&x\theta_{n}(x)
\end{pmatrix}
\begin{pmatrix}
1\\1
\end{pmatrix}=
\begin{pmatrix}
H(x)\\1-H(x)
\end{pmatrix}
$$
where $H(x)=1/(1+x)$;
then, it is simple to deduce that (see Figure~\ref{complicated}-left for illustration)
\begin{eqnarray}\label{approx}
\VV(\omega)=\begin{pmatrix}\mathtt{p}(\omega)\\1-\mathtt{p}(\omega)\end{pmatrix}
\quad\text{where}\quad\mathtt{p}(\omega)=\[1,a_0,a_1,\dots\]
\end{eqnarray}

\subsection{Example 6} We consider that $P_n$ is the $n$-step product whose formal limit (as $n\to+\infty$) is the infinite product 
$RSTST^2ST^3S\cdots ST^{k}S\cdots$,
where
$$
R:=\begin{pmatrix}2&1&1&1\\1&2&1&1\\0&0&0&0\\0&0&0&0\end{pmatrix},\quad S:=\begin{pmatrix}11&0&0&0\\7&4&0&0\\7&4&0&0\\7&2&1&1\end{pmatrix},\quad T:=\begin{pmatrix}11&0&0&0\\7&4&0&0\\7&2&2&0\\7&2&1&1\end{pmatrix}.
$$
Here, condition ${\bf (C)}$ may be checked directly and the (non necessarily injective) map $h\mapsto V_h$ (with $1\le h\le H$ and $V_{h-1}\ne V_h$) given by Theorem~A is defined by
$$
V_1=\begin{pmatrix}1/2\\1/2\\0\\0\end{pmatrix},\quad
V_2=\begin{pmatrix}3/7\\4/7\\0\\0\end{pmatrix},\quad
V_3=\begin{pmatrix}1/2\\1/2\\0\\0\end{pmatrix}\quad\hbox{and}\quad H=3.
$$ 
It is rather simple to illustrate how the number of different (exponential) growth rates of $\Vert P_nU_j\Vert$ may depend on $n$ --~while $H=3$ remains constant~-- together with the fact  that (for $2\le h\le H$) any column of $P_n$ that converges in direction to $V_h$ is negligible with respect to (at least) one column that converges in direction to $V_{h-1}$.

\begin{figure}[H]
\includegraphics[scale=0.8]{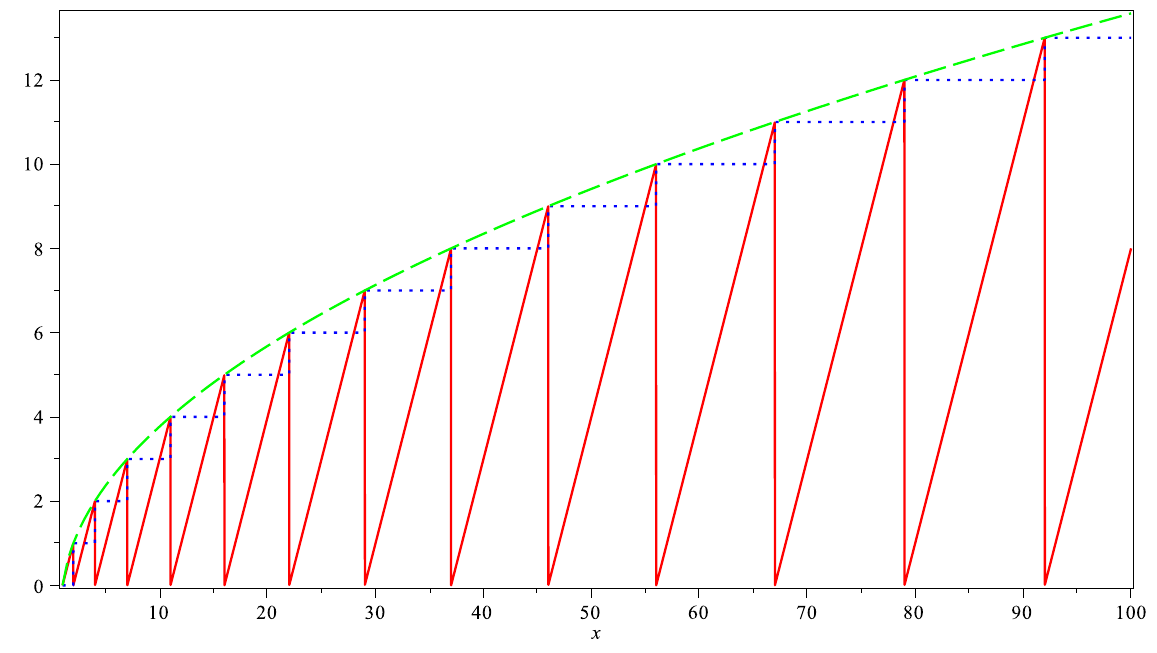}
\caption{\label{scie}\small\it The solid red graph represent the exponential growth rate of $C_3^{(n)}$:
$$
x\mapsto\varphi(x)=
x-\left(1+{1\over 2}\left\lfloor{1\over 2}\sqrt{8x-7}-{1\over 2}\right\rfloor\left\lfloor{1\over 2}\sqrt{8x-7}+{1\over 2}\right\rfloor\right),
$$
while the dashed green and the dotted blue ones represent respectively 
$$
x\mapsto{1\over 2}\sqrt{8x-7}-{1\over 2}\quad\text{and}\quad x\mapsto\left\lfloor{1\over 2}\sqrt{8x-7}-{1\over 2}\right\rfloor.
$$}
\end{figure}

Consider the increasing sequence of integers $n_0=1$, $n_1=n_0+1$, $n_2=n_1+2,\dots$, $n_{k}=n_{k-1}+k,\dots$, so that 
$$
P_{n_0}=R,\quad P_{n_1}=R(S),\quad P_{n_2}=RS(TS),\quad P_{n_3}=RSTS(TTS),\dots
$$
Any integer $n\ge1$, may be written $n=n_k+r$ for some $k\ge0$ and $0\le r\le k$, so that

\begin{equation}\label{Pn}
P_n=\left(
\begin{array}{c|c|c|c}
5\cdot11^{n-1}-3\cdot4^{n-1}&3\cdot4^{n-1}-2^{r+1}&2^{r+1}-1&\ 1\\
5\cdot11^{n-1}-4\cdot4^{n-1}&4\cdot4^{n-1}-2^{r+1}&2^{r+1}-1&\ 1\\
0&0&0&\ 0\\
0&0&0&\ 0
\end{array}
\right).
\end{equation}
For any $n\ge1$, the sets of indices $J_h(n)$ ($1\le h\le H=3$) involved in Theorem~A together with the estimates of $\Vert P_nU_j\Vert$ for $j\in J_h(n)$ are  respectively
$J_{1}(n)=\{1\}$ with $\Vert P_nU_1\Vert \approx 11^n$,
$J_{2}(n)=\{2\}$ with $\Vert P_nU_2\Vert \approx 4^n$
and finally, $J_{3}(n)=\{3,4\}$ with $\Vert P_nU_4\Vert \approx 1$ and 
$\Vert P_nU_3\Vert \approx 2^{\varphi(n)}$. The function $n\mapsto\varphi(n)$ may be obtained by a straightforward computation, which gives
$$
\varphi(n)=r=n-\left(1+{1\over 2}\left\lfloor{1\over 2}\sqrt{8n-7}-{1\over 2}\right\rfloor\left\lfloor{1\over 2}\sqrt{8n-7}+{1\over 2}\right\rfloor\right).
$$
In particular, for $n=n_k+k=1+k(k+1)/2+k$, one gets $\Vert P_nU_3^{(n)}\Vert \approx 2^{\sqrt{2n}}$ (see Figure~\ref{scie} for a representation of $x\mapsto\varphi(x)$).

\section{\bf Rank one property of normalized matrix products}\label{limpoints}
Let $A$ be a matrix in the space $\mathcal{M}_d(\CCC)$ of the $d\times d$ matrices with complex entries. We denote by $A^\star$ the adjoint matrix of $A$ and we write its {\it singular value decomposition} 
$$
A=S\begin{pmatrix}\sigma_1&\dots&0\\\vdots&\ddots&\vdots\\0&\dots&\sigma_d\end{pmatrix}T^\star\;;
$$
here $S$ and $T$ are unitary matrices and 
$\sigma_1\ge\cdots\ge \sigma_d$
is the ordered lists of the {\it singular values} of $A$. 
Recall that by definition $A\mapsto\Vert A\Vert_2:=\sigma_1$ defines the so-called {\it spectral norm} of $A$. 
In this paragraph, we  consider that $\mathcal{P}=(P_1,P_2,\dots)$ is a sequence in  $\mathcal{M}_d(\CCC)$ and  we write 
\begin{equation}\label{SVD}
P_n=S_n\begin{pmatrix}e^{n\chi_1(n)}&\dots&0\\\vdots&\ddots&\vdots\\0&\dots&e^{n\chi_d(n)}\end{pmatrix}T_n^\star
\end{equation}
the singular value decomposition of $P_n$. Here the $i$-th singular value of $P_n$ is denoted by $e^{n\chi_i(n)}$ and $\chi_1(n)\ge\cdots\ge\chi_d(n)$ may be thought as the ordered list of the $n$-step Lyapunov exponents. 
\begin{theorem}\label{r1}Let $\mathcal{P}=(P_1,P_2,\dots)$ a sequence in $\mathcal M_d(\mathbb C)$; then there exists a sequence $\mathcal{B}=(B_1,B_2,\dots)$ of rank $1$ matrices s.t. $\left\Vert{P_n}/{\Vert P_n\Vert}-B_n\right\Vert\to 0$ if and only if ${e^{n\chi_2(n)}/e^{n\chi_1(n)}}\to0$.
\end{theorem}

\begin{proof}[\bf Proof]The norms on the space $\mathcal M_d(\mathbb C)$ being equivalent, and because 
$$
\frac{P_n}{\Vert P_n\Vert}-B_n=\frac{\Vert P_n\Vert_2}{\Vert P_n\Vert}\ \left(\frac{P_n}{\Vert P_n\Vert_2}-\frac{\Vert P_n\Vert}{\Vert P_n\Vert_2}B_n\right),
$$
we work with $\Vert\cdot\Vert_2$ in place of $\Vert\cdot\Vert$.
According to the singular value decomposition  of $P_n$ in (\ref{SVD}), 
\begin{equation}\label{svd}
\frac{P_n}{\Vert P_n\Vert_2}=\frac{P_n}{e^{n\chi_1(n)}}=S_n(U_1U_1^\star)T_n^\star+\sum_{i=2}^d{e^{n\chi_i(n)}\over e^{n\chi_1(n)}}S_n(U_iU_i^\star)T_n^\star.
\end{equation}
Therefore, the condition $e^{n(\chi_2(n)-\chi_1(n))}\to0$ for $i=2,\dots d$ ensures that $\left\Vert{P_n}/{\Vert P_n\Vert_2}-B_n\right\Vert_2\to 0$ holds with $B_n:=S_n(U_1U_1^\star)T_n^\star$.
For the converse implication,
assume the existence of rank $1$ matrices $B_1,B_2,\dots$ such that $\left\Vert{P_n}/{\Vert P_n\Vert_2}-B_n\right\Vert_2\to 0$, and assume --~for a contradiction~-- the existence of a sequence of integers $1\le n_1<n_2<\cdots$ for which $e^{n_k\chi_2(n_k)}/e^{n_k\chi_1(n_k)}\to\alpha_2>0$. By a compactness argument,
it is possible to choose the $n_k$ in such a way that $S_{n_k}U_i$, $T_{n_k}U_i$ and $B_{n_k}$ as well as the reals $e^{n_k\chi_i(n_k)}/e^{n_k\chi_1(n_k)}$ ($i=1,\dots,d$) converge as $k\to+\infty$, with $L_i$, $R_i$, $B$ and $\alpha_i$ being their respective limits. From our assumption that $\left\Vert{P_n}/{\Vert P_n\Vert_2}-B_n\right\Vert_2\to 0$ together with (\ref{svd}) it follows that
$B=\sum_{i=1}^d\alpha_iL_iR_i^\star$ : because $\{L_i\}_{i=1}^d$ and $\{R_i\}_{i=1}^d$ are both orthonormal families, the fact that $BR_i=\alpha_iL_i$, together with  $\alpha_1=1$ and $\alpha_2>0$, means that $B$ is at least of rank $2$: however $B$ must be a rank one matrix, a contradiction.

\end{proof}

\begin{remark} (1) : Let $\sigma:\Omega\to\Omega$ be the full shift map on $\Omega=\{\0,\dots,\aa\}^\NNN$ and let $A:\Omega\to\mathcal{M}_d(\CCC)$ be a borelian map. We already saw that $(n,\omega)\mapsto P_n(\omega)=A(\omega)A(\sigma\cdot \omega)\cdots A(\sigma^{n-1}\cdot \omega)$ is a sub-multiplicative process whenever $\Omega$ is endowed with a $\sigma$-ergodic probability measure, say $\mu$. To emphasize the dependance on $\omega$ we note $\chi_k(n,\omega)$ ($1\le k\le d$) the  $n$-step Lyapunov exponents associated with $P_n(\omega)$ as in (\ref{SVD}).
By Oseledets Theorem, we know that $\chi_k(n,\omega)$ tends $\mu$-a.s. to the $k$-th  Lyapunov exponent $\chi_k(\mu)$. We apply  Theorem~\ref{r1} to $P_n=P_n(\omega)$ for $\omega$ a $\mu$-generic sequence: hence, in that case, a sufficient condition for $e^{n\chi_2(n,\omega)}/e^{n\chi_1(n,\omega)}\to0$ is that $\chi_1(\mu)>\chi_2(\mu)$; therefore, if the top-Lyapunov exponent $\chi_1(\mu)$ is strictly dominant then, there exists a sequence $B_1(\omega),B_2(\omega),\dots$ of rank one matrices s.t.
$$
\lim_{n\to+\infty}\left\Vert\frac{P_n(\omega)}{\Vert P_n(\omega)\Vert}-B_n(\omega)\right\Vert=0.
$$

(2) : Let $\mathcal{A}=(A_1A_2,\dots)$ be a sequence in $\mathcal{M}_d(\CCC)$ and  suppose  that $\left\Vert{P_n}/{\Vert P_n\Vert}-B_n\right\Vert\to 0$, where $P_n=A_1\dots A_n$  and each $B_n$ is a rank $1$ matrix. Then, there is no reason for the projective convergence of either the rows or the columns of $B_n$. Consider for instance  the matrices
$$
A=\begin{pmatrix}1&1&0\\0&1&0\\0&0&1\end{pmatrix},\ B=\begin{pmatrix}0&0&1\\1&0&0\\0&1&0\end{pmatrix},\ C=\begin{pmatrix}1&0&0\\1&1&0\\0&0&1\end{pmatrix},\ D=\begin{pmatrix}0&1&0\\0&0&1\\1&0&0\end{pmatrix}
$$
and consider that $P_n$ is the right product of the $n$ first matrices  of the infinite product 
$$
A^{2^0}BC^{2^1}DA^{2^2}BC^{2^3}D\cdots
$$
Then, for $n_k=2^0+2^1+\dots+2^{2k-1}+2k$ and $m_k=2^0+2^1+\dots+2^{2k}+2k+1$
$$
\lim_{k\to+\infty}B_{n_k}=\begin{pmatrix}0&\frac13&0\\0&0&0\\0&\frac23&0\end{pmatrix}\quad\hbox{and}\quad\lim_{k\to+\infty}B_{m_k}=\begin{pmatrix}\frac23&0&0\\0&0&0\\\frac13&0&0\end{pmatrix}.
$$
\end{remark}

About the divergence of the normalized infinite products of matrices, we prove the following proposition and corollary inspired from Elsner \& Friedland argument in \cite[Theorem~1]{EF97}.
\begin{proposition}\label{flgdhgsjf}Let $\mathcal{A}=(A_1,A_2,\dots)$ be a sequence in $\mathcal{M}_d(\CCC)$ and  $P_n:=A_1\cdots A_n$ be the $n$-step product; if $n\mapsto P_n/\Vert P_n\Vert$ converges (in $\mathcal{M}_d(\CCC)$) then, the matrices that occur infinitely many times in $\mathcal{A}$ (i.e. the $A$ such that $\#\{n\;;\;A_n=A\}=+\infty$) have  a common left-eigenvector.
\end{proposition}
\begin{proof}[{\bf Proof}] Suppose that $P_n/\Vert P_n\Vert\to P$ as $n\to+\infty$ and let $\lambda_n= \Vert P_n\Vert/\Vert P_{n-1}\Vert$. If $A_{n_k}=A$ for $n_1 <n_2 <\cdots$ then, 
$\lambda_{n_k}P_{n_k}/\Vert P_{n_k}\Vert={P_{n_k-1}A/\Vert P_{n_k-1}\Vert}\to PA$ as $k\to+\infty$:
hence $\lambda_{n_k}\to\Vert PA\Vert$  
and taking the limit as $k\to+\infty$ gives $\Vert PA\Vert P=PA$. Since $\Vert P\Vert=1$, there exists at least one $i$ s.t. $0\ne U_i^\star P$ so that $(U_i^\star P)A=U_i^\star(PA)=\Vert PA\Vert (U_i^\star P)$.

\end{proof}

\begin{corollary}Let $A(\0),\dots,A(\aaa)\in \mathcal{M}_d(\CCC)$ having no common left-eigenvector and for any $\omega=\omega_1\omega_2\cdots\in \Omega:=\{\0,\dots,\aaa\}^\NNN$ define $P_n(\omega)=A(\omega_1)\cdots A(\omega_n)$. Then, for any  Borel continuous probability  $\mu$ fully supported by $\Omega$, the sequence $n\mapsto P_n(\omega)/\Vert P_n(\omega)\Vert$ diverges for $\mu$-a.e.~$\omega\in \Omega$.
\end{corollary}
\begin{proof}[{\bf Proof}] 
For $\mu$-a.e. $\omega\in\Omega$ the sequence $\mathcal{A}(\omega)=(A(\omega_1),A(\omega_2),\dots)$ has infinitely many occurrences of each  $A(i)$: hence,  by Proposition~\ref{flgdhgsjf}, the sequence $P_n(\omega)/\Vert P_n(\omega)\Vert$ must be divergent.

\end{proof}

In the case of products of  nonnegative matrices, a straightforward consequence of parts (i) and (ii) of Theorem~A is the following theorem.

\begin{theorem}\label{pointwise}
Let $P_n=A_1\cdots A_n$ be the $n$-step right product of a sequence $\mathcal{A}=(A_1,A_2,\dots)$ made of $d\times d$ matrices with non negative entries, satisfying condition ${\bf (C)}$ and, for $1\le h\le H$, let $J_h(n)$ (resp. $V_h$) be the subset of $\{1,\dots,d\}$ (resp. the probability vector) given by Theorem~A; then, the following propositions hold:

(i) : each limit point of ${P_n}/{\Vert P_n\Vert}$ is a rank one matrix: actualy, for any $n\ge1$, there exists a probability vector $R_n$ such that 
$$
\lim_{n\to+\infty}\left\Vert{P_n\over \Vert P_n\Vert}-V_1R_n^\star\right\Vert=0\;;
$$ 

(ii) : recall that $h\diamond P_n$ is for the $d\times d$ matrix obtained from $P_n$ by replacing by $0$ the entries whose column indices are not in $J_{h}(n)$; then, each  limit point of $n\mapsto h\diamond P_n/\Vert h\diamond P_n\Vert$ is a rank one matrix: actually, for any $n\ge1$, there exists a probability vector $R_{n,h}$ such that 
$$
\lim_{n\to+\infty}\left\Vert{h\diamond P_n\over \Vert h\diamond P_n\Vert}-V_hR_{n,h}^\star\right\Vert=0
$$ 
(moreover, for any $2\le h\le H$ the ratio
$
\Vert h\diamond P_n\Vert/\Vert (h-1)\diamond P_n\Vert\to0$ as $n\to+\infty$).
\end{theorem}

\section{\bf Gibbs properties of linearly representable measures}\label{GibbsSection}

\subsection{Gibbs measures}Fix $\aaa\ge2$ an integer; an element in $\{\0,\dots,\aaa\}^n$ is written as a word $x_1\cdots x_n$ of length $n$, that is the (ordered) concatenation of the $n$ letters/digits $x_i$ in the alphabet $\{\0,\dots,\aaa\}$. By convention $\{\0,\dots,\aaa\}^0$ is reduced to the singleton $\{\emptyword\}$ whose unique element is called the empty word: we use the Kleene star notation for monoid of words on the alphabet $\{\0,\dots,\aaa\}$ (endowed with concatenation): more preciselly,
$$
\{\0,\dots,\aaa\}^*=\bigcup_{n=0}^{+\infty}\{\0,\dots,\aaa\}^n
$$
(the neutral element of $\{\0,\dots,\aaa\}^*$ being the empty word $\emptyword$).

The topology of the product space $\Omega:=\{\0,\dots,\aaa\}^\mathbb{N}$ is compact and given by the metric $(\omega,\xi)\mapsto 1/\aaa^{n(\omega,\xi)}$, where $n(\omega,\xi)=+\infty$ if $\omega=\xi$ while $n(\omega,\xi)=\min\{i\ge1\;;\;\omega_i\ne\xi_i\}$, otherwise.  Actually, for any $1/\aaa^{n+1}\le r<1/\aaa^n$, the (open-closed) ball of radius $r$ and centered at $\omega=\omega_1\omega_2\cdots$ 
is the so called {\it cylinder set} of  the sequences $\xi\in\Omega$ such that $\xi_1\cdots \xi_n=\omega_1\cdots \omega_n$, denoted by $[\omega_1\cdots \omega_n]$. We shall also consider the shift map $\sigma:\Omega\to\Omega$ s.t. $\sigma\cdot\omega=\sigma(\omega_1\omega_2\cdots)=\omega_2\omega_3\cdots$, which is an expanding continuous transformation of $\Omega$. (We write a sequence $\omega$ in $\Omega$ as one-sided infinite word $\omega_1\omega_2\cdots$ whose digits $\omega_i$ are in $\{\0,\dots,\aaa\}$.) 
The notion of Gibbs measure we are concerned with, is the one originally introduced by Bowen in his Lecture Notes \cite{Bow74}; here, it is not necessary to enter into the details of the underlying theory, since only a few elementary facts are needed. Suppose that $\mu$ is a Borel probability on $\Omega=\{\0,\dots,\aaa\}^\NNN$ which, for simplicity, we suppose to be fully supported by $\Omega$. The $n$-step potential $\phi_n:\Omega\to\RRR$ associated with $\mu$ is $\phi_1(\omega)=\log\mu[\omega_1]$ if $n=1$ and for $n\ge2$, 
\begin{equation}\label{nstep}
\phi_n(\omega):=\log\left({\mu[\omega_1\cdots \omega_n]\over \mu[\omega_2\cdots \omega_n]}\right).
\end{equation}
The fact that $\mu$ is fully supported ensures the functions $\phi_n$ to be well defined and continuous. Moreover, it is easily checked that 
$$
\mu[\omega_1\cdots \omega_n]=\exp\left(\sum\nolimits_{k=n}^1\phi_k(\sigma^{n-k}\cdot\omega)\right),
$$
this last identity leading to the following proposition.
\begin{proposition}\label{symbweakGibbs}
Suppose the $\phi_n$ converge uniformly to $\phi$ on $\Omega$, as~$n\to+\infty$: then,
$$
{1\over K_n}\le{\mu[\omega_1\cdots \omega_n]\over \exp\big(\sum_{k=0}^{n-1}\phi(\sigma^k\cdot\omega)\big)}\le K_n,
\quad\text{where}\quad
K_n:=\exp\left(\sum\nolimits_{k=1}^n\Vert \phi-\phi_n\Vert_\infty\right)
$$
and $\mu$ is a weak Gibbs measure (c.f. Yuri \cite{Yur98}), because (use $\Vert \phi-\phi_n\Vert_\infty\to0$ as $n\to+\infty$)
$$
\lim_{n\to+\infty}{1\over n}\log K_n=0\;;
$$
moreover, if $K:=\sup_n\{K_n\}<+\infty$, then $\mu$ is a Gibbs measure  (in the sense of Bowen \cite[Theorem~1.3]{Bow74}), since
$$
{1\over K}\le{\mu[\omega_1\cdots \omega_n]\over \exp\big(\sum_{k=0}^{n-1}\phi(\sigma^k\cdot\omega)\big)}\le K.
$$
\end{proposition}

\subsection{Linearly representable measures}Now, suppose 
that $A(\0),\dots,A(\aaa)$ are $d\times d$ matrices with nonnegative entries; moreover, $R$ is a fixed probability vector and  $\omega=\omega_1\omega_2\cdots$ is a given sequence in $\Omega$; then we usually note 
$$
A(\omega_1\cdots\omega_n):=A(\omega_1)\cdots A(\omega_n)
$$
and by definition, the $n$-step probability vector about $\omega$ and $R$ is 
\begin{equation}\label{projratio}
\Pi_n(\omega,R):={A(\omega_1\cdots\omega_n)R\over\left\Vert A(\omega_1\cdots\omega_n)R\right\Vert}.
\end{equation}
 Actually, the ratio in (\ref{projratio}) is not necessarily defined for any $\omega\in\Omega$ and we are led to introduce
\begin{equation}\label{defOmegaR}
\Omega_{R}:=\Big\{\omega\in\Omega\;;\;\forall n\ge1,\ A(\omega_1\cdots\omega_n)R\ne0\Big\}.
\end{equation}
This set is shift invariant in the sense that $\sigma(\Omega_R)=\Omega_R$ and compact (in $\Omega$), because for fixed $n$, the set of words $\omega_1\cdots\omega_n\in\{\0,\dots,\aaa\}^n$ such that $A(\omega_1\cdots\omega_n)R\ne0$ is obviously finite.
We note $A_*:=A(\0)+\dots+A(\aaa)$ and suppose in addition that $R$ satisfies the eigen-equation $A_*R=R$ (if $A_*R=\alpha R$ for some $\alpha>0$, there is no loss of generality to replace $A_i$ by $A_i/\alpha$). Let $L$ be a column vector with nonnegative entries for which $L^\star R=1$.
Then, we define the map $\mu:\{\0,\cdots,\aaa\}^*=\bigcup_{n=0}^{+\infty}\{\0,\dots,\aaa\}^n\to[0,1]$ by setting 
\begin{equation}\label{MPM}
\mu(\omega_1\cdots \omega_n)=L^\star A(\omega_1\cdots\omega_n)R
\end{equation}
An application of Kolmogorov Extension Theorem ensures $\mu$ to extend to a unique Borel probability measure on $\Omega$, that we abusively note $\mu$: this measure is determined by the condition 
$\mu[\omega_1\cdots \omega_n]=\mu(\omega_1\cdots \omega_n)$ to be satisfied for any word $\omega_1\cdots \omega_n\in\{\0,\cdots,\aaa\}^*$.
A measure defined by means of matrix products, like for instance in (\ref{MPM}), is usually called {\it linearly representable} (see 
\cite{BP11} for more details).   


The measure $\mu$ has support in $\Omega_R$ defined in (\ref{defOmegaR}) and we shall consider it as a measure on $\Omega_R$ instead; moreover,   for the sake of simplicity, we shall now assume  that $\mu$ is actually fully supported by $\Omega_R$, that is $\mu[\omega_1\cdots\omega_n]>0$ for any $\omega\in\Omega_R$ and any $n\ge1$ (by defintion of $\Omega_R$, this holds  for instance if $L$ has positive entries). In view of Gibbs properties, the main point is to start from Proposition~\ref{symbweakGibbs} and look at the  convergence (as $n\to+\infty$) of the $n$-step potential  
$\phi_n:\Omega_R\to\RRR$
such that, for any $\omega=\omega_1\omega_2\cdots\in\Omega_R$,  
$$
\phi_{n}(\omega):=\log\left({\mu[\omega_1\cdots \omega_{n}]\over \mu[\omega_2\cdots \omega_{n}]}\right)=
\log\left({L^\star A(\omega_1)\Pi_{n-1}(\sigma\cdot\omega,R)\over L^\star \Pi_{n-1}(\sigma\cdot\omega,R)}\right).
$$
We emphasize on the fact that existence and  continuity of the limit map $\VV:\Omega_R\to \mathcal{S}_d$ with
\begin{equation}\label{profromA}
\VV(\omega):=\lim_{n\to+\infty}\Pi_n(\omega,R)
\end{equation}
does not prevent from the possibility that $L^\star A(\omega_1)\VV(\sigma\cdot\omega)=0$ for some $\omega\in\Omega_R$ (the map $\omega\mapsto \VV(\omega)$ --~provided it exists~-- should be related to the {\it top Lyapunov direction map} introduced in \S~\ref{toplyapmap}).
This remark leads to the following proposition.
\begin{proposition}\label{symbweakbis}
If $\Pi_n(\omega,R)\to\VV(\omega)$, uniformly over $\Omega_R$ as $n\to+\infty$ and $A(i)\VV(\omega)\ne0$, for any $i\in\{\0,\dots,\aaa\}$ and any $\omega\in\Omega_R$ then, for any $L$ with positive entries,
$$
\lim_{n\to+\infty}\log\left({L^\star A(\omega_1)\Pi_{n-1}(\sigma\cdot\omega,R)\over L^\star \Pi_{n-1}(\sigma\cdot\omega,R)}\right)=\log\left({L^\star A(\omega_1)\VV(\sigma\cdot\omega)\over L^\star\VV(\sigma\cdot\omega)}\right)
$$
uniformly over $\Omega_R$: moreover, if $A_*R=R$ and $L^\star R=1$ then $\mu$ defined in (\ref{MPM}) is a weak Gibbs measure.
\end{proposition}
\begin{proof}[{\bf Proof}]If two real-valued continuous maps $f_n,g_n$ converge uniformly on a compact set, to positive limits $f$ and $g$ respectively, then $\displaystyle\log{f_n}/{g_n}$ converges uniformly to $\displaystyle\log{f}/{g}$.

\end{proof}

\subsection{Conditions for the uniform convergence}\label{u}
Pointwise convergence of $\Pi_n(\omega,R)$ (provided it holds), 
does not imply uniform convergence over $\Omega_R$. To see this, consider for instance 
$$
A(\0)=\begin{pmatrix}1&0&1\\0&1&0\\0&0&0\end{pmatrix},\;A(\1)=\begin{pmatrix}1&0&1\\1/2&0&0\\0&0&0\end{pmatrix}
\quad\text{and}\quad R=\begin{pmatrix}1/3\\1/3\\1/3\end{pmatrix}.
$$ 
Because $A(\0)$ and $A(\1)^2=A(\1\1)$ are both idempotent matrices,  
${\Pi_n(\omega,R)}\to(2/3,1/3,0)$ as $n\to+\infty$, this pointwise convergence being valid 
for any $\omega\in\Omega_R=\Omega=\{\0,\1\}^\mathbb{N}$.
However, the convergence is not uniform on $\Omega$ because $\Pi_n(\omega,R)=(4/5,1/5,0)$, for any $\omega\in[\0^{n-1}\1]$ with $n\ge1$.
Similarly, for $L=(1,1,1)$ the convergence, for $n\to+\infty$, of the potential 
$$
\displaystyle\phi_n(\omega)=\log\frac{L^\star A(\omega_1\dots\omega_n)R}{L^\star A(\omega_2\dots\omega_n)R}
$$ 
is not uniform on $\Omega$ because (for $n\ge2$) either $\phi_n(\omega)=\log(6/5)$ if $\omega\in[\1\0^{n-2}\1]$ or $\phi_n(\omega)=0$ otherwise
and yet, the limit $\phi$ of $\phi_n$ is continuous, since  $\phi\equiv0$.
The following lemma gives a criterion for uniform convergence of $\Pi_n(\omega,R)$. 
\begin{lemma}\label{uniform}
$n\mapsto\Pi_n(\cdot,R)$ is uniformly convergent over  $\Omega_{R}$ if and only if for any $\omega\in\Omega_{R}$
\begin{equation}\label{modifunif}
\lim_{n\to+\infty}\left(\sup
\Big\{\left\Vert \Pi_{n+r}(\xi,R)-\Pi_{n}(\xi,R)\right\Vert\;;\;\xi\in[\omega_1\cdots\omega_n]\cap\Omega_{R}\;\text{and}\;r\ge0\Big\}\right)=0.
\end{equation}
\end{lemma}
\begin{proof}[{\bf Proof}] The direct implication is given by the Cauchy criterion.
Conversely, suppose that (\ref{modifunif})~holds. Given $\varepsilon>0$ and $\omega\in\Omega_{R}$ there exists a rank $n(\omega)\ge1$ such that for~any~$r,s\ge n(\omega)$
$$
\xi\in[\omega_1\cdots\omega_{n(\omega)}]\cap \Omega_{R}\;\Longrightarrow\;\left\Vert \Pi_r(\xi,R)-\Pi_s(\xi,R)\right\Vert\le\varepsilon.
$$
Each cylinder  $[\omega_1\dots\omega_{n(\omega)}]$ is an open set containing $\omega$. Because $\Omega_{R}$ is compact, it is covered by finitely many of such cylinders, say $[\omega^{i}_1\dots\omega^{i}_{n(\omega^i)}]$ for $i=1,\dots,N$. Hence if $\displaystyle r,s\ge\max\nolimits_i\{n(\omega^i)\}$ the inequality $\left\Vert \Pi_r(\xi,R)-\Pi_s(\xi,R)\right\Vert\le\varepsilon$ holds for any $\xi\in\Omega_{R}$: this proves that $\Pi_n(\cdot,R)$ is uniformly Cauchy and converges uniformly.

\end{proof}

For a fixed $R\in\mathcal{S}_d$, we now  consider the problem of uniform convergence of $\Pi_n(\cdot,R)$ toward $\VV(\cdot)$ within the framework of Theorem A. The point is to deal with  conditions ${\bf (C)}$ w.r.t. each sequence $\mathcal{A}(\omega)=(A(\omega_1),A(\omega_2),\dots)$ for $\omega$ running over the whole $\Omega_R$. We shall say that $\omega\in\Omega_R$ is  ${\bf (C)}$-regular (resp. ${\bf (C)}$-singular) if $\mathcal{A}(\omega)$ satisfies (resp. does not satisfy) condition ${\bf (C)}$. Actually, one possible --~consistent~--  difficult point (as we shall see for instance in Section~\ref{7X7}) is the existence of $\sss\in \{\0,\cdots,\aaa\}$ so that  $\bar\ss=\sss\sss\cdots$ (and thus each sequence in  $\sigma^{-n}\{\bar\sss\}$) is  ${\bf (C)}$-singular. For $\omega=\omega_1\omega_2\cdots\in\Omega_R$ and $1\le i\le j$, we use the notation 
$$
\omega_{i,j}=\omega_{i+1}\dots\omega_j\;;
$$
in particular $\omega_{i,i}$ is the empty word $\emptyword$.
\goodbreak

\begin{proposition}[Pointwise convergence]\label{Simplecvg}
For $R\in\mathcal{S}_d$ and $\omega\in\Omega_R$, the  convergence of $\Pi_n(\omega,R)$ to a limit vector $\VV(\omega)\in\mathcal{S}_d$ holds  if $\omega$ satisfies at least one of the following conditions:

(i) : $\omega\in\Omega_R$ is ${\bf (C)}$-regular, there exist an integer $N$, a real $\Lambda\ge1$ and a sequence of integers $n\mapsto\psi(n)$, with  $\psi(n)\le n$ and $\psi(n)\to+\infty$ as $n\to+\infty$, s.t. for any $n\ge N$, any $r\ge0$,
$$
{\rm (S1)}\;:\;\Delta\Big(A(\omega_{\psi(n),n+r})R\Big)=\Delta\Big(A(\omega_{\psi(n),n})R\Big).
$$

(ii) : $\omega=w\bar\sss$, for $\sss\in\{\0,\dots,\aaa\}$ (with $\bar\sss$ possibly ${\bf (C)}$-singular) and there exists a $d\times d$ matrix $B_\sss$ with $\Vert B_\sss\Vert=1$, s.t.
$$
{\rm (S2.1)}\;:\;\lim_{n\to+\infty}\left\Vert {A(\ss)^n\over \Vert A(\sss)^n\Vert}-B_\sss\right\Vert=0,\quad
{\rm (S2.2)}\;:\;A(w)B_sR\ne0\;.
$$
\end{proposition}

\begin{proof}[{\bf Proof}](i) : Let $\omega\in\Omega_R$ be ${\bf (C)}$-regular and satisfy ${\rm (S1)}$ w.r.t. the sequence $n\mapsto\psi(n)$. Using the notations of Definition \ref{C}, let $k=\kk(\psi(n))$ be the integer such that $s_{k+1}\le\psi(n)<s_{k+2}$. We can replace the sequence $n\mapsto\psi(n)$ by the sequence $n\mapsto s_{\kk(\psi(n))}<\psi(n)$: this sequence also satisfies ${\rm (S1)}$ because the equality $\Delta\Big(A(\omega_{\psi(n),n+r})R\Big)=\Delta\Big(A(\omega_{\psi(n),n})R\Big)$ implies $\Delta\Big(A(\omega_{s_k,n+r})R\Big)=\Delta\Big(A(\omega_{s_k,n})R\Big)$. Moreover with this new definition of the sequence $n\mapsto\psi(n)$, by (\ref{asterix}) of Lemma \ref{increasing} the column vector $X_r:=A(\omega_{\psi(n),n+r})R/\Vert A(\omega_{\psi(n),n+r})R\Vert$ belongs to $\mathcal H_2(\Lambda')$ with $\Lambda'=\Lambda/(1-\lambda)$. Take $n\ge N$, $r\ge0$; from ${\rm (S1)}$, we know that $\Delta(X_r)=\Delta(X_0)$ and thus
$$
\min\Big\{h\;;\;\mathcal{I}(X_r)\cap J_h(\psi(n))\ne\emptyset\Big\}
=\min\Big\{1\le h\le H\;;\;\mathcal{I}(X_0)\cap J_h(\psi(n))\ne\emptyset\Big\}=:h_n(\omega),
$$ 
that is (with notations in Theorem~A) $1\le \hhh_{\psi(n)}(X_r)=\hhh_{\psi(n)}(X_0)=h_n(\omega)\le H$: we emphasize that $h_n(\omega)$  depends on $n$ and $\omega$ but not on~$r$ (this is due to the synchronization condition ${\rm (S1)}$). 
 For $\omega$ being supposed ${\bf (C)}$-regular, we can apply part (iii) of Theorem~A which gives real numbers $\varepsilon_n\to0$ (depending on $\omega$ but not on $r$)~s.t. 
$$
\left\Vert{A(\omega_{0,\psi(n)})X_r\over \Vert A(\omega_{0,\psi(n)})X_r\Vert}-V_{h_n(\omega)}\right\Vert\le \varepsilon_{\psi(n)}\Lambda'
$$
(here $V_1,\dots,V_H$ are the probability vectors associated with $\mathcal{A}(\omega)=(A(\omega_1),A(\omega_2),\dots)$ by Theorem~A). Therefore, using the triangular inequality, 
\begin{equation}\label{newequation}
\begin{array}{rcl}\Vert\Pi_{n+r}(\omega,R)-\Pi_{n}(\omega,R)\Vert
&\le &\displaystyle\left\Vert{A(\omega_{0,\psi(n)})X_r\over \Vert A(\omega_{0,\psi(n)})X_r\Vert}-V_{h_n(\omega)}\right\Vert+\\
&&\displaystyle\left\Vert{A(\omega_{0,\psi(n)})X_0\over \Vert A(\omega_{0,\psi(n)})X_0\Vert}-V_{h_n(\omega)}\right\Vert\le 2\varepsilon_{\psi(n)}\Lambda'
\end{array}
\end{equation}
and thus $n\mapsto\Pi_n(\omega,R)$ is a Cauchy sequence of $\mathcal{S}_d$.

(ii) : Let  $\omega=w\bar\sss$ where $w$ is either the empty word $\emptyword$ or  
 $w=\omega_1\cdots\omega_N$ for a $N\ge1$. Under condition ${\rm (S2.2)}$, we know that $A(w)B_\sss R\ne0$ (with the convention that $A(\emptyword)$ is the $d\times d$ identity matrix) and it follows from  ${\rm (S2.1)}$ that
$\Pi_n(w\bar\sss,R)\to{A(w)B_\sss R/\Vert (A(w)B_\sss R\Vert}$ as $n\to+\infty$.

\end{proof}

The following Proposition~\ref{uuuuuu} is an extrapolation of the above Proposition~\ref{Simplecvg}, in view of uniform convergence of $\Pi_n(\cdot,R)$. It gives (among many others) two possible situations for which the Cauchy condition (\ref{modifunif}) of Lemma~\ref{uniform} holds. It is intended to  illustrate how Theorem~A applies to establish (when it holds) the uniform convergence of $\Pi_n(\cdot,R)$. Apart~its~il\-lus\-tra\-tive content, we shall make a more specific  use of Proposition~\ref{uuuuuu} when dealing with the uniform convergence part of  Theorem~\ref{3matrices} considered in \S~\ref{unif7X7}. (We stress that uniform convergence of $\Pi_n(\cdot, R)$ over $\Omega_R$ is  anyway a technical~question.) The two items of Proposition \ref{uuuuuu} also correspond to the (inevitably technical) ideas we already use in some previous papers, to prove the uniformity of the convergence. For proving this, the main difficulty is related to the set $H_2(\Lambda)$, that is why this set is involved in the conditions of Proposition \ref{uuuuuu} and not in Proposition~\ref{Simplecvg}.

\begin{proposition}[Uniform convergence]\label{uuuuuu}
Let $R\in\mathcal{S}_d$ for which the set  of the ${\bf (C)}$-regular sequences in $\Omega_R$ is 
$\XX=\Omega_R\setminus\bigcup_{\sss\in S}\bigcup_{n=0}^{+\infty}\sigma^{-n}\{\bar\sss\}$,
where $S\subset\{\0,\dots,\aaa\}$; then,  the Cauchy condition  (\ref{modifunif}) of Lemma~\ref{uniform} holds at $\omega\in\Omega_R$ if at least one of the following conditions is satisfied:

(i) : $\omega\in\XX$, there exist an integer $N$, a real $\Lambda\ge1$ and a sequence of integers $n\mapsto\psi(n)$, with  $\psi(n)\le n$ and $\psi(n)\to+\infty$ as $n\to+\infty$, s.t. for any $n\ge N$, any $\xi\in[\omega_1\cdots\omega_n]\cap\Omega_R$ and $r\ge0$,
$$
{\rm (U1.1)}\;:\;A(\xi_{\psi(n),n+r})R\in\mathcal H_2(\Lambda),
\quad{\rm (U1.2)}\;:\;\Delta\Big(A(\xi_{\psi(n),n+r})R\Big)=\Delta\Big(A(\xi_{\psi(n),n})R\Big).
$$

(ii) : $\omega\in\Omega_R$ and $\omega=w\bar\sss$, where $w$ is a (possibly empty) word on $\{0,\dots,\aaa\}$ and $\sss\in S$ is a digit for which {\rm (U2)}  : there  exist  $C_\sss,D_\sss\in\mathcal{S}_d$ s.t.
$$
{\rm (U2.1)}\;:\;\lim_{n\to+\infty}\left\Vert{A(\ss)^n\over \Vert A(\sss)^n\Vert}-C_\sss D_\sss^\star\right\Vert=0,\quad
{\rm (U2.2)}\;:\;A(w)C_\sss\ne0\;;
$$
and {\rm (U3)} : there exist an integer $N$, a real $\Lambda\ge1$ and a sequence of integers $n\mapsto\varphi(n)$, with  $\varphi(n)\le n$ and $\varphi(n)\to+\infty$ as $n\to+\infty$ and such that for any $n\ge N$, any $\xi\in[\omega_1\cdots\omega_n]\cap\Omega_R$ and $r\ge0$, there exist $m\in\{\varphi(n),\dots,n+r\}$ for which the following three propositions hold:
$$
{\rm (U3.1)}\;:\;\xi\in[\omega_1\cdots\omega_m],\quad
{\rm (U3.2)}\;:\;A(\xi_{m,n+r})R\in\mathcal H_2(\Lambda),\quad
{\rm (U3.3)}\;:\;D_\sss^\star A(\xi_{m,n+r})R\ne0.
$$
\end{proposition}

\begin{proof}[\bf Proof](i) : Let $\omega\in\XX$ for which ${\rm (U1.1)}$ and ${\rm (U1.2)}$ are satisfied with the integer $N$, the constant $\Lambda$ and the sequence $n\mapsto\psi(n)$. Take $n\ge N$,  $\xi\in[\omega_1\cdots\omega_n]\cap\Omega_R$, $r\ge0$ and define 
$X_r(\xi):={A(\xi_{\psi(n),n+r})R/\Vert A(\xi_{\psi(n),n+r})R\Vert}$. 
Notice that $X_0(\xi)=X_0(\omega)=:X_0$; then, from  ${\rm (U1.2)}$ we know that $\Delta(X_r(\xi))=\Delta(X_0)$ and thus
$$
\min\Big\{h\;;\;\mathcal{I}(X_r(\xi))\cap J_h(\psi(n))\ne\emptyset\Big\}
=\min\Big\{h\;;\;\mathcal{I}(X_0)\cap J_h(\psi(n))\ne\emptyset\Big\}=:h_n(\omega),
$$ 
that is (with notations in Theorem~A) $1\le \hhh_{\psi(n)}(X_r(\xi))=\hhh_{\psi(n)}(X_0)=h_n(\omega)\le H$: we emphasize that $h_n(\omega)$ does depend on $n$ and $\omega$ but neither on $\xi\in[\omega_1\cdots\omega_n]\cap \Omega_R$ nor on~$r$ (this is due to the synchronization condition ${\rm (U1.2)}$). 
 For $\omega$ being supposed ${\bf (C)}$-regular, we can apply part (iii) of Theorem~A which gives real numbers $\varepsilon_n\to0$ (depending on $\omega$ but not on $r$)~s.t. 
$$
\Vert\Pi_{n+r}(\xi,R)-V_{h_n(\omega)}\Vert=\left\Vert{A(\omega_{0,\psi(n)})X_r(\xi)\over \Vert A(\omega_{0,\psi(n)})X_r(\xi)\Vert}-V_{h_n(\omega)}\right\Vert\le \varepsilon_{\psi(n)}\Lambda_{X_r(\xi)}
$$
(here $V_1,\dots,V_H$ are the probability vectors associated with $\mathcal{A}(\omega)=(A(\omega_1),A(\omega_2),\dots)$ by Theorem~A). Therefore, using the triangular inequality, 
\begin{align*}
\Vert\Pi_{n+r}(\omega,R)-\Pi_{n}(\omega,R)\Vert
\le &\left\Vert{A(\omega_{0,\psi(n)})X_r(\xi)\over \Vert A(\omega_{0,\psi(n)})X_r(\xi)\Vert}-V_{h_n(\omega)}\right\Vert+\\
&\left\Vert{A(\omega_{0,\psi(n)})X_0\over \Vert A(\omega_{0,\psi(n)})X_0\Vert}-V_{h_n(\omega)}\right\Vert\le \varepsilon_{\psi(n)}\Big(\Lambda_{X_r(\xi)}+\Lambda_{X_0}\Big)\;;
\end{align*}
however, by ${\rm (U1.1)}$ one has $\Lambda_{X_r(\xi)},\Lambda_{X_0}\le \Lambda$ and thus the uniform Cauchy condition (\ref{modifunif}) of Lemma~\ref{uniform} holds for $\omega\in X$ for which ${\rm (U1.1)}$-${\rm (U1.2)}$ are satisfied: this proves part (i). 

(ii) :  Let $\omega=w\bar\sss\in\Omega_R$, where  $\sss\in S$ and $w$ is either the empty word $\emptyword$ (and $k_0:=0$) or  
$w=\omega_1\cdots\omega_{k_0}$ for a $k_0\ge1$. Consider $n\ge1$ and $r\ge0$ arbitrary given; for $\xi\in\Omega_R$ s.t. $\xi_1\cdots\xi_n=\omega_1\cdots\omega_n=w\sss^{n'}$ we note $\zeta=\xi$ if $w=\emptyword$ or $\zeta=\sigma^{k_0}\cdot\xi$ otherwise. We suppose that ${\rm (U2)}$ and ${\rm (U3)}$ are satisfied by $\omega$ and, assuming that $\varphi(n)\ge k_0$, we put $n'=n-k_0$, $m'=m-k_0$. By definition of $C_\sss$ and $D_\sss$,  conditions in  ${\rm (U2)}$ ensure $A(\sss)^{m'}=C_\sss D_\sss^\star+M_{m'}$ with $\Vert M_n\Vert\to0$ as $n\to+\infty$, while $A(w)C_\sss\ne0$: in particular, it is licit to put $B(w):=A(w)/\Vert A(w)C_\sss\Vert$. Defining   
$Y:=A(\zeta_{m',n'+r})R/\Vert A(\zeta_{m',n'+r})R\Vert$ and $U:=U_1+\cdots+U_d$, one gets the following identities:
\begin{align*}
{A(\xi_{0,n+r})R\over \Vert A(\xi_{0,n+r})R\Vert}-B(w)C_\sss
&={B(w)A(\zeta_{0,m'})Y\over \Vert B(w)A(\zeta_{0,m'})Y\Vert}-
B(w)C_\sss\\
&=B(w)
\left({C_\sss D_\sss^\star Y+M_{m'}Y\over U^\star B(w)C_\sss D_\sss^\star Y+U^\star M_{m'}Y}-C_\sss\right)\\
\text{(recall that $U^\star B(w)C_\sss=1$)}\quad&=B(w)
\left({C_\sss D_\sss^\star Y+M_{m'}Y
-D_\sss^\star YC_\sss-U^\star M_{m'}YC_\sss)\over D_\sss^\star Y+U^\star M_{m'}Y}\right)\\
&=B(w)
\left(\displaystyle{M_{m'}Y-U^\star M_{m'}YC_\sss\over D_\sss^\star Y+U^\star M_{m'}Y}\right).
\end{align*}
On the one hand,  $\Vert M_{m'}Y-U^\star M_{m'}YC_\sss\Vert \le 2\Vert M_{m'}\Vert$. On the other hand, to lower bound $D_\sss^\star Y+U^\star M_{m'}Y$ we use conditions ${\rm (U3.2)}$ and {\rm (U3.3)} ensuring
$Y\in\mathcal{H}_2(\Lambda)$ and $D_\sss^\star Y\ne0$. Indeed, $Y(i)\ne0$ implies $Y(i)\ge 1/\Lambda$ and $D_\sss^\star Y\ge \underline{D}_\sss/\Lambda$, where $\underline{D}_\sss=\min\{D_\sss(i)\ne0\}$. Finally   
$D_\sss^\star Y+U^\star M_{m'}Y\ge\underline{D}_\sss/\Lambda-\Vert M_{m'}\Vert,
$
and thus (for $n$ large enough so that $\underline{D}_\sss/\Lambda-\Vert M_{m'}\Vert>0$),
\begin{align*}
\Vert\Pi_{n+r}(\xi,R)-\Pi_{n}(\xi,R)\Vert=\left\Vert {A(\xi_{0,n+r})R\over \Vert A(\xi_{0,n+r})R\Vert}-
{A(\xi_{0,n})R\over \Vert A(\xi_{0,n})R\Vert}\right\Vert
&\le {4\Vert B(w)\Vert\cdot \Vert M_{m'}\Vert\over \underline{D}_\sss/\Lambda-\Vert M_{m'}\Vert}.
\end{align*}
By definition, $\varphi(n)-k_0\le m'$ while $\varphi(n)\to +\infty$ as $n\to+\infty$: because $\Vert M_{m'}\Vert\to0$ as $m'\to+\infty$, the Cauchy  condition (\ref{modifunif}) of Lemma~\ref{uniform} is satisfied.

\end{proof}

\section{\bf The Kamae measure\protect\footnote{This section may be skipped by readers not familiar with the symbolic dynamics of self-affine sets.}}\label{KAMAESEC}

\begin{figure}[H]
  \begin{center}
       \includegraphics[scale=0.7]{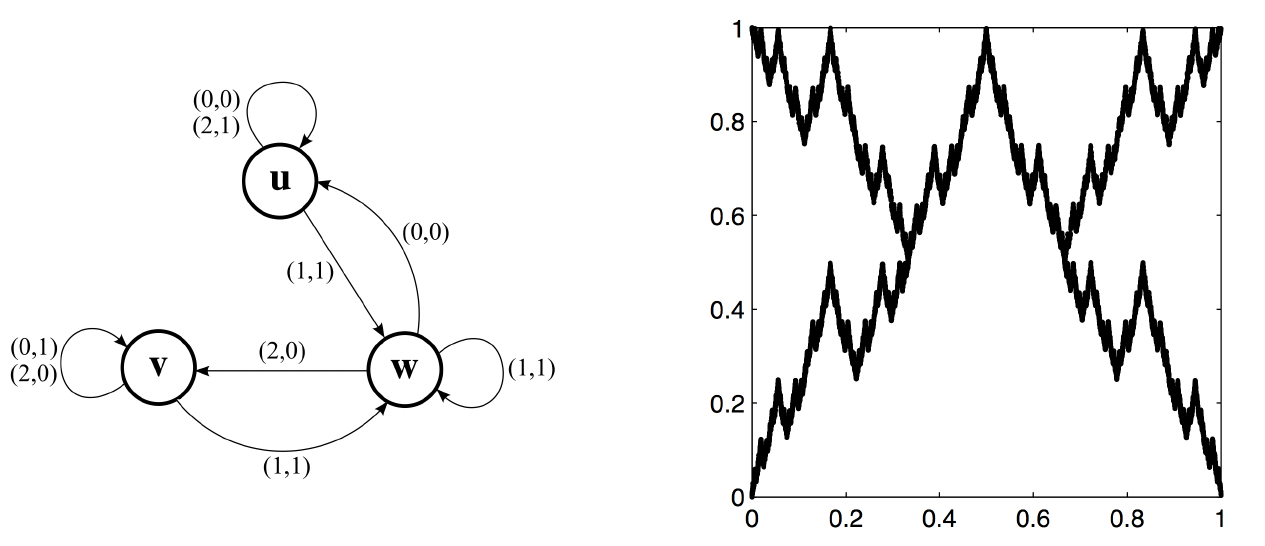}
\caption{\label{KAMAEBIS}\leftskip=10mm \rightskip=10mm
\small  \it The adjacency graph (left) defining the sofic system related to  Kamae's sofic affine subset of $\TTT^2$ (right).} 
  \end{center}
\end{figure}

Let $\TTT^2=\RRR/\ZZZ\times\RRR/\ZZZ$ (identified with $[0\,;1[\times[0\,;1[$) be the $2$-torus and  $T:\TTT^2\to\TTT^2$ s.t. $T(x,y)=(\{3x\},\{2x\})\;;$ (here $\{x\}$ stands for the fractional part of the real number $x$). A {\it $T$-invariant carpet of $\TTT^2$} is a compact subset $\KK$ of $\TTT^2$ which is $T$-invariant in the sense that $T(\KK)=\KK$. The symbolic model of $T$ is given by the full shift map  $\sigma:\Sigma\to\Sigma$, where
$$
\Sigma=\Big\{(\0,\0),(\1,\0),(\2,\0),(\0,\1),(\1,\1),(\2,\1)\Big\}^\NNN.
$$
It is associated with the representation map ${\rm Rep}:\Sigma\to\TTT$ such that
$$
{\rm Rep}\Big((x_1,y_1)(x_2,y_2)\cdots\Big)=\left(\left\{\sum\limits_{k=1}^{+\infty}{x_k\over 3^{k}}\right\},\left\{\sum\limits_{k=1}^{+\infty}{y_k\over 2^{k}}\right\}\right).
$$
Any subshift $\ZZ\subset \Sigma$ (i.e. $\ZZ$ is a compact subset of $\Sigma$ which is shift-invariant in the sense that $\sigma(\ZZ)=\ZZ$) is associated with the $T$-invariant carpet $\KK={\rm Rep}(\ZZ)$: we call $\ZZ$ the symbolic model of $\KK$. (Conversely any $T$-invariant carpet of $\TTT^2$, has a symbolic model.) The fractal geometry of these subsets of the $2$-torus, has been widely studied in the framework of the so-called {\it Variational Principle for Dimension} (see for instance \cite{McM84,Bed84,KP96a}). In this paragraph we shall focus our attention on what we call the {\it  Kamae carpet} (see \cite{Kam86} and \cite[\S~6]{Oli10b}), that is  $\KK={\rm Rep}(\ZZ)$ where $\ZZ$ is the sofic system specified by the adjacency scheme 
\begin{align*}
\gotA=&\left\{A_{(\0,\0)}=
\begin{pmatrix}1&0&0\cr0&0&0\cr1&0&0\end{pmatrix},\; 
A_{(\1,\0)}=
\begin{pmatrix}0&0&0\cr0&0&0\cr 0&0&0\end{pmatrix},\; 
A_{(\2,\0)}=
\begin{pmatrix}0&0&0\cr0&1&0\cr 0&1&0\end{pmatrix},\right.\\
&\left.\;\;\;
A_{(\0,\1)}=
\begin{pmatrix}0&0&0\cr0&1&0\cr0&0&0\end{pmatrix},\; 
A_{(\1,\1)}=
\begin{pmatrix}0&0&1\cr0&0&1\cr 0&0&1\end{pmatrix},\; 
A_{(\2,\1)}=
\begin{pmatrix}1&0&0\cr0&0&0\cr 0&0&0\end{pmatrix}\right\}
\end{align*}
(An equivalent form of $\gotA$ is given by the adjacency graph in Fig.~\ref{KAMAEBIS} : see \cite{DGS76} and \cite{LM95} for the theory of shift systems). By the theory of sofic systems, 
$\xi=(x_1,y_1)(x_2,y_2)\cdots\in\ZZ$ if and only if 
$\Vert A_{(x_1,y_1)}\cdots A_{(x_n,y_n)}\Vert>0$ for any $n\ge1$.
The adjacency matrix of $\gotA$, that is 
$$
A_*=\sum\nolimits_{i,j}A_{(i,j)}=\begin{pmatrix}2&0&1\cr0&2&1\cr 1&1&1\end{pmatrix},
$$
has a spectral radius $\rho_{A_*}=3$. Hence, the topological entropy of $\sigma:\ZZ\to\ZZ$ is  $\mathrm{h}_\sigma(\ZZ)=\log3$: the Parry-Bowen measure over $\ZZ$  is the unique shift-ergodic measure $\bar\eta$ supported by $\ZZ$ and whose Kolmogorov-Sinaï entropy $\hh_\sigma(\bar\eta)$ coincides with the topological entropy $\hh_\sigma(\ZZ)$: in other words:
$$
\hh_\sigma(\bar\eta)=\hh_\sigma(\ZZ)=\log 3
$$
(for the concepts of entropy --~not decisive in the sequel~-- see \cite{DGS76}). Let $\yy:\Sigma\to\{\0,\1\}^\NNN=:\YY$ be the map 
$$
(x_1,y_1)(x_2,y_2)\cdots\mapsto \yy\big((x_1,y_1)(x_2,y_2)\cdots\big)=y_1y_2\cdots\;;
$$
it is easily seen that $\yy$ makes\footnote{We abusively use $\sigma$ for the shift maps over symbolic spaces as soon as the digits are well specified.} $\sigma:\YY\to\YY$ a factor of $\sigma:\Sigma\to\Sigma$ ($\yy$ is continuous surjective and $\yy\circ\sigma=\sigma\circ\yy$). Because $\yy(\ZZ)=\YY$, the restriction $\yy:\ZZ\to\YY$ is also a factor map.
Actually,  the Gibbs properties of the $\yy$-axis projection $\bar\eta_\yy=\bar\eta\circ \yy^{-1}$ of $\bar\eta$ over $\YY$, plays a crucial role w.r.t. the fractal geometry of $\KK$ (see \cite[Theorem A]{Oli10a}). The measure $\bar\eta_\yy$ --~that we call {\it the Kamae measure}~-- turns out to be linearly representable, since for any binary word $\xi_1\cdots \xi_n$, 
\begin{equation}\label{Kamae}
\bar\eta_\yy[y_1\cdots y_n]={1\over 3^n}\Vert A(y_1\cdots y_n)\Vert,
\end{equation}
where $A(\0)=A_{(\0,\0)}+A_{(\1,\0)}+A_{(\2,\0)}$ and $A(\1)=A_{(\0,\1)}+A_{(\1,\1)}+A_{(\2,\1)}$, so that
$$
A(\0)=\begin{pmatrix}1&0&0\cr0&1&0\cr 1&1&0\end{pmatrix}
\quad\text{and}\quad
A(\1)=
\begin{pmatrix}1&0&1\cr0&1&1\cr 0&0&1\end{pmatrix}.
$$
To analyze the Gibbs property of $\bar\eta_\yy$, we need  prove the following proposition.
\begin{proposition}\label{gfdjsk} The continuous functions $\Pi_n(\cdot,U):\YY\to\mathcal{S}_3$ ($n\ge1$), where, for any $y\in\YY$, 
$$
\Pi_n(y,U):={A(y_1\cdots y_n)U\over \Vert A(y_1\cdots y_n)U\Vert}
\quad\text{and}\quad
U=\begin{pmatrix}1\\1\\1\\\end{pmatrix}.
$$ 
converge  toward a limit function  $\VV:\YY\to\mathcal{S}_3$ uniformly over $\YY$. 
\end{proposition}
Proposition~\ref{gfdjsk} may be establishes by elementary direct methods. However, to illustrates the technics developed in the present paper, we shall remain in the framework of  Theorem~A, using in particular Lemma~\ref{uniform} and the technics involved in Propositions~\ref{Simplecvg}~and~\ref{uuuuuu}.
\begin{proposition}
The set of sequences  $y\in\YY$ ($=\{\0,\1\}^\NNN$) which are ${\bf (C)}$-regular  --~in the sense that $\mathcal{A}(y)=(A(y_1),A(y_2),\dots)$ satisfies condition ${\bf (C)}$~-- is 
$$
\YY':=\left(\YY\setminus\bigcup_{n=0}^{+\infty}\sigma^{-n}\{\bar\0\}\right)\cap\left(\YY\setminus\bigcup_{n=0}^{+\infty}\sigma^{-n}\{\bar\1\}\right).
$$
\end{proposition}
\begin{proof}[{\bf Proof}]To begin with, $A(\0)$ is an idempotent matrix and for $a\ge 0$ one has the matrix identities 
\begin{equation}\label{KAMH1}
A(\0)A(\1)^aA(\0)=\begin{pmatrix}a+1&a&0\cr a&a+1&0\cr 2a+1&2a+1&0\end{pmatrix}=:B_a;
\end{equation}
moreover, a simple induction shows that for $a_1,a_2,\dots$ a given sequence of positive integers   
\begin{equation}\label{fhgkds}
B_{a_1}\cdots B_{a_n}=B_{\alpha_n}
\quad\text{where}\quad
\alpha_n={1\over 2}(2a_1+1)\cdots(2a_n+1)-{1\over 2}
\end{equation}
and this proves $\mathcal{A}(y)$ satisfies ${\bf (C)}$, for any $y\in\YY'$. The matrix $A(\0)$ is idempotent and does not belong to $\mathcal{H}_1$: hence, $\bar\0$ is a ${\bf (C)}$-singular sequence. Likewise, $\bar\1$ is ${\bf (C)}$-singular, since 
\begin{equation}\label{powern}
A(\1)^n=\begin{pmatrix}1&0&n\cr0&1&n\cr 0&0&1\end{pmatrix}.
\end{equation}

\end{proof}

The simple convergence of $\Pi_n(y,U)$ toward a limit  $\VV(y)\in\mathcal{S}_d$, holds for any $y\in \YY$. For $y\in\YY'$, the existence of the limit vector $\VV(y)$ may be obtained  by application of Theorem~A (actually Corollary~A). A more  straightforward  approach gives  the expression of $\VV(y)$: indeed, since $\0\1$ must factorize $y$ infinitely many times, it is necessary~that 
$$
y=\1^{a_0}\0^{b_0}\1^{a_1}\0^{b_1}\1^{a_2}\cdots=\1^{a_0}\0^{*}\1^{a_1}\0^{*}\1^{a_2}\cdots,
$$
where  $a_i$ and $b_i$ are positive integers, with the possible exception that $a_0\ge0$ (and  the convention that $\1^0=\emptyword$); here, we also write $\0^*$ for a non empty arbitrary block of $\0$. Decomposing $y$ in blocks of the form $\0^*\1^{a_i}\0^*$ (for $a_i>0$) and using the fact that $A(\0)$ is idempotent together with (\ref{KAMH1})-(\ref{fhgkds})-(\ref{powern}), one finds: 
\begin{equation}\label{1a00}
\VV(y)={1\over a_0+1}\cdot\begin{pmatrix}1&0&a_0\cr0&1&a_0\cr 0&0&1\end{pmatrix}
\begin{pmatrix}
1/4\\1/4\\1/2
\end{pmatrix}.
\end{equation}
The simple convergence also holds about sequence of the form $y=w\bar\0$ and  $y=w\bar\1$, since 
$$
\VV(w\bar\0)={A(w)A(0)U\over \Vert A(w)A(0)U\Vert}\quad\text{and}\quad \VV(w\bar\1)={A(w)(U_1+U_2)\over \Vert A(w)(U_1+U_2)\Vert}.
$$

\begin{proof}[{\bf Proof of Proposition~\ref{gfdjsk}}] Let $a\ge0$ (and recall that $\1^0=\emptyword$): then, for any $y\in[\1^a\0]$
\begin{equation}\label{powernbis}
n\ge a+1\Longrightarrow\Pi_n(y,U)=\begin{pmatrix}\theta(a)\\\theta(a)\\1-2\theta(a)\end{pmatrix}
\quad\text{where}\quad \theta(a)={1\over 4}\left({1+2a\over 1+a}\right).
\end{equation}
The Cauchy condition (\ref{modifunif}) of Lemma~\ref{uniform} holds for $y\in[\1^{a_0}\0]$, since for $n\ge a_0+1$ and~$r\ge0$,
$$
z\in[y_1\cdots y_n]\Longrightarrow\Pi_{n+r}(z,U)-\Pi_n(z,U)=0.
$$
For the remaining case of $y=\bar\1$, take $z\in\YY$ s.t. $z_1\cdots z_n=\1\cdots\1$. Then,  either $z_1\cdots z_{n+r}=\1\cdots\1$ and we use (\ref{powern}) to write
\begin{align*}
\Vert\Pi_{n+r}(z,U)-\Pi_n(z,U)\Vert\le 2\left\Vert{2n+2\over 2n+3}\begin{pmatrix}1/2\\1/2\\1/(2n+2)\end{pmatrix}-\begin{pmatrix}1/2\\1/2\\0\end{pmatrix}\right\Vert,
\end{align*}
or $z_1\cdots z_{n+r}=\1^a\0w$, so that using (\ref{powernbis})
\begin{align*}
\Vert \Pi_{n+r}(z,U)-\Pi_n(z,U)\Vert
&\le \left\Vert\begin{pmatrix}\theta(n)\\\theta(n)\\1-2\theta(n)\end{pmatrix}-\begin{pmatrix}1/2\\1/2\\0\end{pmatrix}\right\Vert+\left\Vert{2n+2\over 2n+3}
\begin{pmatrix}1/2\\1/2\\1/(2n+2)\end{pmatrix}-\begin{pmatrix}1/2\\1/2\\0\end{pmatrix}\right\Vert\;;
\end{align*}
the Cauchy condition (\ref{modifunif}) of Lemma~\ref{uniform} is again satisfied about the ${\bf (C)}$-singular sequence $y=\bar\1$.

\end{proof}

The function $\VV$ being continuous, the value of $\VV(y)$ we have obtained in (\ref{1a00}) for $y=\1^{a_0}\0^{*}\1^{a_1}\0^{*}\1^{a_2}\cdots$ remains valid for any $y\in[\1^{a_0}\0]$. Finally, notice that for any $y\in\YY$ and $n\to+\infty$, one has
$$
\phi_n(y)=\log\left({\bar\eta_\yy[y_1\cdots y_n]\over \bar\eta_\yy[y_2\cdots y_n]}\right)\to \phi(y)=
\log\left(\frac13\cdot{U^\star A(y_1)\VV(\sigma\cdot y)\over U^\star\VV(\sigma\cdot y)}\right)=\log\Big(\frac13\cdot U^\star A(y_1)\VV(\sigma\cdot y)\Big).
$$

\begin{proposition}[see \cite{Oli10a}]
The Kamae measure $\bar\eta_\yy$ is an ergodic weak Gibbs measure of  $\phi:\{\0,\1\}^\NNN\to\mathbb{R}$ such that
$$
\phi(y)=
\begin{cases}
\log{1\over 3} & \text{if $y\in[{\0}^2]$ or $y=\bar\1$},\cr
\log{1\over 3}+\log\left({2a+1\over a+1}\right)&\text{if $y\in[\0{\1}^{a}\0]$, $a\ge 1$,}\cr
\log{1\over 3}+\log\left({a+1\over a}\right)&\text{if $y\in[{\1}^{a}\0]$, $a\ge 1$,}\cr
\log{2\over 3}&\text{if $y=\0\bar\1$.}
\end{cases}
$$

\end{proposition}

\begin{figure}[H]
  \begin{center}
       \includegraphics[scale=0.5]{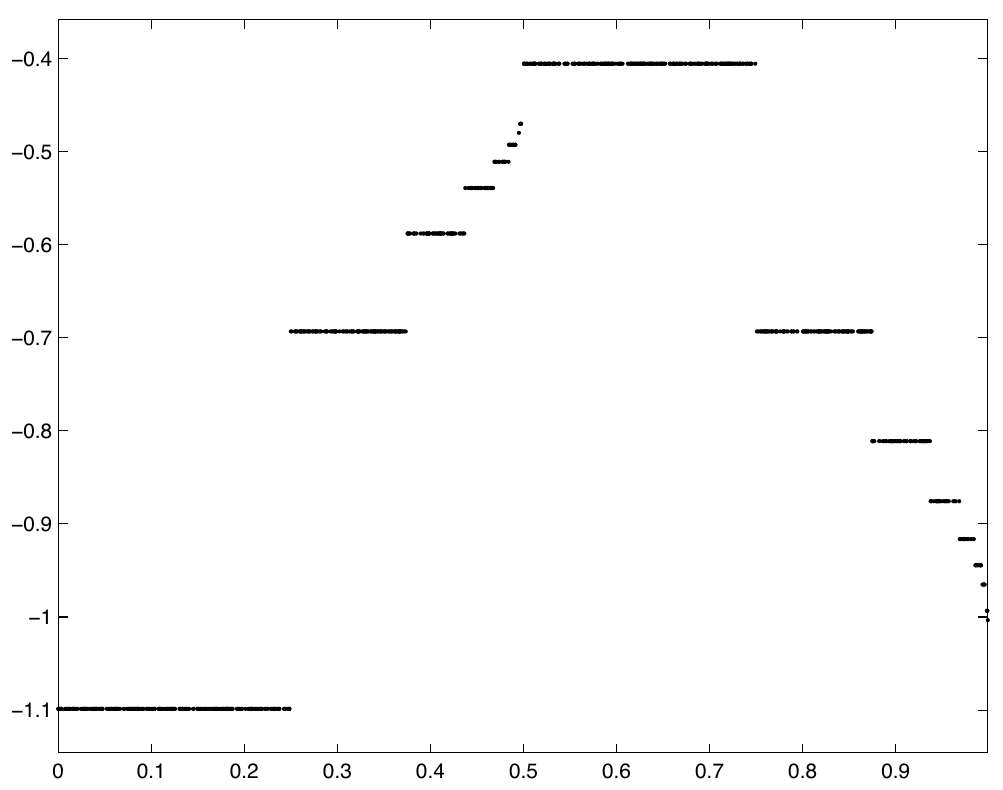}
\caption{\label{KAMAEPOT}\leftskip=10mm \rightskip=10mm
\small  \it Graph plot of $[0\,;1[\ni y\mapsto \phi(y_1y_2\cdots)$, where $0.y_1y_2\cdots$ is the dyadic expansion of $y$.} 
  \end{center}
\end{figure}

\section{\bf Gibbs structure of a special Bernoulli convolution}\label{ICBM}

\subsection{Bernoulli convolutions}Let $1<\beta<2$ be an arbitrary real number and 
define $\mu$ to be the distribution of the random variable $X:\{\0,\1 \}^\NNN\to\RRR$ such that
$$
X(\xi_1\xi_2\cdots) =(\beta-1)\sum_ {n = 1}^{+\infty} {\xi_n\over\beta^{n}}
$$
where $\{\0,\1 \}^\NNN$ is weighted by the Bernoulli probability $\PPP$ of parameter $ (p_0, p_1) $ (the normalization coefficient $(\beta-1)$ is introduced so that $0\le X(\xi_1\xi_2\cdots)\le1$). Equivalently,  $\mu$ may be defined as a {\it Bernoulli convolution}, that is an {\it Infinite Convolved Bernoulli Measure} (ICBM). Indeed, the law of the random variable $ X_n:\{\0,\1 \}^\NNN\to\RRR$ such that $ X_n (\xi) =(\beta-1)\xi_n/\beta^{n} $ is $p_0\cdot\delta_ {0} +p_1\cdot \delta_ {x_n}$, where $x_{n}=(\beta-1)/\beta^n$ (and $\delta_a$ is the Dirac mass concentrated on  $\{a\}$):  because $ X_1,X_2,\dots$ form a sequence of independent random variables, one gets:
\begin{equation*}\label{ICBMeq}
\mu:=
\mathop{\raise-2pt\hbox{\huge $*$}}_{n=1}^{+\infty}\Big(p_0\cdot \delta_{0}+p_1\cdot \delta_{x_n}\Big).
\end{equation*}

The present Section is based on previous works \cite{FO03,OST05,OT10,Oli12}, concerned with the Gibbs properties of Bernoulli convolutions when $\beta$ has specific algebraic properties. The framework of these papers  stands the arithmetic properties of the  Pisot-Vijayaraghavan (PV) numbers:  $\beta$ is a PV-Number if it is a  real algebraic integer with $\beta>1$, each ones of its Galois conjugates having a modulus stricly bounded by $1$. We shall develop here the approach in \cite{OST05} and \cite{OT10} used to deal with the multinacci numbers: the multinacci (number)  of order $n$ (for $n\ge2$) is the PV-number whose minimal polynomial is $X^n-(X^{n-1}+\cdots+X+1)$; the large Golden Number $\Phi$ is the multinacci  of order $2$, that is  $\Phi>1$ and $\Phi^2=\Phi+1$ --~the opposite of its Galois conjugate being the small Golden Number $1/\Phi$. The so-called {\it Erd\H os measure}  corresponds to the case of $\beta=\Phi$ and $(p_0,p_1)=(1/2,1/2)$: it is issued from {\it An Erd\H os Problem} dating back to the 1939 paper by Erd\H os \cite{Erd39} (see \cite{PSS00} for a review, historical notes and detailed references).

The question of the Gibbs nature of the Bernoulli convolutions was raised by Sidorov \& Vershik in \cite{SV98}: partial answers do  exist (for instance for the multinacci numbers \cite{OST05}\cite{OT10}) but, even for $\beta$ within the class of the PV numbers, a general approach is lacking. In other to enlighten the difficulties of the question, we shall consider a   PV-number --~not belonging to the multinacci family~-- that is $\beta\approx1.755$ s.t. 
\begin{equation}\label{minimalpol}
\beta^3=2\beta^2-\beta+1.
\end{equation}
In the subsequent analysis, the algebraic properties of the number $\beta$ are not apparent: actually they are implicit w.r.t. the fact that the vertex set (namely $\mathcal{V}$) of a certain adjacency graph (see Lemma~\ref{vertex}) is finite: let's mention the underlying result ensuring that property for a general PV-number is known as {\it Garsia separation Lemma} (see \cite{Gar63} and \cite[Proposition~2.6]{OST05}).

\subsection{$\beta$-shifts}For the sake of simple exposition, we  consider in this paragraph that $1<\beta<2$. 
The {\em Parry expansion} \cite{Par60} of the real number $x\in[0\,;1[$ in base $\beta$ is the unique sequence $\epsilon_1\epsilon_2\cdots\in\{\0,\1 \}^\NNN$ such that, for any $n\ge1$, 
\begin{equation}\label{expan}
s_n\le x <s_n+{1\over\beta^n}\quad\hbox{where}\quad s_n=\sum_{k=1}^n{\epsilon_k\over\beta^k}\;;
\end{equation}
therefore, it is licit to define $\varepsilon_k(x):=\epsilon_k$. The set $\Omega_\beta$ of the {\em$\beta$-admissible sequences} is the closure (in $\{\0,\1 \}^\mathbb{N}$) of the set of Parry expansions of real numbers in $[0\,;1[$. The set $\Omega_\beta$ is compact and left invariant by the full shift map $\sigma:\{\0,\1 \}^\mathbb{N}\to\{\0,\1 \}^\mathbb{N}$ in the sense that $\sigma(\Omega_\beta)=\Omega_\beta$: $\Omega_\beta$ is called the $\beta$-shift. The set $\Omega_\beta^{*}$ (resp. $\Omega_\beta^{(n)}$) of $\beta$-admissible words (resp. $\beta$-admissible words of length $n$) is made of finite (resp. of length $n$) factors  of $\beta$-admissible sequences (considered as one-side infinite words in $\Omega_\beta$).
We shall use the partition of $[0\,;1[$ by {\it $\beta$-adic intervals of order $n$}, that we write 
$[0\,;1[=\bigsqcup\nolimits_{w} I_w$,
where $w$ runs over $\Omega_\beta^{(n)}$ and for $\epsilon_1\dots\epsilon_n\in \Omega_\beta^{(n)}$, 
$$
I_{\epsilon_1\dots\epsilon_n}:=\Big\{x\in[0\,;1[\;;\;\varepsilon_1(x)\cdots\varepsilon_n(x)=\epsilon_1\dots\epsilon_n\Big\}.
$$

\subsection{A special case} We now assume that $\beta\approx1.755$ is  the PV-number defined by (\ref{minimalpol}) and $\mu$ stands for the uniform (i.e. $(p_0,p_1)=(1/2,1/2)$) Bernoulli convolution associated with this $\beta$.
\begin{lemma}\label{adm}The $\beta$-shift is the subshift of $\{\0,\1\}^\mathbb{N}$ determined by the finite type condition:
\begin{equation}\label{lhs}
\epsilon_1\epsilon_2\cdots\in\Omega_\beta\iff\Big(\epsilon_k\epsilon_{k+1}=\1\1\;\Longrightarrow\;\epsilon_{k+2}\epsilon_{k+3}=\0\0\Big)
\end{equation}
\end{lemma}
\begin{proof}[{\bf Proof}] Let $\epsilon=\epsilon_1\epsilon_2\cdots$ be the Parry expansion of a given real number $x$; the condition (\ref{expan}) is equivalent to $\beta^n(x-s_n)=\beta^n(x-s_{n-1})-\epsilon_n\in[0\,;1[$, hence $\epsilon_n=\lfloor\beta^n(x-s_{n-1})\rfloor$ -- with the convention that $s_0=0$. Now (\ref{expan}) implies
$$
{\epsilon_{n+1}\over\beta^{n+1}}+{\epsilon_{n+2}\over\beta^{n+2}}+{\epsilon_{n+3}\over\beta^{n+3}}+{\epsilon_{n+4}\over\beta^{n+4}}\le x-s_n<{\1\over\beta^n}={\1\over\beta^{n+1}}+{\1\over\beta^{n+2}}+{\1\over\beta^{n+4}}
$$
(the last equality being equivalent to (\ref{minimalpol})); hence, the condition $\epsilon_{n+1}\epsilon_{n+2}=\1\1$ with $\epsilon_{n+3}\epsilon_{n+4}\ne \0\0$ leads to an impossible strict inequality in any cases.

Conversely, let $\epsilon=\epsilon_1\epsilon_2\cdots\in\{\0,\1\}^\mathbb{N}$ be a sequence satisfying the r.h.s. conditions in (\ref{lhs}); for any $n\ge0$,  the sequence $\sigma^n\cdot\epsilon$ is either equal to $\overline{\1\1\0\0}$ or  begins by either $(\1\1\0\0)^i\0$ or $(\1\1\0\0)^i\1\0$ for some integer $i\ge0$ (recall that given $w$ any word, $w^k:=w\cdots w$ with $k$ factors and  $w^0$ stands for the empty word $\emptyword$). For $m> n$ the real numbers  
$s_m=\sum_{k=1}^m{\epsilon_k/\beta^k}$ and $s_n=\sum_{k=1}^n{\epsilon_k/\beta^k}
$ satisfy
$0\le s_m-s_n<1/\beta^{n}\sum_{k=1}^{+\infty}{\xi_k/\beta^k}$
where 
$$
\xi_1\xi_2\cdots\in\Big\{\overline{\1\1\0\0},(\1\1\0\0)^i\0\bar1,(\1\1\0\0)^i\1\0\bar\1\Big\}.
$$
A simple computation yields $\displaystyle0\le s_m-s_n<{1/\beta^n}$, so that $\epsilon_1\dots\epsilon_m\bar\0$ is the Parry expansion of $s_m$; the condition $m>n$ being arbitrary means that  $\epsilon=\epsilon_1\epsilon_2\cdots$ is $\beta$-admissible.

\end{proof}

Lemma~\ref{adm} implies that each $\beta$-admissible sequence $\epsilon=\epsilon_1\epsilon_2\cdots$ may be decomposed from left to right in an infinite concatenation of words, say
$\epsilon=w_1w_2\cdots$, where 
\begin{equation}\label{hgfkdsjhfgk}
w_i\in\Big\{\ww(\0)=\0,\ww(\1)=\1\0, \ww(\2)=\1\1\0\0\Big\}.
\end{equation}
\begin{proposition}\label{threeintervals}Let $J_\emptyword:=[0\,;1[$ and for any $\xi_1\cdots\xi_n\in\{\0,\1,\2\}^n$, the concatened word $\ww(\xi_1)\cdots\ww(\xi_n)$ being a binary $\beta$-admisible word, define
$$
J_{\xi_1\cdots\xi_n}:=I_{\ww(\xi_1)\cdots\ww(\xi_n)}\;;
$$
then, for any word $w\in\{\0,\1,\2\}^*$ one has the partition
$$
J_w=J_{w\0}\sqcup J_{w\1}\sqcup J_{w\2}.
$$
\end{proposition}

The two subsequent lemmas give a way to compute the $\mu$-measure of the intervals $J_w$, for any $w\in\{\0,\1,\2\}^*$ by means of the matrices 
\begin{equation}\label{matrix7X7bis}
\begin{matrix}
A(\0):=\begin{pmatrix}1&0&0&0&0&0&0\\0&0&1&0&0&0&0\\0&0&0&1&1&0&0\\0&0&0&0&0&0&0\\1&0&0&0&0&0&1\\0&0&0&0&1&0&0\\0&1&0&0&0&0&0\end{pmatrix}&
A(\1):=\begin{pmatrix}0&0&1&1&0&0&0\\0&0&0&0&0&1&0\\0&0&0&1&1&0&0\\1&0&0&0&0&0&0\\0&0&1&0&0&0&0\\0&0&0&0&0&0&0\\0&0&0&0&0&0&0\end{pmatrix}\\\\

A(\2):=\begin{pmatrix}1&0&0&0&1&0&1\\0&0&0&0&0&0&0\\1&0&0&0&0&0&1\\0&0&0&1&1&0&0\\0&0&0&0&1&0&0\\0&0&0&0&0&0&0\\0&0&0&0&0&0&0\end{pmatrix}&
\end{matrix}
\end{equation}
(considered in Section~\ref{7X7} in view of an application of Theorem~A).
\begin{lemma}\label{vertex}
Consider the vertex set 
\begin{align*}
\mathcal{V}=\Big\{v_1=0,\;v_2=1,&\;v_3=1-(\beta-1)^2,\;
\\&v_4=-(\beta-1)^2,\;v_5=\beta-1,\;v_6=\beta-(\beta-1)^2,\;v_7=\beta(\beta-1)\Big\}
\end{align*}
together with the matrices 
\begin{align*}
M(\0):={1\over2}A(\0),\quad 
M(\1):={1\over4}A(\1),\quad M(\2):={1\over 16}A(\2)\quad\text{and}\quad R:=\begin{pmatrix}3/5\\2/5\\13/20\\1/5\\3/5\\3/10\\1/5\end{pmatrix}.
\end{align*}
Then, for any $n\ge1$ and any word  $\xi_1\cdots\xi_n\in\{\0,\1,\2\}^n$
\begin{equation}\label{translatedcylinders}
\begin{pmatrix}\mu\left({1\over\beta}\Big(v_1+J_{\xi_1\cdots \xi_n}\Big)\right)\\\vdots\\\mu\left({1\over\beta}\Big(v_7+J_{\xi_1\cdots \xi_n}\Big)\right)\end{pmatrix}=M(\xi_1\cdots \xi_n)R.
\end{equation}
\end{lemma}
\begin{proof}[{\bf Proof}] For any real number $\gamma$ we evaluate 
$$
\displaystyle\mu\left({1\over \beta}\Big(\gamma+J_{\xi_1\cdots \xi_n}\Big)\right)
$$ 
in the three cases when $\xi_1$ is either $\0,\1$ or $\2$. First, consider that $\xi_1=\0$, so that
$$
\displaystyle(\beta-1)\sum_{k=1}^{+\infty}{\omega_k\over\beta^k}\in\displaystyle{1\over\beta}\Big(\gamma+J_{\xi_1\cdots \xi_n}\Big)\iff\displaystyle(\beta-1)\sum_{k=1}^{+\infty}{\omega_{k+1}\over\beta^k}\in{1\over\beta}\Big(\gamma'+J_{\xi_1\cdots \xi_n}\Big)
$$ 
with $\gamma'=\gamma\beta-\omega_1\beta(\beta-1)$ and since $\omega_1\in\{\0,\1 \}$,
\begin{equation}\label{firstcase}
\mu\left({1\over\beta}\Big(\gamma+J_{\xi_1\cdots \xi_n}\Big)\right)={1\over2}\sum_{\gamma'\in\mathcal{V}_0(\gamma)}\mu\left({1\over\beta}\Big(\gamma'+J_{\xi_1\cdots \xi_n}\Big)\right),
\end{equation}
where
$$
\mathcal{V}_{\0}(\gamma):=\Big\{\gamma\beta-\uuu\beta(\beta-1)\;;\;\uuu\in\{\0,\1 \}\Big\}.
$$
We proceed in the same way for $\xi_1=\1$:
\begin{equation}\label{secondcase}
\mu\left({1\over\beta}\Big(\gamma+J_{\xi_1\cdots \xi_n}\Big)\right)={1\over4}\sum_{\gamma'\in\mathcal{V}_1(\gamma)}\mu\left({1\over\beta}\Big(\gamma'+J_{\xi_1\cdots \xi_n}\Big)\right),
\end{equation}
where
$$
\mathcal{V}_{\1}(\gamma):=\Big\{\gamma\beta^2+\beta-(\uuu\beta+\vvv)\beta(\beta-1)\;;\;\uuu,\vvv\in\{\0,\1 \}\Big\}.
$$
Finally, in the case of  $\xi_1=\2$, one gets
\begin{equation}\label{thirdcase}
\mu\left({1\over\beta}\Big(\gamma+J_{\xi_1\cdots \xi_n}\Big)\right)={1\over16}\sum_{\gamma'\in\mathcal{V}_2(\gamma)}\mu\left({1\over\beta}\Big(\gamma'+J_{\xi_1\cdots \xi_n}\Big)\right),
\end{equation}
where
$$
\mathcal{V}_\2(\gamma)=\Big\{\gamma\beta^4+\beta^3+\beta^2-(\uuu\beta^3+\vvv\beta^2+\www\beta+\xxx)\beta(\beta-1)\;;\;\uuu,\vvv,\www,\xxx\in\{\0,\1 \}\Big\}.
$$
However, the measure $\mu$ being supported by the unit interval $[0\,;1]$, one has:

\begin{lemma}\label{truc}
$\gamma\not\in]-1\,;\beta[\;\Longrightarrow\;\displaystyle\mu\left({1\over \beta}\Big(\gamma+[0\,;1[\Big)\right)=0$.
\end{lemma}

Because of Lemma~\ref{truc}, we define for each $\eee\in\{\0,\1,\2\}$ the relation
$$
\gamma\mathop{\longrightarrow}^\eee\gamma'\iff\gamma'\in\mathcal{V}_\eee(\gamma)\cap]-1\,;\beta[
$$
and by definition, $\mathcal{V}$ is the subset of $]-1\,;\beta[$ containing $\0$ and the $\gamma\in ]-1\,;\beta[$ for which there exists a sequence $\eee_0,\dots,\eee_n$ in $\{\0,\1,\2\}$ and $0=\gamma_0,\cdots,\gamma_n=\gamma$ in $]-1\,;\beta[$  (with $n\ge0$) s.t.
$$
\gamma_0\mathop{\longrightarrow}^{\eee_0}\gamma_1,\dots,\gamma_{n-1}\mathop{\longrightarrow}^{\eee_n}\gamma_n\;;
$$
then, a direct computation gives 
\begin{align*}
\mathcal{V}=\Big\{v_1=0,\;v_2=1,&\;v_3=1-(\beta-1)^2,\;
\\&v_4=-(\beta-1)^2,\;v_5=\beta-1,\;v_6=\beta-(\beta-1)^2,\;v_7=\beta(\beta-1)\Big\}.
\end{align*}
For $\eee\in\{\0,\1,\2\}$, we define $A'(\eee)$ to be the $7\times 7$ incident matrix  s.t. for any $1\le i,j\le 7$, the $(i,j)$ entry of 
$A'(\eee)$ is $1$ if and only if $v_i\mathop{\longrightarrow}\limits^{\eee}v_j$ (the corresponding incidence relations are represented by the adjacency graph in Figure~\ref{GRAPH3}).

\begin{figure}[H]
\begin{center}
       \includegraphics[scale=0.34]{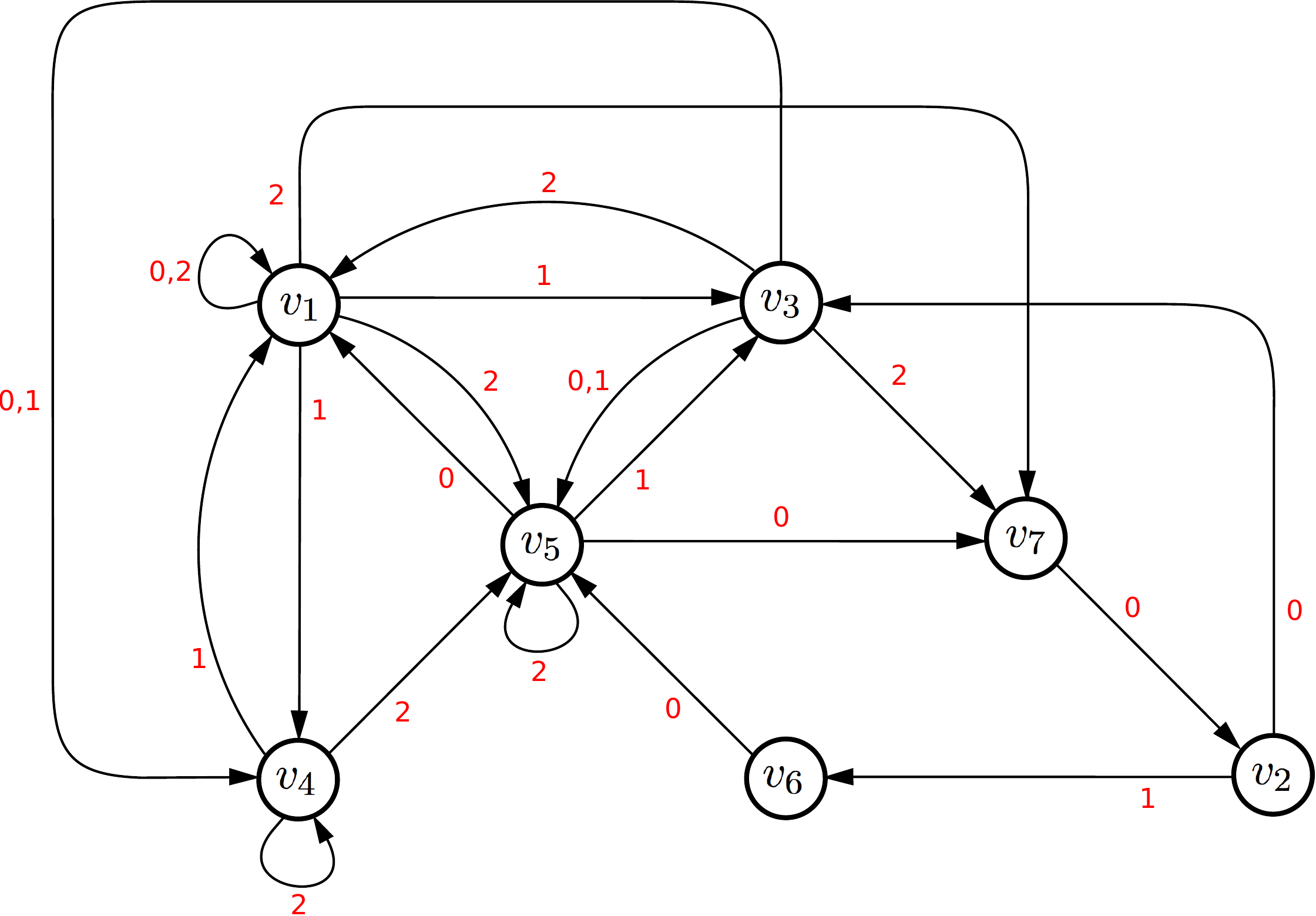}
\caption{\label{GRAPH3}\footnotesize \it Adjacency graph corresponding to the relation $v_i\mathop{\longrightarrow}\limits^\ee v_j$ ($\ee=\0,\1,\2$) between the elements in $\mathcal{V}=\{v_1,\dots,v_7\}$, which are the $7$ vertices of the graph.} 
\end{center}
\end{figure}

The key point is that 
$$
A'(\0)=A(\0),\quad A'(\1)=A(\1),\quad 
A'(\2)=A(\2)
$$
and  one deduces from (\ref{firstcase}), (\ref{secondcase}) and (\ref{thirdcase}) that for any $\xi\in\{\0,\1,\2\}^\NNN$ and any $n\ge1$
$$
\begin{pmatrix}\mu\left({1\over\beta}\Big(v_1+J_{\xi_1\cdots \xi_n}\Big)\right)\\\vdots\\\mu\left({1\over\beta}\Big(v_7+J_{\xi_1\cdots \xi_n}\Big)\right)\end{pmatrix}=M(\xi_1\cdots \xi_n)R,\quad\hbox{where }R=\begin{pmatrix}\mu\left({1\over\beta}\Big(v_1+[0\,;1[\Big)\right)\\\vdots\\\mu\left({1\over\beta}\Big(v_7+[0\,;1[\Big)\right)\end{pmatrix}.
$$
In order to compute $R$ we consider the case $n=1$ and we sum for $\xi_1\in\{\0,\1,\2\}$. Since (Proposition \ref{threeintervals}) the sets $J_{\xi_1}$ form a partition of $[0\,;1[$, we obtain that $R$ is an eigenvector of $M(\0)+M(\1)+M(\2)$. Moreover, because the sum of the two first entries in $R$ is 
$$
\mu\Big(0+{1/\beta}[0\,;1[\Big)+\mu\Big(1/\beta+{1/\beta}[0\,;1[\Big)=\mu\Big({1/\beta}[0\,;2[\Big)=\mu([0\,;1])=1,
$$
the computation of this eigenvector gives the expected value for~$R$.

\end{proof}

A direct consequence is the exact value of  the $\mu$-measure of the intervals $J_w$.
\begin{lemma}\label{cylinders}For any $\xi\in\{\0,\1,\2\}^\NNN$ and any $n\ge1$,
\begin{equation}\label{*values}
\mu\left(J_{\xi_1\cdots \xi_n}\right)=\left\{\begin{array}{ll}{U^\star_1}M(\xi_2\cdots \xi_n)R&\hbox{if }\xi_1=\0\\
{U^\star_2}M(\0)M(\xi_2\cdots \xi_n)R&\hbox{if }\xi_1=\1\\
{U^\star_2}M(\1)M(\0)M(\xi_2\cdots\xi_n)R&\hbox{if }\xi_1=\2.\end{array}\right.
\end{equation}
\end{lemma}
It remains to prove that $\mu$ is a weak-Gibbs measure.

\begin{theorem}\label{Gibbs}
$\mu$ is weak-Gibbs in the sense that there exists a continuous  $\Psi:\{\0,\1,\2\}^\NNN\to\RRR$ and a sub-exponential sequence of positive constants $K_1,K_2,\dots$ s.t for any $\xi\in \{\0,\1,\2\}^\NNN$, 
\begin{equation}\label{A*}
{1\over K_n}\leq{\displaystyle\mu\left(J_{\xi_1\cdots\xi_n}\right)\over
\exp\Big(\sum_{k=0}^{n-1} \Psi(\sigma^k\cdot \xi)\Big)}\leq K_n.
\end{equation}
\end{theorem}
We shall need  the following theorem whose proof (see Section~\ref{7X7}) depends on Theorem~A.
\begin{theorem}\label{3matrices}Let $A(\0),A(\1),A(\2)$ be the $7\times7$ matrices in (\ref{matrix7X7bis}): if the vector $V$ has positive entries  then, the sequence of the probability  vectors
$$
\Pi_n(\omega,V)={A(\omega_1)\cdots A(\omega_n)V\over \Vert A(\omega_1)\cdots A(\omega_n)V\Vert},
$$
converges to a probability vector $\VV(\omega)$, uniformly over  $\omega\in\{\0,\1,\2\}^\NNN$, as $n\to+\infty$; moreover, 
$$
\Delta(\VV(\omega))\ge U_1+U_3+U_5,
$$
except for $\omega\in\{\bar\0=\0\0\cdots,\bar\2=\2\2\cdots\}$ for which
$$
\VV(\bar\0)=1/5(U_2+U_3+U_5+U_6+U_7)\quad\text{and}\quad
\VV(\bar\2)=1/3(U_1+U_3+U_4).
$$

\end{theorem}

\begin{proof}[{\bf Proof of Theorem~\ref{Gibbs}}] In view of Proposition~\ref{symbweakGibbs}, we stress that the uniform convergence of the $n$-step potential of $\mu$ does not hold: to overcome this difficulty we shall introduce an intermediate measure. Let $U=U_1+\cdots+U_7$; by  Kolmogorov Extension  Theorem, there exists a unique measure $\nu$ on $[0\,;1]$ such that, for any $\xi_1,\dots,\xi_n\in{{\mathcal W}_0}$,
\begin{equation}\label{umeas}
\nu\left(J_{\xi_1\cdots \xi_n}\right)=\left\Vert M(\xi_1\cdots \xi_n)R\right\Vert=U^\star M(\xi_1\cdots \xi_n)R.
\end{equation}
Recall that for any $\xi\in\{\0,\1,\2\}^\NNN$ and any $n\ge1$
$$
\Pi_n(\xi,R)={M(\xi_1\cdots\xi_n)R\over \Vert M(\xi_1\cdots\xi_n)R\Vert}={A(\xi_1\cdots\xi_n)R\over \Vert A(\xi_1\cdots\xi_n)R\Vert},
$$
so that the  $n$-step potential of the measure $\nu$,  that is $\phi_n:\{\0,\1,\2\}^\NNN\to\RRR$ s.t.
$$
\phi_n(\xi)=\log\left({U^\star M(\xi_1) \Pi_{n-1}(\sigma\cdot \xi,R)\over U^\star \Pi_{n-1}(\sigma\cdot \xi,R)}\right),
$$
converges uniformly on $\{\0,\1,\2\}^\mathbb N$ toward a continuous  $\Psi:\{\0,\1,\2\}^\NNN\to\RRR$: indeed, by Theorem~\ref{3matrices} the  functions $\xi\mapsto \Pi_n(\xi,R)$ form a sequence which converges uniformly for $\xi\in\{\0,\1,\2\}^\NNN$ toward $\VV(\xi)\in\mathcal{S}_7$ s.t. $M(\aaa)\VV(\xi)\ne0$, for any $\xi\in\{\0,\1,\2\}^\NNN$ and $\aaa\in\{\0,\1,\2\}$.
Therefore (Proposition~\ref{symbweakbis}), the measure $\nu$ is weak Gibbs w.r.t. $\Psi$ in the sense that there exists a sub-exponential sequence of positive constants $K_1,K_2,\dots$ s.t for any $\xi\in \{\0,\1,\2\}^\NNN$, 
\begin{equation*}
{1\over K_n}\leq{\displaystyle\nu\left(J_{\xi_1\cdots\xi_n}\right)\over
\exp\Big(\sum_{k=0}^{n-1} \Psi(\sigma^k\cdot \xi)\Big)}\leq K_n.
\end{equation*}
Hence, the weak Gibbs property of $\mu$ (w.r.t. to $\Psi$) will be established if one is able to show that
\begin{equation}\label{ratiospower1/k}
\displaystyle\lim_{n\to+\infty}\left({\nu\left(J_{\xi_1\cdots\xi_n}\right)\over\mu\left(J_{\xi_1\cdots\xi_n}\right)}\right)^{1/n}=1
\end{equation} 
uniformly  for $\xi=\xi_1\xi_2\cdots\in\{\0,\1,\2\}^\NNN$. First, from Lemma~\ref{cylinders} together with the inequalities 
\begin{align*}
U^\star M(\0)\ge{U^\star_1},\quad 
U^\star M(\1)\ge{U^\star_2}M(\0),\quad 
U^\star M(\2)\ge{U^\star_2}M(\1)M(\0),
\end{align*}
one deduces that for any $\xi\in\{\0,\1,\2\}^\NNN$ and any $n\ge1$, 
\begin{equation}\label{stroumf0}
1\le {\nu(J_{\xi_1\cdots\xi_n})\over \mu(J_{\xi_1\cdots\xi_n})}.
\end{equation}
Concerning the upper bound we consider the different cases for $\xi_1$ either equal to $\0,\1$ or $\2$.

$\bullet$ Suppose that  $\xi_1\in[\1]$. From (\ref{umeas}) and (\ref{*values}) and the fact that $U^\star_2M(\0)=1/2\cdot U^\star_3$ one has:
$$
{\nu\left(J_{\xi_1\cdots\xi_n}\right)\over\mu\left(J_{\xi_1\cdots\xi_n}\right)}={U^\star M(\1)\Pi_{n-1}(\sigma\cdot \xi,R)\over U^\star_2M(\0)\Pi_{n-1}(\sigma\cdot\xi,R)}\le 
{2\cdot \Vert M(\1)\Vert \over U^\star_3\Pi_{n-1}(\sigma\cdot\xi,R)}={2\cdot \Vert M(\1)\Vert \over F_{n-1}(\sigma\cdot\xi)}
$$
where $F_k:\{0,\1,\2\}^\NNN\to\RRR$ is the continuous function s.t. $F_k(\zeta)=U^\star_3\Pi_{k}(\zeta,R)$ (and the convention $\Pi_0(\cdot,R)=R/\Vert R\Vert$). By Theorem~\ref{3matrices}, $\Pi_k(\zeta,R)$ tends uniformly (as $k\to+\infty$) to a probability vector $\VV(\zeta)$ for any $\zeta\in\{\0,\1,\2\}^\mathbb N$: hence, the $F_k$ are continuous positive functions converging uniformly (on $\{0,\1,\2\}^\NNN$) toward $F:\{0,\1,\2\}^\NNN\to\RRR$ s.t. $F(\zeta)=U_3^\star\VV(\zeta)$;  however, (Theorem~\ref{3matrices} again)  $\Delta(\VV(\zeta))\ge U_3$ ensuring the continuous function $F$ to be also positive:  in particular the infimum $\inf\{F_*\}$ of the $F_k(\zeta)$ for $\zeta\in\{\0,\1,\2\}^\mathbb N$ and $k\ge0$ is positive. Finally one deduces 
\begin{equation}\label{bothratios}
\xi\in[\1]\Longrightarrow{\nu\left(J_{\xi_1\cdots\xi_n}\right)\over\mu\left(J_{\xi_1\cdots\xi_n}\right)}\le {2\cdot \Vert M(\1)\Vert \over \inf\{F_*\}}.
\end{equation}

$\bullet$  Suppose that $\xi\in[\0]$; the functions $G_k:[\1]\cup[\2]\to\RRR$ such that $G_k(\zeta)=U_1^\star\Pi_k(\zeta,R)$ are continuous, positive and converge uniformly (on $[\1]\cup[2]$) toward $G:[\1]\cup[2]\to\RRR$ s.t. $G(\zeta)=U_1^\star\VV(\zeta)$: by Theorem~\ref{3matrices}, $G$ is continuous and positive:  in particular the infimum $\inf\{G_*\}$ of the $G_k(\zeta)$ for $\zeta\in[\1]\cup[\2]$ and $k\ge0$ is positive. Because $\xi_1=\0$, we consider the largest $1\le m\le n$ for which $\0^m$ is a prefix of $\xi_1\dots\xi_n$ so that
$\xi_1\dots\xi_k=\0^m\xi_{m+1}\dots\xi_n\quad\hbox{where }\xi_{m+1}\ne\0\hbox{ if }m<n$.
From (\ref{umeas}) and (\ref{*values})  and the facts that 
$$
U^\star_1M(\0)^{m-1}={1\over 2^{m-1}}\cdot U^\star_1
\quad\text{and}\quad
U^\star M(\0)^m\le\frac{2m}{2^m}\cdot U^\star\le\frac n{2^{m-1}}\cdot U^\star
$$
one deduces 
\begin{equation*}\label{ratiosbis}
{\nu\left(J_{\xi_1\cdots\xi_n}\right)\over\mu\left(J_{\xi_1\cdots\xi_n}\right)}=\frac{U^\star M(\0)^m \Pi_{n-m}(\sigma^m\cdot\xi,R)}{U^\star_1M(\0)^{m-1}\Pi_{n-m}(\sigma^m\cdot \xi,R)}\le 
{n\cdot U^\star \Pi_{n-m}(\sigma^m\cdot\xi,R)\over U_1^\star\Pi_{n-m}(\sigma^m\cdot\xi,R)}.
\end{equation*}
Because  $\Vert \Pi_{n-m}(\sigma^m\cdot\xi,R)\Vert=1$ and the definition of both $G_{n-m}$ and $\inf\{G_*\}$
\begin{equation}\label{stroumf1}
\xi\in[\0]\Longrightarrow{\nu\left(J_{\xi_1\cdots\xi_n}\right)\over\mu\left(J_{\xi_1\cdots\xi_n}\right)}\le 
\begin{cases}
\displaystyle{n\cdot U^\star R\over U_1^\star R}&\text{if $m=n$}\\\\
\displaystyle{n\cdot U^\star \Pi_{n-m}(\sigma^m\cdot\xi,R)\over G_{n-m}(\sigma^m\cdot\xi)}\le {n\over \inf\{G_*\}}
&\text{if $m<n$}
\end{cases}
\end{equation}

$\bullet$  Suppose that $\xi\in[\2]$; the functions $H_k:[\1]\cup[\0]\to\RRR$ such that $H_k(\xi)=U_5^\star\Pi_k(\xi,R)$ are continuous, positive and converge uniformly (on $[\1]\cup[\0]$) toward $H:[\1]\cup[\0]\to\RRR$ s.t. $H(\xi)=U_5^\star\VV(\xi)$: by Theorem~\ref{3matrices}, $H$ is continuous and positive and  the infimum $\inf\{H_*\}$ of the $H_k(\xi)$ for $\xi\in[\1]\cup[\0]$ and $k\ge0$ is positive. Because $\xi_1=\2$, we consider the largest $1\le m\le n$ for which $\2^m$ is a prefix of $\xi_1\dots\xi_n$, so that $\xi_1\cdots\xi_n=\2^m\xi_{m+1}\dots\xi_n\quad\hbox{where }\xi_{m+1}\ne\2\hbox{ if }m<n$.
From (\ref{umeas}) and (\ref{*values}) together with
$$
\displaystyle U^\star M(\2)^m\le\frac{3m}{2^{4m}}\cdot U^\star\le\frac{3n}{2^{4m}}\cdot U^\star
\quad\text{and}\quad\displaystyle U^\star_2M(\1)M(\0)M(\2)^{m-1}=\frac1{2^{4m-1}}\cdot U^\star_5,
$$
one deduces
$$
{\nu\left(J_{\xi_1\cdots\xi_n}\right)\over\mu\left(J_{\xi_1\cdots\xi_n}\right)}=
\frac{U^\star M(\2)^m \Pi_{n-m}(\sigma^m\cdot\xi,R)}{U^\star_2M(\1)M(\0)M(\2)^{m-1}\Pi_{n-m}(\sigma^m\cdot\xi,R)}
\le {3n\over 2}\cdot {U^\star\Pi_{n-m}(\sigma^m\cdot\xi,R)\over U_5^\star\Pi_{n-m}(\sigma^m\cdot\xi,R)},
$$
that is, with  $\Vert \Pi_{n-m}(\sigma^m\cdot\xi,R)\Vert=1$ and the definition of both $H_{n-m}$ and $\inf\{H_*\}$
\begin{equation}\label{ratioster}
\xi\in[\2]\Longrightarrow{\nu\left(J_{\xi_1\cdots\xi_n}\right)\over\mu\left(J_{\xi_1\cdots\xi_n}\right)}\le 
\begin{cases}
\displaystyle{3n\cdot U^\star R\over 2\cdot U_5^\star R}&\text{if $m=n$}\\\\
\displaystyle{3n\cdot U^\star \Pi_{n-m}(\sigma^m\cdot\xi,R)\over 2\cdot H_{n-m}(\sigma^m\cdot\xi)}\le {3n\over 2\cdot\inf\{H_*\}}
&\text{if $m<n$}
\end{cases}
\end{equation}
It follows from (\ref{stroumf0})-(\ref{bothratios})-(\ref{stroumf1})-(\ref{ratioster}) that  (\ref{ratiospower1/k}) holds, proving the theorem.

\end{proof}

\section{\bf An advanced application of Theorem~A: proof of Theorem~\ref{3matrices}}\label{7X7}

This section is devoted to the full proof of Theorem~\ref{3matrices} stated and used in the previous section. The main problem is concerned with the application of Theorem~A w.r.t. the matrix products of the three $7\times 7$ matrices, $A(\0),A(\1),A(\2)$ given in (\ref{matrix7X7bis}).
In view of  Theorem~A, the major ingredient we shall need is the following proposition.
\begin{proposition}\label{hyp}The sequences  $\omega\in\{\0,\1,\2\}^\NNN$ which are ${\bf (C)}$-regular  --~in the sense that $\mathcal{A}(\omega)=(A(\omega_1),A(\omega_2),\dots)$ satisfies condition ${\bf (C)}$~-- form the dense set 
\begin{equation}\label{defXX}
\XX:=
\left(\{\0,\1,\2\}^\NNN\setminus\bigcup_{n=0}^{+\infty}\sigma^{-n}\{\bar\0\}\right)\cap\left(\{\0,\1,\2\}^\NNN\setminus\bigcup_{n=0}^{+\infty}\sigma^{-n}\{\bar\2\}\right).
\end{equation}
\end{proposition}
The proof of Proposition~\ref{hyp}  depends on the introduction of a language $\mathcal{W}\subset\{\0,\1,\2\}^*$.

\subsection{The language $\mathcal{W}$}\label{W}
The initial idea is to find (if possible) a language $\mathcal{W}\subset\{\0,\1,\2\}^*$ (i.e. each $W\in\mathcal{W}$ is a word with digits in $\{\0,\1,\2\}$ and whose length $|W|$ is finite) with two properties as follows. First (see Lemma~\ref{split} for the exact statement), $\mathcal{W}$ generates the sequences in $\XX$ (i.e. the $\omega\in\{\0,\1,\2\}^\NNN$ for which $\sigma^k\cdot \omega$ never belong to $\{\bar\0,\bar\2\}$)   
in the sense that such a  $\omega$ may be written as a concatenation of words, that is 
$\omega=W_{0}W_{1}W_{2}\dots$,
where $W_{1}W_{2}\dots$ are in $\mathcal{W}$ and $W_{0}$ is a (possibly empty) strict suffix of a word in $\mathcal{W}$. The second condition is the existence of $0\le\lambda_0<1\le\Lambda_0<+\infty$ such that $A(W)\in\mathcal{H}_2(\Lambda_0)$ for any $W\in\mathcal{W}$ (Lemma~\ref{H2}), and $A(W)\in\mathcal{H}_1\cap\mathcal{H}_3(\lambda_0)$ when $W$ is a concatenation of a (fixed) number of words in $\mathcal{W}$ (Lemma~\ref{H3}). To construct $\mathcal{W}$, notice that for $\sss$ either equal to $\0$ or $\2$, the condition that $A(\sss^k)=A(\sss)^k\in\mathcal{H}_2(\Lambda)$ for any $k$ implies $\Lambda=+\infty$ (this is the reason why $\bar\0$ and $\bar\2$ are ${\bf (C)}$-singular);
taking this remark into account, we begin to fix $\1\in\mathcal{W}$, the other words $W\in\mathcal{W}$ being of the form $W'\sss^k$, for some words $W'\ne\emptyword$ specified later, while $\sss=\0$ or $\2$, and for any $k\ge 1$: 
more precisely,  the (non empty) prefixes $W'$ of the words $W\in\mathcal{W}$ will be determined in order to ensure that any finite word not ending by $\1$ has a suffix of the form $W'\sss^k$, or $W''\sss^k$ with $W''$ suffix of $W'$, while $A(W)=A(W'\sss^k)$ belongs to $\mathcal{H}_2(\Lambda_0)$, for a $\Lambda_0<+\infty$ independent of $W$ (the determination of $\lambda_0$ being obtained in a second step).
\begin{figure}[H]
\begin{center}
\includegraphics[scale=0.5]{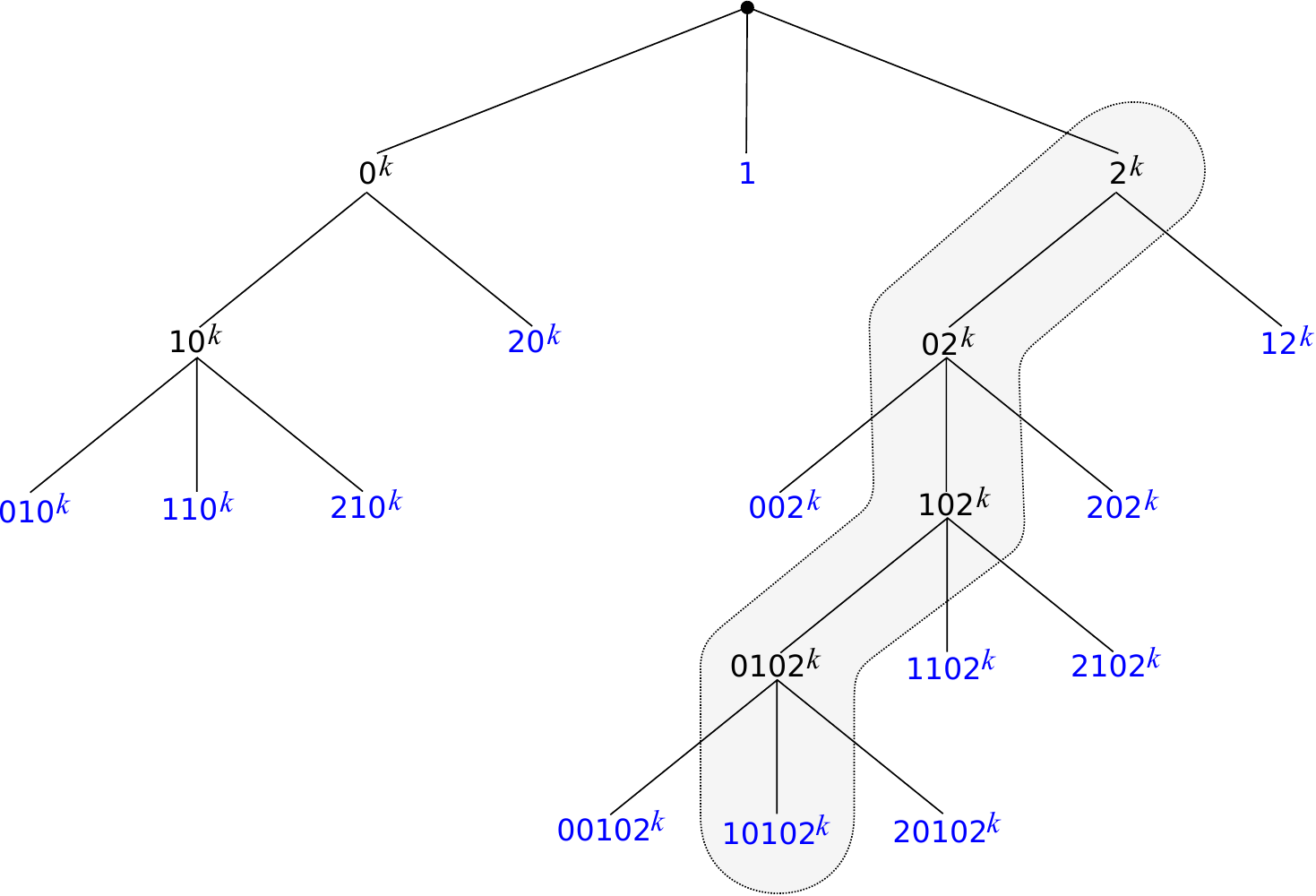}
\caption{\label{ARBRE}\footnotesize \it Lexicographical construction (from right to left) of the words in $\mathcal{W}$ (here $k\ge1$).} 
\end{center}
\end{figure}
For instance, let's compute successively $A(\2^k),A(\0\2^k),A(\1\0\2^k),A(\0\1\0\2^k),\dots$, which gives 
\begin{align*}
 A(\2^k)=&
\begin{pmatrix}1&0&0&0&k&0&1\\0&0&0&0&0&0&0\\1&0&0&0&{k-1}&0&1\\0&0&0&1&k&0&0\\0&0&0&0&1&0&0\\0&0&0&0&0&0&0\\0&0&0&0&0&0&0\end{pmatrix}\qquad
 A(\0\2^k)=
\begin{pmatrix}1&0&0&0&k&0&1\\1&0&0&0&{k-1}&0&1\\0&0&0&1&{k+1}&0&0\\0&0&0&0&0&0&0\\1&0&0&0&k&0&1\\0&0&0&0&1&0&0\\0&0&0&0&0&0&0\end{pmatrix}\\
 A(\1\0\2^k)=&
\begin{pmatrix}0&0&0&1&{k+1}&0&0\\0&0&0&0&1&0&0\\1&0&0&0&k&0&1\\1&0&0&0&k&0&1\\0&0&0&1&{k+1}&0&0\\0&0&0&0&0&0&0\\0&0&0&0&0&0&0\end{pmatrix}\quad
 A(\0\1\0\2^k)=
\begin{pmatrix}0&0&0&1&{k+1}&0&0\\1&0&0&0&k&0&1\\1&0&0&1&{2k+1}&0&1\\0&0&0&0&0&0&0\\0&0&0&1&{k+1}&0&0\\0&0&0&1&{k+1}&0&0\\0&0&0&0&1&0&0\end{pmatrix}\\
 A(\1\0\1\0\2^k)=&
\begin{pmatrix}1&0&0&1&{2k+1}&0&1\\0&0&0&1&{k+1}&0&0\\0&0&0&1&{k+1}&0&0\\0&0&0&1&{k+1}&0&0\\1&0&0&1&{2k+1}&0&1\\0&0&0&0&0&0&0\\0&0&0&0&0&0&0\end{pmatrix}.
\end{align*}
Observe that if $A(W)\in\mathcal{H}_2(\Lambda)$, for any $W=\2^k$, with  $k\ge1$ (and likewise for $W=\0\2^k$, $\1\0\2^k$ or $\0\1\0\2^k$) then $\Lambda=+\infty$, while for  $W\in\bigcup_{k=1}^{+\infty}\{\1\0\1\0\2^k\}$ it is simple to check that $A(W)\in\mathcal{H}_2(7)$ (for indeed
$2\le 7\cdot \1$, $5\le 7\cdot\1$, $7k+5\le 7\cdot (k+\1)$, $7k+5\le 7\cdot(2k+\1)$
are valid inequalities, for any $k\ge 1$). Using this algorithm, leads to the language:\hfil\break
$\mathcal{W}:=\Big\{\1\Big\}\cup\bigcup\limits_{k=1}^{+\infty}\Big\{\0\1\0^k,\1\1\0^k,\2\1\0^k,\2\0^k,\0\0\2^k,\0\0\1\0\2^k,\1\0\1\0\2^k,\2\0\1\0\2^k,\1\1\0\2^k,\2\1\0\2^k,\2\0\2^k,\1\2^k\Big\}$,
together with the fact that $A(\1)\in\mathcal{H}_2(2)$, while for any $k\ge\1$
$$
\begin{matrix}
\hfill A(\0\1\0^k)&\!\!\!\!\in\mathcal{H}_2(8)\hfill
&\hfill A(\1\1\0^k)&\!\!\!\!\in\mathcal{H}_2(8)\hfill
&\hfill A(\2\1\0^k)&\!\!\!\!\in\mathcal{H}_2(6)\hfill
&\hfill A(\2\0^k)&\!\!\!\!\in\mathcal{H}_2(7)\hfill\\
\hfill A(\0\0\2^k)&\!\!\!\!\in\mathcal{H}_2(12)\hfill
&\hfill A(\0\0\1\0\2^k)&\!\!\!\!\in\mathcal{H}_2(13)\hfill
&\hfill A(\1\0\1\0\2^k)&\!\!\!\!\in\mathcal{H}_2(7)\hfill
&\hfill A(\2\0\1\0\2^k)&\!\!\!\!\in\mathcal{H}_2(6)\hfill\\
\hfill A(\1\1\0\2^k)&\!\!\!\!\in\mathcal{H}_2(8)\hfill
&\hfill A(\2\1\0\2^k)&\!\!\!\!\in\mathcal{H}_2(6)\hfill
&\hfill A(\2\0\2^k)&\!\!\!\!\in\mathcal{H}_2(5)\hfill
&\hfill A(\1\2^k)&\!\!\!\!\in\mathcal{H}_2(9)\hfill
\end{matrix} 
$$
These computations also ensure that  $A(W)\in{\mathcal H}_3(5)$, for any $W\in\mathcal{W}$.
\begin{lemma}\label{H2}$A(W)\in{\mathcal H}_2(13)\cap{\mathcal H}_3(5)$, for any $W\in\mathcal{W}$.
\end{lemma}

\begin{lemma}\label{split}Given $\omega=\omega_1\omega_2\cdots\in\{\0,\1,\2\}^\NNN$ one has 

(i)~:~for any   $n\ge1$, either $\omega_1\cdots\omega_n$ is a strict suffix $W_0(n)$ of a word in $\mathcal{W}$ or there exists a possibly empty strict suffix $W_0(n)$ of a word in $\mathcal W$ and $W_1(n),\dots,W_{k_n}(n)\in\mathcal W$ s.t. 
\begin{equation*}\label{wk}
\omega_1\cdots\omega_n=W_0(n)W_1(n)\cdots W_{k_n}(n)\;;
\end{equation*}

(ii)~:~if $\omega\in\XX$ (i.e. $\sigma^k\cdot\omega$ never belongs to $\{\bar\0,\bar\2\}$) then,  there exist $W_1,W_2,\dots$ in $\mathcal W$ and $W_0$ a possibly empty strict suffix of a word in $\mathcal W$, s.t.
\begin{equation*}\label{wkbis}
\omega_1\omega_2\dots=W_0W_1W_2\cdots
\end{equation*}
\end{lemma}

\begin{proof}[{\bf Proof}](i) : Let $n\ge 1$ : by definition of $\mathcal{W}$, either $\omega_1\dots\omega_n$ is a strict suffix of a word of $\mathcal W$ or $\omega_1\dots\omega_n$ has a suffix $W'\in\mathcal W$. In the later case, there exists a word $W$ such that $\omega_1\dots\omega_n=WW'$.  If $W$ is neither empty nor a strict suffix of a word in $\mathcal{W}$, we repeat the procedure with $W$ and so on and so forth with a finite induction. 

(ii) : We use the decomposition in part (i) for any $n\ge 1$: for the word $W_0(n)$ is a strict suffix of a word of $\mathcal W$, there exists $k\ge 1$ such that
$W_0(n)\in\{\emptyword,\1\0^k,\0^k,\0\1\0\2^k,\1\0\2^k,\0\2^k,\2^k\}$.
Since $\omega\in\XX$, the set $\{W_0(n)\;;\;n\ge1\}$ is necessarily  finite. Hence, there exists a word $W_0$ and an infinite set $E_0\subset\{1,2,\dots\}$ such that $W_0(n)=W_0$, for any $n$ in $E_0$.
The words $W_0,W_1,\dots$ are defined by induction: given $k\ge0$, let $W_0,\dots,W_k$ be k+1 words  for which $W_0(n)=W_0,\dots,W_k(n)=W_k$ for any $n$ in an infinite set $E_k\subset\{1,2,\dots\}$. Since $\omega\in\XX$, there cannot exist infinitely many words $W\in\mathcal W$ such that $W_0\dots W_kW$ is a prefix of $\omega_1\dots\omega_n$ for a $n\ge1$. Hence, there exists  $W_{k+1}$ s.t. $W_0(n)=W_0,\dots,W_{k+1}(n)=W_{k+1}$ for any $n$ in an infinite set $E_{k+1}\subset E_k$: the induction holds, so that $\omega_1\omega_2\cdots=W_0W_1W_2\cdots$.

\end{proof}

The language $\mathcal{W}$ has two  important drawbacks w.r.t. condition ${\bf (C)}$. First, for many $W$ in $\mathcal{W}$ the  matrices $A(W)$ belongs to $\mathcal{H}_1$, but (unfortunately) there are words in $\mathcal{W}$ for which this is not true: for instance, $A(\0\1\0)=
\begin{pmatrix}
0&0&0&1&1&0&0\\
1&0&0&0&0&0&1\\
1&0&0&1&1&0&0\\
0&0&0&0&0&0&0\\
0&0&0&1&1&0&0\\
0&0&0&1&1&0&0\\
0&0&0&1&1&0&0\end{pmatrix}$
(we shall see (Lemma~\ref{13}) that $A(W)\in\mathcal{H}_1$ as soon as $W$ is factorized by a concatenation of at least 13 words in $\mathcal{W}$). Secondly, the  constant $\lambda=5$ in Lemma~\ref{H2} is larger than $1$. The following lemma is crucial (its proof is given in \S~\ref{proofH3} below).

\begin{lemma}\label{H3} If $W\in\{\0,\1,\2\}^*$ is factorized by a concatenation of $130$ ($=10\cdot 13$) words in $\mathcal W$,~then: 
$$
A(W)\in\mathcal{H}_1\cap{\mathcal H}_3\left(3/4\right).
$$
\end{lemma}

\begin{proof}[{\bf Proof of Proposition~\ref{hyp}}] The fact that $\bar\0$ and $\bar\2$ are ${\bf (C)}$-singular may be checked directly. To prove that an arbitrary given $\omega\in\XX$ is necessarily ${\bf (C)}$-regular we begin with Lemma~\ref{split} to write $\omega=W_0W_1W_2\cdots$, where $W_0$ is a possibly empty strict suffix of a word in $\mathcal{W}$ and each $W_1,W_2,\dots$ are words in $\mathcal{W}$. Let  $0=s_0=s_1$
and for $k\ge1$ define:  
$$
s_{k+1}=\sum\nolimits_{i=0}^{130k}\vert W_i\vert.
$$
This leads to the recoding 
$\omega=\tilde{W}_{s_2}\tilde{W}_{s_3}\cdots$, where $\tilde{W}_{s_2}=W_0W_1\cdots W_{130}=\omega_1\cdots\omega_{s_2}$, while for $k\ge2$,
$$
\tilde{W}_{s_{k+1}}=W_{130(k-1)+1}\cdots W_{130k}=\omega_{s_{k}+1}\cdots\omega_{s_{k+1}}.
$$
The latter notations are consistent with setting $\tilde{W}_{0}=\emptyword$ and for $n\ge1$, $k\ge0$ and $s_{k+1}\le n<s_{k+2}$
$$
\tilde{W}_n=\omega_{s_{k}+1}\cdots\omega_n=\tilde{W}_{s_{k+1}}\tilde{W}_n'
\quad\text{where}\quad\tilde{W}_n':=\omega_{s_{k+1}+1}\cdots\omega_n\;;
$$
in other words, $\omega_1\cdots\omega_n=\tilde{W}_n$ when $n<s_3$, while for $n\ge s_3$,
$$
\omega_1\cdots\omega_n=\tilde{W}_{s_2}\cdots\tilde{W}_{s_k}\tilde{W}_n=\tilde{W}_{s_2}\cdots\tilde{W}_{s_k}\tilde{W}_{s_{k+1}}\tilde{W}_n'.
$$
(For $s_{k+1}\le n<s_{k+2}$, the word $\tilde{W}_n'$ a strict prefix  of $\tilde{W}_{s_{k+2}}$.)
\begin{figure}[H]
\begin{center}
\includegraphics[trim=0 150 70 150,clip,scale=0.58]{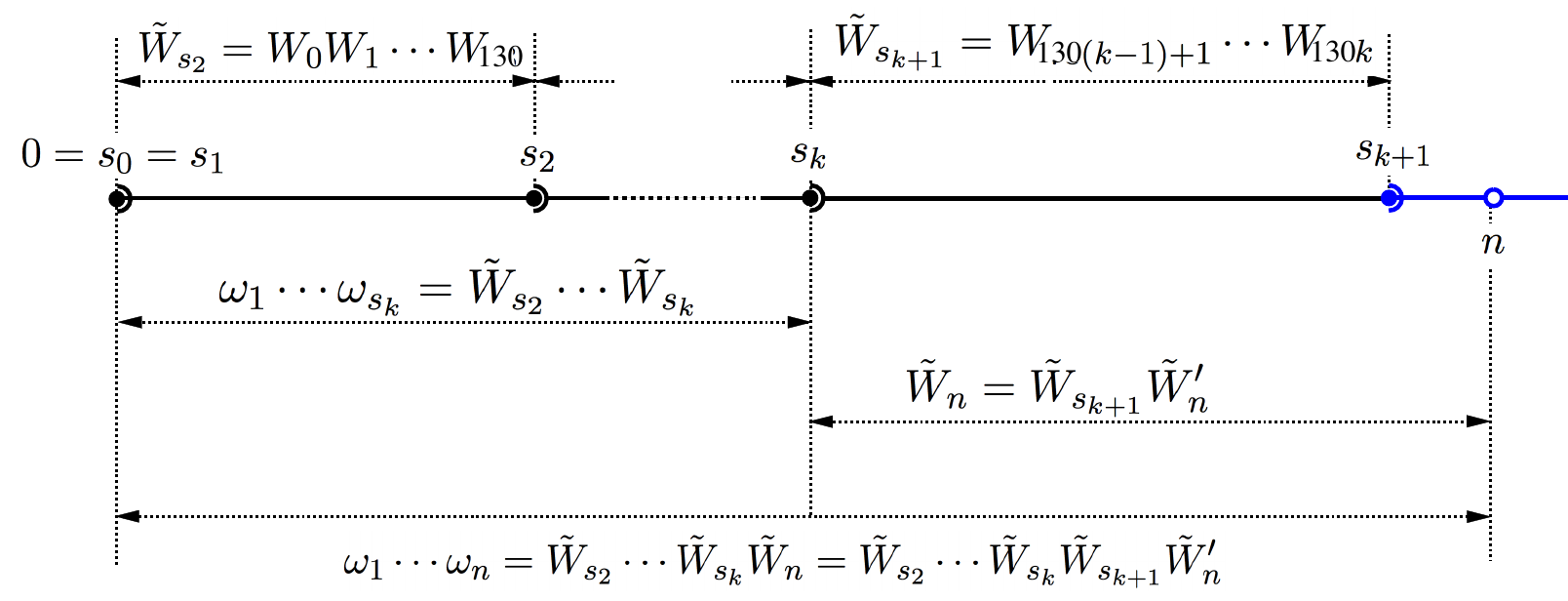}
\caption{\label{scales117}\footnotesize\it Each $\omega\in\XX$ may be recoded by writing $\omega_1\omega_2\cdots=\tilde{W}_{s_2}\tilde{W}_{s_3}\cdots$, where each word $\tilde{W}_{s_k}$ is factorized by a concatenation of $130$ words in $\mathcal{W}$; the sequence $0=s_0=s_1<s_2<\cdots$ for condition ${\bf (C)}$ is defined in consequence, that is $s_k=|\tilde{W}_2\cdots \tilde{W}_{s_k}|$ and for each $s_{k+1}\le n<s_{k+2}$ (with $k\ge0$)
$\tilde{W}_n=\omega_{s_k+1}\cdots\omega_n=\tilde{W}_{s_{k+1}}\tilde{W}_n'$
is also factorized by a concatenation of $130$ words in $\mathcal{W}$; the importance of the number $130$ comes from crucial Lemma~\ref{H3}.} 
\end{center}
\end{figure}
With the notations of Definition~\ref{C}, one has $P_n=A(\omega_1\cdots\omega_n)=P_{s_{k}}Q_n$, where we stress on the key correspondence
\begin{equation}\label{correspondence} 
Q_n=A(\tilde{W}_n).
\end{equation}
The first and the third parts of condition ${\bf (C)}$ are satisfied since by definition of $\tilde{W}_n$ and of the sequence $(s_0,s_1,\dots)$, Lemma~\ref{H3} ensures that  
\begin{equation}\label{mmmmm} 
n\ge s_2\Longrightarrow Q_n=A(\tilde{W}_{n})\in{\mathcal H}_1\cap {\mathcal H}_3\left(3/4\right)
\end{equation}
It remains to establish the existence of a constant $\Lambda\ge1$ (depending on $\omega$) s.t.
\begin{equation}\label{lllll} 
n\ge s_2\Longrightarrow Q_n=A(\tilde{W}_n)\in{\mathcal H}_2(\Lambda).
\end{equation}
Suppose  $s_2\le s_{k+1}\le n<s_{k+2}$; by construction, $\tilde{W}_n=\tilde{W}_{s_{k+1}}\tilde{W}_n'$, where $\tilde{W}_n'$ is a (possibly empty) strict prefix of $\tilde{W}_{s_{k+2}}$: actually, either $\tilde{W}_n'=W_n'$ or 
$\tilde{W}_n'=W_{130k+1}\cdots W_{130k+m}W_n'$,
where $1\le m<130$, while $W_n'$ is a strict prefix of a word in $\mathcal{W}$. However (Lemma~\ref{H2}) we know  that $A(W_i)\in\mathcal{H}_2(13)\cap\mathcal{H}_3(5)$, for any $i>0$ and by part (ii) of Lemma~\ref{hypotheses}),
$$
A(\tilde{W}_n)\in
\mathcal{H}_2\Big(\left(1+5+\cdots+5^{2\cdot130-2}\right)\cdot 13+5^{2\cdot 130-1}\cdot\Lambda_{A(W_n')}\Big)\;;
$$ 
therefore $A(\tilde{W}_n)\in\mathcal{H}_2\big(5^{260}(4+\Lambda')\big)$, 
where $\Lambda'$ stands for  the maximum of the $\Lambda_{A(W')}$ over $W'$ being any prefix of a word in $\mathcal{W}$. The inequality $\Lambda'<+\infty$ is valid, because either $W'\in\mathcal{W}$ and $A(W')\in\mathcal{H}_2(13)$ or $W'\not\in\mathcal{W}$ which is only possible for
$$
W'\in\{\emptyword,\0,\0\0,\0\0\1,\0\0\1\0,\0\1,\1\1,\1\0,\1\0\1,\1\0\1\0,\2,\2\0\1,\2\0\1\0,\2\1\}.
$$

\end{proof}

\subsection{Proof of Lemma~\ref{H3}}\label{proofH3}For Proposition~\ref{hyp} to be completely established, it remains to prove Lemma~\ref{H3}. The argument --~given in \S~\ref{Argument}~-- depends on a key lemma (Lemma~\ref{13}) established in \S~\ref{Key} together with several properties of  two adjacency graphs --~$(\Gamma_1)$ and $(\Gamma_2)$~-- that we  consider in  $\S~\ref{Gamma1}$ and $\S~\ref{Gamma2}$ respectively.  

\subsubsection{\bf Key Lemma}\label{Key}The following lemma shows --~in part (iii)~-- that $A(W)\in{\mathcal H}_1$ as soon as $W$ is factorized by a concatenation of at least $13$ words in $\mathcal W$. Parts (i) and (ii) will be determinant in the final argument proving Lemma~\ref{H3} in \S~\ref{Argument}, and use the following set of four column vectors (that we represent by convenience by four $7$-uplets):
\begin{equation}\label{defmathcalT2}
\mathcal{T}_2:=\Big\{(1,0,1,1,1,0,0),(1,1,1,0,1,1,0),(1,1,1,0,1,1,1),(1,1,1,1,1,0,0)\Big\}.
\end{equation}
 
\begin{lemma}\label{13}If $W\in\{\0,\1,\2\}^*$ is factorized by a concatenation of $13$ words in $\mathcal W$, then 

(i) : $\big\{\Delta(A(W)U_1),\Delta(A(W)U_3),\Delta(A(W)U_5)\big\}\cap\mathcal{T}_2\ne\emptyset$;

(ii) : $\#{\bf Col}(A(W))\le 2$ ;

(iii) : $A(W)\in\mathcal{H}_1$.
\end{lemma}

\begin{proof}[{\bf Proof}]
We shall begin to verify that any $W=W_1\cdots W_{13}$ with $W_i\in\mathcal{W}$ has a factor $W'$ s.t. $A(W')$ satisfies (i--ii--iii).
First, suppose that $\0\1$ is a factor of $W_2\dots W_{10}$ and write $W_2\dots W_{10}=X\0\1Y$.
On the one hand, it is necessary that $W_1X\ne\0$ and thus, $W_1X$ has suffix in $\mathcal{L}':=\{\0\0, \1\0, \2\0, \1, \2\}$; on the other hand, it is also necessary that $YW_{11}W_{12}W_{13}$ has prefix in 
$$
\mathcal{L}'':=\{\0\0, \0\1, \0\2, \1\0, \1\1\0, \1\1\1, \1\1\2, \1\2, \2\0, \2\1, \2\2\}.
$$
Hence (see Figure~\ref{ARBRE2} -- left) $W$ is factorized by a word 
belonging to $\mathcal{L}'\0\1\mathcal{L}''$ and (see \S~\ref{firstcaseofAppendix} in Appendix~\ref{H1})
each word in $\mathcal{L}'\0\1\mathcal{L}''$ satisfies (i--ii--iii). For the second case, suppose that $\2\1$ is a factor of $W_2\dots W_{10}$ that is, $W_2\dots W_{10}=X\2\1Y$; the word $YW_{11}W_{12}W_{13}$ has necessarily a prefix in 
$$
\mathcal{L}''':=\{\0, \1\0, \1\1\0, \1\1\1, \1\1\2, \1\2, \2\0, \2\1, \2\2\}\;;
$$
like in the first case (see Figure~\ref{ARBRE2} -- right) $W$ is factorized by a word 
$W'\in \mathcal{L}'\2\1\mathcal{L}'''$ and (see \S~\ref{secondcaseofAppendix} in Appendix~\ref{H1})
each word in $\mathcal{L}'\2\1\mathcal{L}'''$ satisfies (i--ii--iii).
\begin{figure}[H]
\begin{center}
       \includegraphics[scale=0.6]{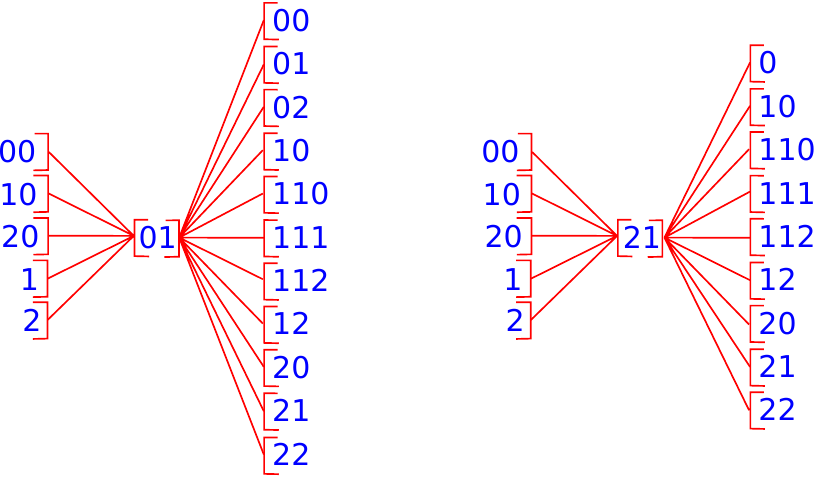}
\caption{\label{ARBRE2}\footnotesize \it How to obtain the words in $\mathcal{L}'\0\1\mathcal{L}''$ (left) and in $\mathcal{L}'\2\1\mathcal{L}'''$ (right).} 
\end{center}
\end{figure}
Finally, it remains to consider the case where neither $\0\1$ nor $\2\1$ is a factor of $W_2\dots  W_{10}$: in particular, this means that $W_2\cdots W_{10}=\1^kW'$ with $k\ge0$ and $W'\in\{\0,\2\}^*$  (it is here  where the $13$ words in $\mathcal{W}$ will prove to be needed).
Notice that  $A(\1),\dots,A(\1\1\1\1)$ do not satisfy (i--ii--iii), while 
$$
A(\1\1\1\1\1)
=\begin{pmatrix}3&0&3&2&1&0&0\\0&0&0&0&0&0&0\\1&0&3&3&1&0&0\\1&0&1&3&2&0&0\\2&0&1&1&1&0&0\\0&0&0&0&0&0&0\\0&0&0&0&0&0&0\end{pmatrix}
$$
does satisfy this property: hence, $W$ has a factor $W'$ ($=\1\1\1\1\1$) with $A(W')$ satisfying (i--ii--iii), as soon as $k\ge5$. Now, suppose that  $k\le 4$, so that $W_6\cdots W_{10}\in\{\0,\2\}^*$ (indeed the worst case arises for $k=4$ with $W_2W_3W_4W_5=\1\1\1\1$).
The definition  of $\mathcal{W}$ implies that each word $W_6$, $W_7$, $W_8$, $W_9$, $W_{10}$ must be either $\2\0^k,\0\0\2^k$ or $\2\0\2^k$ (for some $k\ge1$). Then, it is necessary that $\0\2$ is a factor of both $W_6W_7$ and $W_8W_9$: this implies that $W_6\cdots W_{10}$ has a factor of the form $\0^x\2^y\0^z\2^t\0^u$ with $x,y,z,t,u\ge1$. On the one hand, if $t=1$ then $z=1,2,3,\dots$ and $\0^x\2^y\0^z\2^t\0^u$ is factorized by either $\2\0\2\0$, $\2\0\0\2\0$ or $\0\0\0\2\0$; on the other hand, if $t\ge\2$, the word $\0^x\2^y\0^z\2^t\0^u$ is factorized by $\2\2\0$. One concludes this case, since $W$ is  factorized by a word $W'$ which is either equal to $\2\0\2\0$, $\2\0\0\2\0$, $\0\0\0\2\0$ or $\2\2\0$ and  which (see \S~\ref{thirdcaseofAppendix} in  Appendix~\ref{H1}) satisfies (i--ii--iii).

We now conclude with the general case of $W\in\{\0,\1,\2\}^*$ which is factorized by $W_1\cdots W_{13}$ with $W_i\in\mathcal{W}$. Let $W'$ be the factor of $W_1\cdots W_{13}$ s.t. $A(W')$ satisfies (i--ii--iii) and notice  that $A(W')\in\mathcal{H}_1$. Because $A(W')$ satisfies (i), a direct verification shows that both $A(W'i)$ and $A(iW')$ satisfy (i) for any $i=\0,\1,\2$: hence, by a finite induction, $A(W)$ also satisfies (i). The fact that (ii--iii) holds for $A(W)$ follows from Lemma~\ref{hypothesesbis}: because $A(W')\in\mathcal{H}_1$ and $\#{\bf Col}(A(W'))\le 2$, part (e) of Lemma~\ref{hypothesesbis} ensures that $\#{\bf Col}(A(W))\le 2$ and part (f) that $A(W)\in\mathcal{H}_1$.

\end{proof}

\subsubsection{\bf The adjacency graph $(\Gamma_1)$}\label{Gamma1}We begin to define the infinite adjacency graph $(\Gamma_0)$ whose vertex set ${\rm Vert}(\Gamma_0)$ is made of the nonzero $7\times 1$ vectors (represented for convenience by) the $7$-uplets $(x_1,\dots,x_7)$ having nonnegative integral entries and obtained from one of the basis vector $U_1,\dots,U_7$ by successive left multiplication with the  matrices $A(\0),A(\1)$ and $A(\2)$. For instance, 
$$
(0,0,1,1,0,0,0)\mathop{\longrightarrow}\limits^{\omega_n}\cdots
\mathop{\longrightarrow}\limits^{\omega_1}(x,0,y,z,t,0,0)
$$
is a path in $(\Gamma_0)$ so that $A(\omega_1\cdots\omega_n)(0,0,1,1,0,0,0)=(x,0,y,z,t,0,0)$. Setting $V\sim V'$ whenever $\Delta(V)=\Delta(V')$ gives an equivalent relation on ${\rm Vert}(\Gamma_0)$: the quotient space ${\rm Vert}(\Gamma_0)/_\sim$ is made of finitely many classes of equivalence, each ones being uniquely represented by a $7$-uplet $(a_1,\dots,a_7)\in\{*,0,1,2,\dots\}^7$: here $(a_1,\dots,a_7)$ represents the set $\{V_k\}_{k=1}^N$ ($1\le N\le+\infty$) of all the nonzero $7\times 1$ vectors in a given class of equivalence, with $a_i$ either equal to $\sup_k\{U_i^\star V_k\}$ or $*$ if either $\sup_k\{U_i^\star V_k\}$ is finite or infinite respectfully. This leads to introduce the adjacency graph $(\Gamma_1)$ whose vertex set is ${\rm Vert}(\Gamma_1)={\rm Vert}(\Gamma_0)/_\sim$ an whose paths are obtained by quotient projection from the paths in $(\Gamma_0)$: the representation of $(\Gamma_1)$ is given in  Figure~\ref{GRAPH1GOOD} of Appendix~\ref{ADJACENCY}. The set ${\rm Vert}(\Gamma_1)$ may be written as a partition $\mathcal{V}_1\sqcup\mathcal{V}_2$, where 
$$
\mathcal{V}_2:=\Big\{(*,0,*,*,*,0,0),\;(*,*,*,0,*,*,0),\;(*,*,*,0,*,*,*),\;(*,*,*,*,*,0,0)\Big\}
$$
By an abuse of notations the vertices in ${\rm Vert}(\Gamma_1)$ may be considered as usual $7\times 1$ vectors, with for instance $\Vert (*,0,*,*,*,0,0)\Vert=+\infty$ and $\Delta(*,0,*,*,*,0,0)=(1,0,1,1,1,0,0)$; but we may also note $(x_1,\dots,x_7)\in(y_1,\dots,y_7)$, when $(y_1,\dots,y_7)$ is considered as a class of equivalence. The set of column vectors $\mathcal{T}_2$ considered in Lemma \ref{13} is $\Delta(\mathcal{V}_2)$.

\begin{lemma}[Synchronization lemma]\label{synchronization} The following proposition holds

(a) : for any $W\in\{\0,\1,\2\}^*$ and any $1\le i\le 7$ 
$$
A(W)U_j\ne 0\Longrightarrow \exists V\in\mathcal{T}_2,\;\Delta(A(W)U_j)\le V\;;
$$

(b) : if $V\in\mathcal{V}_1$ then $A(i)V\le (2,2,2,2,2,2,2)$, for any $i=\0,\1,\2$;

(c) : if $V\in\mathcal{V}_2$ then $A(i)V\in\mathcal{V}_2$, for any $i=\0,\1,\2$;

(d) : in the subgraph $(\Gamma_1')$ of $(\Gamma_1)$ with vertices in $\mathcal{V}_2$ (see Figure~\ref{GRAPH1SYNCHRO}), the words $W\in\{\0,\1,\2\}^*$ of length $|W|\ge3$ are synchronizing: in other words, if $V$ and $V'$ are two $d\times 1$ vectors with nonnegative entries s.t. both $\Delta(V)$ and $\Delta(V')$ belong to $\mathcal{T}_2$, then  
$$
|W|\ge3\Longrightarrow \Delta(A(W)V)=\Delta(A(W)V')\;;
$$

(e) : let $V$ be a vector having positive entries; then, for any word $W\in\{\0,\1,\2\}^*$,
$$
W\ne\emptyword\Longrightarrow\Delta(A(W)V)\in\mathcal{T}_2\;;
$$

(f) : let $V$ be a vector having positive entries; then, for any word $W,W'\in\{\0,\1,\2\}^*$,
$$
|W|\ge3\Longrightarrow \Delta(A(W)V)=\Delta(A(WW')V).
$$
\end{lemma}

\begin{proof}[{\bf Proof}] Parts (a), (b) and (c) are verified directly on $(\Gamma_1)$ while (d) is clear from the subgraph  $(\Gamma_1')$ (see Figure~\ref{GRAPH1SYNCHRO}): the words in $\{\0\0,\0\1,\1\1,\2\}$ being synchronizing, each word $W\in\{\0,\1,\2\}^*$ with $|W|\ge3$ is synchronizing as well, i.e. if $W=\omega_1\cdots\omega_n\in\{\0,\1,\2\}^n$ for $n\ge3$~and~if 
$$
V_1\mathop{\longrightarrow}\limits^{\omega_1}\cdots\mathop{\longrightarrow}\limits^{\omega_n}V_2
\quad\text{and}\quad
V_1'\mathop{\longrightarrow}\limits^{\omega_1}\cdots\mathop{\longrightarrow}\limits^{\omega_n}V_2'
$$
are two paths in $(\Gamma_1')$ then it is necessary that $V_2=V_2'$. For parts (e) and (f), take  a vector $V$ with positive entries: then by direct verification one gets $\Delta(A(i)V)\in\mathcal{T}_2$, for $i=\1,\2,\3$ and  (e) follows from (c); part (f) is a consequence (e) and of the synchronization property in (d).

\begin{figure}[H]
\includegraphics[scale=0.47]{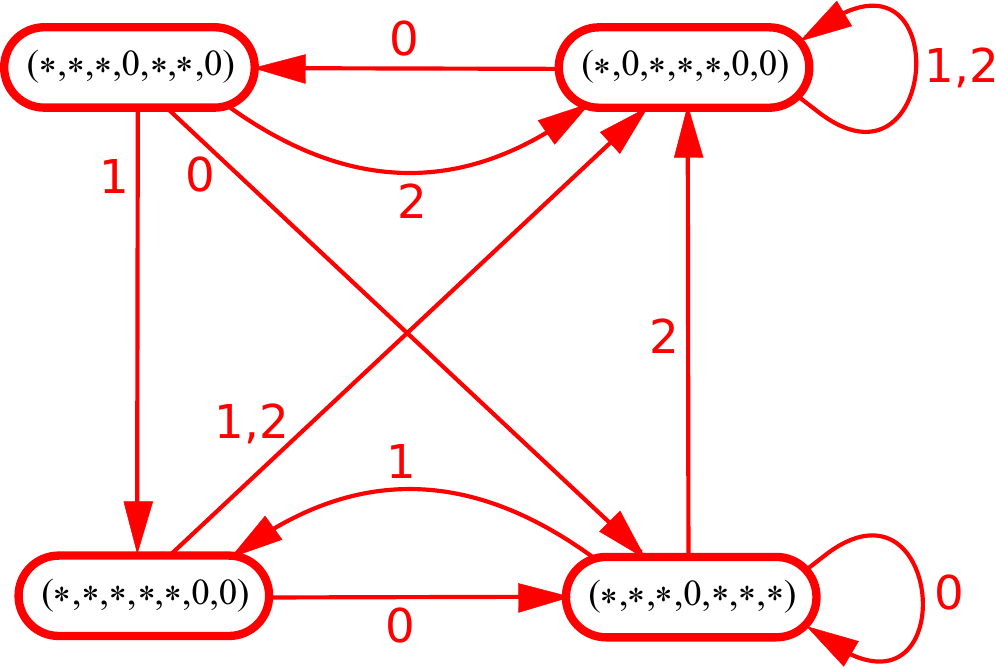}
\caption{\label{GRAPH1SYNCHRO}\small\it Subgraph $(\Gamma_1')$ of $(\Gamma_1)$ with vertices $(a_1,\dots,a_7)\in\mathcal{V}_2$:  each word in $W\in\{\0\0,\0\1,\1\1,\2\}$ is synchronizing (in this subgraph) in the sense that each path coded by $W$ ends on the same vertex: it follows that a word $W'\in\{\0,\1,\2\}^*$ with $|W'|\ge3$ is synchronizing as well.}
\end{figure}

\end{proof}

\subsubsection{\bf The adjacency graph $(\Gamma_2)$}\label{Gamma2}The final argument leading to Lemma~\ref{H3} stands on a {\it doubling property} (Lemma~\ref{H30}) displayed by a second adjacency graph  $(\Gamma_2)$, related to the graph $(\Gamma_1)$ introduced in \S~\ref{Gamma1}, but  defined in a slightly different way.
Each $V\in{\rm Vert}(\Gamma_2)$ is a $d\times 1$ nonzero vector with nonnegative integral entries and
represents the set of  $d\times 1$ vectors $V'\ne0$ with nonnegative integral entries such that
$\Delta(V')=\Delta(V)$ and $V'\ge V$: then, the labelled edge $V\mathop{\longrightarrow}\limits^{i}V'$ in $(\Gamma_2)$ means 
$A(i)V\ge V'\quad\text{and}\quad\Delta(A(i)V)=\Delta(V')$.
By definition $\mathcal{T}_2\subset{\rm Vert}(\Gamma_2)$, all the labeled edges of $(\Gamma_2)$ are obtained by starting from one vertex in $\mathcal{T}_2$ and making all the possible (left) multiplication by either $A(\0),A(\1)$ or $A(\2)$. Among the infinitely many possibilities of an adjacency graph satisfying the previous specifications, we have chosen one --~namely $(\Gamma_2)$~-- with finitely many vertices with the additional property that each  $V\in{\rm Vert}(\Gamma_2)$ satisfies $V\le (2,\dots,2)$. (A representation of $(\Gamma_2)$ is given in Figure~\ref{GRAPH2GOOD} in  Appendix~\ref{ADJACENCY}).   
We shall use the important fact that {\it most} of the paths of $(\Gamma_2)$ starting from a vertex in $\mathcal{T}_2$
terminate on a vertex in $2\mathcal{T}_2$, while all paths starting from a vertex in $2\mathcal{T}_2$ terminate on a vertex in $2\mathcal{T}_2$.

\begin{lemma}[Doubling property]\label{H30} 
Let $W=\omega_1\cdots\omega_n\in\{\0,\1,\2\}^n$ be factorized by a concatenation of~\footnote{Actually $12$ words are sufficient: we state the lemma with 13 instead because of the 13 words in Lemma~\ref{13}.} $13$ words in $\mathcal W$ but not factorized by $\0\0\1\0\0\1$;  
then, for any nonnegative~vector~$V$, 
$$
\Delta(V)\in\mathcal{T}_2\Longrightarrow \Big(\Delta(A(W)V)\in \mathcal{T}_2\quad\text{and}\quad V\ge a\cdot \Delta(V)\Longrightarrow A(W)V\ge2a\cdot\Delta(A(W)V)\Big).
$$
\end{lemma}
 
\begin{proof}[{\bf Proof}] 
Suppose that $W=\omega_1\cdots\omega_n\in\{\0,\1,\2\}^*$ is factorized by  $\omega_i\cdots\omega_j=W_1\cdots W_{13}$ with $W_i\in\mathcal W$ but not factorized by $\0\0\1\0\0\1$. Hence, the mirror word $\omega_j\cdots\omega_i$ is a concatenation of at least $13$ words of the form $\0^k$, $\1$ or $\2^k$ with $k\ge1$, without being  factorized by $\1\0\0\1\0\0$. 
We claim  that a path in $(\Gamma_2)$, labelled by $\omega_j\cdots\omega_i$ and starting from a vertex in $\mathcal{T}_2$ terminates in $2\mathcal{T}_2$:  
indeed (see Figure~\ref{GRAPH2GOOD}) the  longest paths 
(longest w.r.t. the number of blocks of the form $\0^k$, $\1$ or $\2^k$),   
not factorized by $\1\0\0\1\0\0$ and whose terminal vertex  does not belong to $2\mathcal{T}_2$,
are those joining $(1,0,1,1,1,0,0)$ to $(2,2,2,0,2,2,1)$ 
and labelled by words of the form $\0\0\0\1\0\0\1\0\2\1\2^k\0\1\0$; these words are concatenations of the $11$ words $\0\0\0$, $\1$, $\0\0$, $\1$, $\0$, $\2$, $\1$, $\2^k$,~$\0$,~$\1$,~$\0$.

\end{proof}

\subsubsection{\bf Proof of Lemma~\ref{H3}}\label{Argument}Consider $W:=\omega_1\cdots\omega_n=XW_1\cdots W_{10}Y$, where $X,Y\in\{\0,\1,\2\}^*$ and  each $W_1,\dots,W_{10}$ is a concatenation of $13$ words in $\mathcal W$ (making $W_1\cdots W_{10}$ a concatenation of $130$ words in $\mathcal{W}$). The assertion $A(W)\in\mathcal{H}_1$ is proved in Lemma~\ref{13}. From this lemma, $\#{\bf Col}(A(W_{10}Y))\le 2$ and  there exists at least one column index $j_0$ such that $\Delta\big(A(W_{10}Y)U_{j_0}\big)\in\mathcal{T}_2$ and $\Delta\big(A(W_{10}Y)U_{j}\big)\le\Delta\big(A(W_{10}Y)U_{j_0}\big)$ for any $j$:
hence, by part (c) of Lemma~\ref{synchronization} 
\begin{equation}\label{j1}
\Delta\big(A(W_{10}Y)U_{j_0}\big),\dots,\Delta\big(A(XW_1\cdots W_{9}W_{10}Y)U_{j_0}\big)=\Delta\big(A(W)U_{j_0}\big)\in\mathcal{T}_2.
\end{equation}
If $\#{\bf Col}(A(W))=1$ then $A(W)\in{\mathcal H}_3\left(\lambda\right)$ for any~$\lambda\ge0$. We now  assume that  $\#{\bf Col}(A(W))=2$ so that it is licit to consider any arbitrary column index, say ${j'}$, s.t. $0<\Delta\big(A(W)U_{{j'}}\big)<\Delta\big(A(W)U_{j_0}\big)$. 
One deduce from (\ref{j1}) and the synchronization property (part (d) of Lemma \ref{synchronization}) that
\begin{equation}\label{j2}
\mathcal{T}_2\not\ni\Delta\big(A(\omega_4\cdots\omega_n)U_{{j'}}\big)<\Delta\big(A(\omega_4\cdots\omega_n)U_{j_0}\big)\in\mathcal{T}_2.
\end{equation}
However --~in view of graph $(\Gamma_1)$~-- this is possible only if 
$A(\omega_4\cdots\omega_n)U_{{j'}}\in{\rm Vert}(\Gamma_1)\setminus\mathcal{V}_2=\mathcal{V}_1$, 
so that (part (b) of Lemma \ref{synchronization})
\begin{equation}\label{maj1}
A(\omega_4\cdots\omega_n)U_{{j'}}\le (2,2,2,2,2,2,2)\quad\hbox{and}\quad A(W)U_{{j'}}\le 3^3\cdot(2,2,2,2,2,2,2).
\end{equation}
We claim that the following implication holds for any $j$:
\begin{equation}\label{maj2}
\Delta\big(A(W_{10}Y)U_j\big)\in\mathcal{T}_2\Longrightarrow A(W)U_{j}\ge 2^{9}\cdot \Delta\big(A(W)U_{j}\big).
\end{equation}
Take any $j$ s.t. $\Delta(A(W_{10}Y)U_j)\in\mathcal{T}_2$ and ${j'}$ with $0<\Delta(A(W)U_{{j'}})<\Delta(A(W)U_j)$. It follows from (\ref{maj1}) and (\ref{maj2}) that
$\Vert A(W)U_{{j'}}\Vert\le 7\cdot 3^3\cdot2\le \left(7\cdot 3^3\cdot2/2^9\right)A(W)U_j(i)$, for any index $1\le i\le 7$ s.t. $A(W)U_j(i)\ne 0$: this means 
$A(W)\in{\mathcal H}_3\left(7\cdot 3^3\cdot2/2^9\right)\subset{\mathcal H}_3\left(3/4\right)$.

To complete the proof of Lemma~\ref{H3} it remains to establish (\ref{maj2}), for which the graph $(\Gamma_2)$ is crucial. The key point is $\#{\bf Col}(A(W))=2$ implies  $W$ cannot be factorized by  $\0\0\1\0\0\1$: indeed, on the contrary --~for $A(\0\0\1\0\0\1)$ being  a rank one matrix~-- the matrix $A(W)$ would be also of rank one and  $\#{\bf Col}(A(W))$ would be equal to $1$. In particular each word $W_9,W_8,\dots,XW_1$ is not factorized by $\0\0\1\0\0\1$ and applying  Lemma~\ref{H30} yields: 
\begin{eqnarray*}
A(W_{10}Y)U_{j}&\ge& 1\cdot \Delta\big(A(W_{10}Y)U_{j}\big)\in\mathcal{T}_2\\
A(W_{9})A(W_{10}Y)U_{j}=A(W_{9}W_{10}Y)U_{j}&\ge& 2\cdot\Delta\big(A(W_{9}W_{10}Y)U_{j}\big)\in\mathcal{T}_2\\
A(W_{8})A(W_{9}W_{10}Y)U_{j}=A(W_{8}W_{9}W_{10}Y)U_{j}&\ge& 2^2\cdot\Delta\big(A(W_{8}W_{9}W_{10}Y)U_{j}\big)\in\mathcal{T}_2\\
&\vdots&\\
A(XW_1)A(W_2\cdots W_{10}Y)U_{j}=A(XW_1\cdots W_{10}Y)U_{j}&\ge& 2^{9}\cdot\Delta\big(A(XW_1\cdots W_{10}Y)U_{j}\big)\in\mathcal{T}_2.
\end{eqnarray*}
Hence (\ref{maj2}) is proved as well as Lemma~\ref{H3}.

\hfill\qed

\subsection{Proof of Theorem~\ref{3matrices}}\label{unif7X7} Let $V\in\mathcal{S}_7$ with positive entries. By part (e) of Lemma~\ref{synchronization} the subshift  $\Omega_V$ --~i.e. the compact shift-invariant subset of $\{\0,\1,\2\}^\NNN$ made  of the sequences $\omega$ s.t. $A(\omega_1\cdots\omega_n)V\ne0$, for any $n\ge1$~-- coincides with  $\{\0,\1,\2\}^\NNN$. Moreover, according to Proposition~\ref{hyp}, the  ${\bf (C)}$-regular $\omega$  --~in the sense that $\mathcal{A}(\omega)=(A(\omega_1),A(\omega_2),\dots)$ satisfies condition ${\bf (C)}$~-- form the dense subset  $\{\0,\1,\2\}^\NNN$ that is
$$
\XX:=
\left(\{\0,\1,\2\}^\NNN\setminus\bigcup_{n=0}^{+\infty}\sigma^{-n}\{\bar\0\}\right)\cap\left(\{\0,\1,\2\}^\NNN\setminus\bigcup_{n=0}^{+\infty}\sigma^{-n}\{\bar\2\}\right).
$$
To establish the uniform convergence of $\Pi_n(\omega,V)$ for $\omega\in\{\0,\1,\2\}^\NNN$, we want to apply Propositon~\ref{uuuuuu}, which is the suitable form of Theorem A for the uniform convergence (Recall that we note  $\omega_{i,j}=\omega_{i+1}\dots\omega_j$ for any $\omega=\omega_1\omega_2\cdots\in\{0,1,2\}^{\mathbb N}$ and $0\le i\le j$,  in particular $\omega_{i,i}=\emptyword$.)

$\bullet$ First, consider $\omega\in \XX$.  To check the conditions {\rm (U1.1)} and {\rm (U1.2)} in part (i) of Propositon~\ref{uuuuuu}, let $n\ge1$ and consider  $\xi\in[\omega_1\cdots\omega_n]$ (i.e. $\xi_{0,n}=\omega_{0,n}$). Part (ii) of Lemma~\ref{split} ensures that $\omega=\omega_1\omega_2\cdots=W_0W_1W_2\cdots$, where $W_0$  is a possibly empty  strict suffix of a word in $\mathcal W$, and $W_1,W_2,\dots$ are words in  $\mathcal{W}$. Let $p=p(n)\ge0$ be the maximal rank such that $W_0\cdots W_p=\omega_{0,\psi'(n)}$ is a prefix of $\omega_{0,n}$. For $n$ large enough, there exists  an integer $\psi(n)$, with ~$0\le \psi(n)\le \psi'(n)\le n$ and~s.t.
$$
W_0\dots W_{p-3}=\omega_{0,\psi(n)}\quad\hbox{and}\quad W_{p-2}W_{p-1}W_p\xi_{\psi'(n),n+r}=\xi_{\psi(n),n+r},
$$
where $r\ge0$ is arbitrarily given. Notice that $n\mapsto\psi(n)$ only depends on $\omega$ and that $\psi(n)\to+\infty$ as $n\to+\infty$. Consider the vector
\begin{equation}\label{defXr}
X_r=A(\xi_{\psi(n),n+r})V=A(W_{p-2}W_{p-1}W_p\xi_{\psi'(n),n+r})V\;;
\end{equation}
here, the key point is the synchronization property (f) in Lemma~\ref{synchronization} which implies 
$\Delta(X_r)=\Delta(X_0)$, for any $r\ge0$, proving that {\rm (U1.2)} holds. As for {\rm (U1.1)}, it remains to prove the following lemma.
\begin{lemma}\label{sss}
The exists a constant $\Lambda\ge1$, such that, $A(\xi_{\psi(n),n+r})\in\mathcal{H}_2(\Lambda)$, for any $r\ge0$. 
\end{lemma}
\begin{proof}[{\bf Proof}]According to part (i) of Lemma~\ref{split}), it is possible to write
$\xi_{\psi(n),n+r}=W'_0\dots W'_{q}$
where $W_0'$  is a possibly empty strict suffix of a word in $\mathcal W$, and $W_1',\dots,W_q'$   are words in  $\mathcal{W}$.  For $W'_0$ being a strict suffix of a word of $\mathcal W$, there exists $k\ge1$ s.t. 
$$
W'_0\in\{\emptyword,\1\0^k,\0^k,\0\2^k,\2^k,\0\1\0\2^k,\1\0\2^k\}
$$
and since $W'_0$ is also  a prefix of $W_{p-2}W_{p-1}W_p\xi_{n,n+r}$, it is necessary that
$$
W'_0\in\{\emptyword,\1\0,\1\0\0,\0,\0\0,\2,\0\1\0\2\}.
$$
We note $\lambda_0$ and $\Lambda_0$ s.t. $A(W)\in\mathcal{H}_2(\Lambda_0)\cap\mathcal{H}_3(\lambda_0)$, for any $W\in\{\emptyword,\1\0,\1\0\0,\0,\0\0,\2,\0\1\0\2\}$. Moreover (Lemma~\ref{H2}) $A(W_i')\in\mathcal{H}_2(13)\cap\mathcal{H}_3(5)$ for any $1\le i\le q$. With  $q=130\cdot a+b$ ($0\le b<130$) being the euclidean division of $q$ by $130$, one  writes
$$
W_0'\cdots W_q'=W_0'W_1'\cdots W_b'\tilde{W}_1'\cdots \tilde{W}_a',
$$
where each word $\tilde{W}_i'$ is a concatenation of $130$ words in $\mathcal{W}$: part (ii) of Lemma~\ref{hypotheses} implies $W_1'\cdots W_b'\in\mathcal{H}_2(\Lambda')$ and $\tilde{W}_i'\in\mathcal{H}_2(\Lambda')$ for 
$\Lambda'={13}\sum_{k=0}^{129}5^k$ while Lemma~\ref{H3} ensures $\tilde{W}_i'\in\mathcal{H}_3(3/4)$. 
Finally (part (ii) of Lemma~\ref{hypotheses} again), one concludes that  $A(W'_0\dots W'_{q})\in\mathcal{H}_2(\Lambda)$, where 
\begin{equation}\label{universal}
\Lambda=\Lambda_0+\lambda_0\Lambda'+\lambda_05^b\Lambda'\sum_{k=0}^{+\infty}\left({3\over 4}\right)^k.
\end{equation}

\end{proof}

$\bullet$ Consider now $\omega=w\bar\0$, where $w$ is a (possibly empty) word on $\{0,\dots,\aaa\}$, we shall verify the assertions ${\rm (U2)}$ and ${\rm (U3)}$ in part (ii) of Proposition~\ref{uuuuuu}. First, ${\rm (U2.1)}$ is satisfied, for one checks directly that 
$$
\lim_{n\to+\infty}\left\Vert {A(\0^n)\over \Vert A(\0^n)\Vert}-C_\0U_1^\star\right\Vert=0,
\quad\text{where}\quad
C_\0=(0,1/5,1/5,0,1/5,1/5,1/5),
$$ 
and ${\rm (U2.2)}$ is satisfied because $A(\0)C_\0=C_\0$ and $A(\sss_1\sss_2)C_\0\in\mathcal{V}_2$ for $\sss_1\sss_2\in\{\0,\1,\2\}\times\{\1,\2\}$.

To verify ${\rm (U3.1)}-{\rm (U3.2)}-{\rm (U3.3)}$,  consider $n$ large enough in order that $\sigma^{n-2}\cdot\omega=\bar\0$, let $\xi\in[\omega_1\cdots\omega_n]$ and $r\ge0$. Let $m$ be the integer s.t. $\xi_{n-2,m+2}=\0^{m-n+4}$ with $n\le m+2\le n+r$ and maximal w.r.t this identity (in the notations of Proposition~\ref{uuuuuu} we take $\varphi(n)=n-2$): hence condition ${\rm (U3.1)}$ --~i.e. $\xi\in[\omega_1\cdots\omega_m]$~-- is satisfied by definition. Then, there exists a word $W$, either empty or beginning by $\1$ or $\2$, such that
$\xi_{n-2,n+r}=\0^{m-n+4}W$.
According to part (i) of Lemma \ref{split}, we write $W=W_0\dots W_{p}$
where $W_0$ is a possibly empty strict suffix of a word in $\mathcal W$, and $W_1,\dots,W_p$   are words in  $\mathcal{W}$. If $W_0\ne\emptyword$ we claim that $W_0\in\bigcup_{k=1}^{+\infty}\{\1\0^k,\2^k,\1\0\2^k\}$; hence $W_0'$ defined  to be either equal to $\0W_0$ or $\0\0W_0$, belongs to $\mathcal W$. If $W_0=\emptyword$ we put $W_0'=\emptyword$, so there exists $0\le j\le 2$ ~s.t. $\xi_{n-2,n+r}=\0^{m-n+2}\0^jW_0'W_1\cdots W_p$~that is:
$$
\quad \xi_{m,n+r}=\0^jW_0'W_1\cdots W_p.
$$
Condition ${\rm (U3.3)}$ holds by the following  lemma (obtained by a direct systematic verification) 
\begin{lemma}
$\forall w\in\{\0,\1,\2\}^*,\; U_1^\star A(w) \ne0$.
\end{lemma}
Condition ${\rm (U3.2)}$ follows from the following lemma, that can be proved in the same way as (\ref{universal}):
\begin{lemma}\label{ssss}
There exists a constant $\Lambda\ge1$ s.t. for any $0\le j\le 2$, $q\ge0$, $W_0''\in\{\emptyword\}\cup\mathcal{W}$ and $W_1'',\dots,W_q''\in\mathcal{W}$, one has $A(\0^jW_0''\cdots W_q'')\in\mathcal{H}_2(\Lambda)$. 
\end{lemma}

$\bullet$ The verification of   the conditions ${\rm (U2)}$ and ${\rm (U3)}$ about $\omega=w\bar\2$ is similar to the case of $\omega=w\bar\0$, using the fact that 
$$
\lim_{n\to+\infty}\left\Vert {A(\2^n)\over \Vert A(\2^n)\Vert}-C_\2U_5^\star\right\Vert=0,
\quad\text{where}\quad
C_\2=(1/3,0,1/3,1/3,0,0,0).
$$ 

The proof of Theorem~\ref{3matrices} is complete.

\hfill\qed

\newpage

\section{\bf Appendix: The adjacency graphs $(\Gamma1)$ and $(\Gamma2)$}\label{ADJACENCY}

\begin{figure}[H]
\includegraphics[scale=0.32]{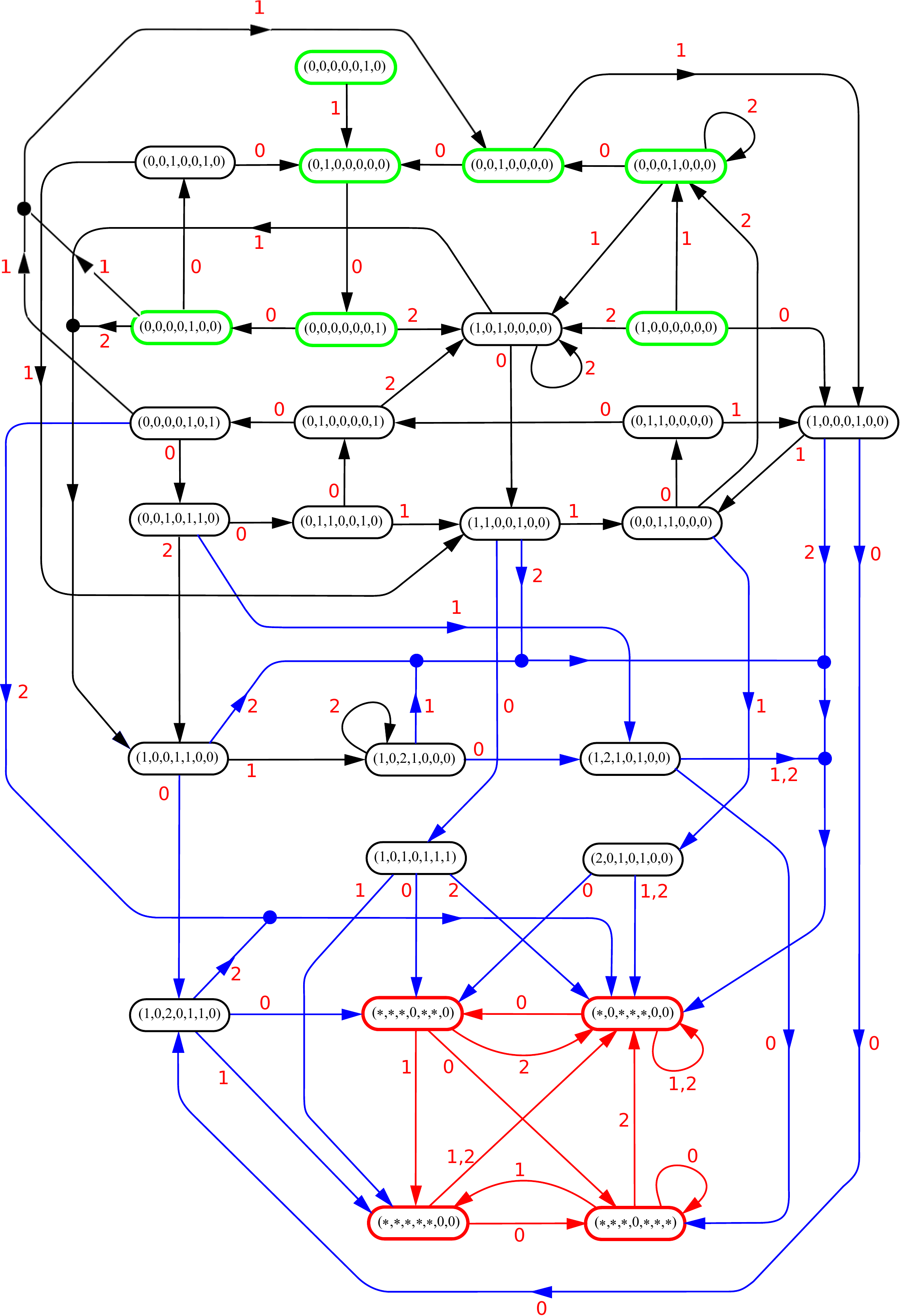}
\caption{\label{GRAPH1GOOD}\small\it A representation of $(\Gamma_1)$; for instance the vertex  
$
(1,0,2,1,0,0,0)$ corresponds to $\big\{(1,0,1,1,0,0,0),(1,0,2,1,0,0,0)\big\}
$
while 
$$
\!\!\!\!\!\!\!\!\!\!\!\!\!\!\!\!\!\!\!\!\!\!\!\!\!(*,0,*,*,*,0,0)=\big\{(x,0,y,z,t,0,0)\;;\;x,y,z,t=1,2,\dots\big\}\;;
$$
one has the partition ${\rm Vert}(\Gamma_1)=\mathcal{V}_1\sqcup\mathcal{V}_2$, where 
$$
\!\!\!\!\!\!\!\!\!\!\!\!\!\!\!\!\!\!\!\!\!\!\!\!\!\mathcal{V}_2:=\Big\{(*,0,*,*,*,0,0),\;(*,*,*,0,*,*,0),\;(*,*,*,0,*,*,*),\;(*,*,*,*,*,0,0)\Big\}\;;
$$
notice that  $(x_1,\dots,x_7)\in \mathcal{V}_1={\rm Vert}(\Gamma_1)\setminus\mathcal{V}_2$ implies $(x_1,\dots,x_7)\le(2,2,2,2,2,2,2)$.}
\end{figure}

\begin{figure}[H]
\includegraphics[scale=0.37]{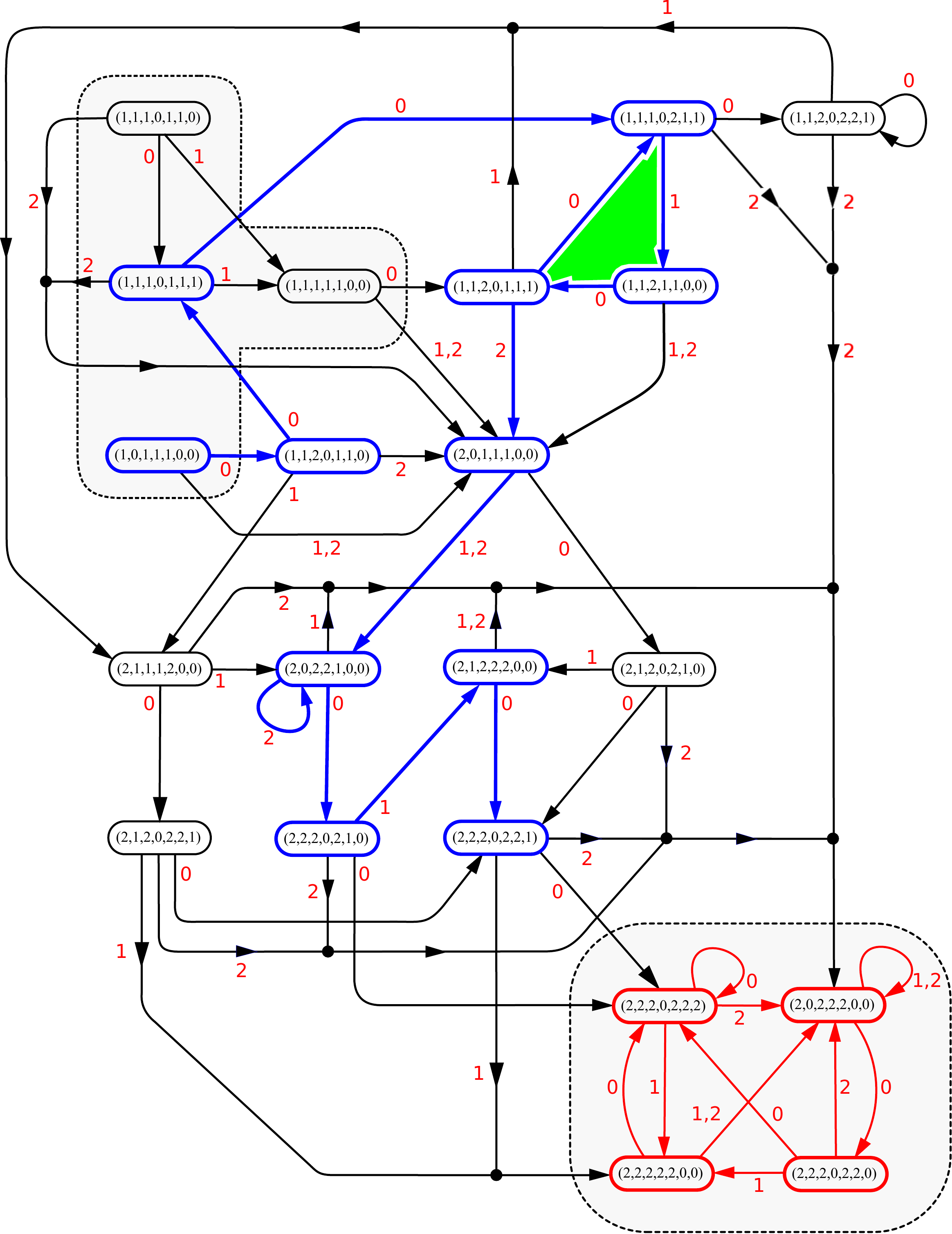}
\caption{\label{GRAPH2GOOD}\small\it Representation of the adjacency graph denoted by $(\Gamma_2)$ : the word-codes $W\in\{\0,\1,\2\}^*$ of the "longest" paths not entering $2\mathcal{T}_2$ (red vertices) and not factorized by $\1\0\0\1\0\0$ are of the form 
$W=\0\0\0\1\0\0\1\0\2\1\2^k\0\1\0$ ($k\ge1$);
they are represented in blue in the adjacency graph. The green triangle shows the  "forbidden" cycle $\1\0\0\1\0\0\1\0\0\cdots.$}
\end{figure}

\newpage

\section{\bf Appendix: Verification for matrix products in Lemma~\ref{13}}\label{H1}

\subsection{First case}\label{firstcaseofAppendix}We check that $A(W)$ satisfies  parts (i-ii-iii) of Lemma~\ref{13}, for any $W\in\mathcal{L'}\0\1\mathcal{L}''$. 

$\bullet$ for $w\in \0\0\0\1\{\0\0,\0\1,\0\2,\1\0\}$

\begin{align*}
A(\0\0\0\1\0\0)
=&\begin{pmatrix}1&0&0&0&0&0&1\\1&0&0&0&0&0&1\\2&0&0&0&0&0&2\\0&0&0&0&0&0&0\\2&1&0&0&0&0&1\\2&0&0&0&0&0&2\\2&0&0&0&0&0&1\end{pmatrix}
\quad
A(\0\0\0\1\0\1)
=\begin{pmatrix}1&0&1&0&0&0&0\\1&0&1&0&0&0&0\\1&0&2&0&0&0&0\\0&0&0&0&0&0&0\\1&0&2&1&0&0&0\\1&0&2&0&0&0&0\\1&0&2&1&0&0&0\end{pmatrix}\\
A(\0\0\0\1\0\2)
=&\begin{pmatrix}0&0&0&1&2&0&0\\0&0&0&1&2&0&0\\0&0&0&1&3&0&0\\0&0&0&0&0&0&0\\1&0&0&1&3&0&1\\0&0&0&1&3&0&0\\1&0&0&1&3&0&1\end{pmatrix}
\quad
A(\0\0\0\1\1\0)
=\begin{pmatrix}2&0&0&0&0&0&1\\2&0&0&0&0&0&1\\2&0&0&0&0&0&1\\0&0&0&0&0&0&0\\3&0&0&1&1&0&1\\2&0&0&0&0&0&1\\1&0&0&1&1&0&1\end{pmatrix}
 \end{align*}

$\bullet$ $\0\0\0\1\1\1$ is a prefix of each word in $\0\0\0\1\{\1\1\0,\1\1\1,\1\1\2\}$ and 

$$
A(\0\0\0\1\1\1)
=\begin{pmatrix}1&0&2&1&0&0&0\\1&0&2&1&0&0&0\\1&0&2&1&0&0&0\\0&0&0&0&0&0&0\\1&0&3&3&1&0&0\\1&0&2&1&0&0&0\\2&0&1&1&1&0&0\end{pmatrix}
$$

$\bullet$ for $w=\0\0\0\1\1\2$, one has
$$
A(\0\0\0\1\1\2)
=\begin{pmatrix}1&0&0&1&3&0&1\\1&0&0&1&3&0&1\\1&0&0&1&3&0&1\\0&0&0&0&0&0&0\\3&0&0&1&4&0&3\\1&0&0&1&3&0&1\\1&0&0&2&3&0&1\end{pmatrix}
$$

$\bullet$ $\0\0\0\1\2$ is a prefix of each word in $\0\0\0\1\{\2\0,\2\1,\2\2\}$ and 
$$
A(\0\0\0\1\2)
=\begin{pmatrix}1&0&0&1&1&0&1\\1&0&0&1&1&0&1\\1&0&0&1&1&0&1\\0&0&0&0&0&0&0\\1&0&0&2&3&0&1\\1&0&0&1&1&0&1\\2&0&0&0&1&0&2\end{pmatrix}
$$

$\bullet$ $\1\0\0\1$ is prefix of each word in $\1\0\0\1\mathcal{L}''$ and 
$$
A(\1\0\0\1)
=\begin{pmatrix}0&0&1&1&0&0&0\\0&0&1&1&0&0&0\\0&0&1&1&0&1&0\\0&0&1&1&0&0&0\\0&0&1&1&0&0&0\\0&0&0&0&0&0&0\\0&0&0&0&0&0&0\end{pmatrix}
$$

$\bullet$ each word in $\2\0\0\1\mathcal{L}''$ has a prefix in $\2\0\0\1\{\0,\1,\2\}$
\begin{align*}
A(\2\0\0\1\0)
=&\begin{pmatrix}1&0&0&2&3&0&1\\0&0&0&0&0&0&0\\1&0&0&1&1&0&1\\0&0&0&1&2&0&0\\0&0&0&1&2&0&0\\0&0&0&0&0&0&0\\0&0&0&0&0&0&0\end{pmatrix}
\quad 
A(\2\0\0\1\1)
=\begin{pmatrix}3&0&1&2&2&0&0\\0&0&0&0&0&0&0\\2&0&1&1&1&0&0\\1&0&0&1&1&0&0\\1&0&0&1&1&0&0\\0&0&0&0&0&0&0\\0&0&0&0&0&0&0\end{pmatrix}\\
A(\2\0\0\1\2)
=&\begin{pmatrix}2&0&0&3&4&0&2\\0&0&0&0&0&0&0\\1&0&0&2&3&0&1\\1&0&0&1&1&0&1\\1&0&0&1&1&0&1\\0&0&0&0&0&0&0\\0&0&0&0&0&0&0\end{pmatrix}
\end{align*}

$\bullet$ $\1\0\1\0$ is a prefix of each word in $\1\0\1\{\0\0,\0\1,\0\2\}$ and 
$$
A(\1\0\1\0)
=\begin{pmatrix}1&0&0&1&1&0&0\\0&0&0&1&1&0&0\\0&0&0&1&1&0&0\\0&0&0&1&1&0&0\\1&0&0&1&1&0&0\\0&0&0&0&0&0&0\\0&0&0&0&0&0&0\end{pmatrix}
$$

$\bullet$ for the remaining words in $\1\0\1\mathcal{L}''$ one has
\begin{align*}
A(\1\0\1\1\0)
=&\begin{pmatrix}1&0&0&1&1&0&1\\1&0&0&0&0&0&1\\2&0&0&0&0&0&1\\2&0&0&0&0&0&1\\1&0&0&1&1&0&1\\0&0&0&0&0&0&0\\0&0&0&0&0&0&0\end{pmatrix}\quad
A(\1\0\1\1\1\0)
=\begin{pmatrix}3&0&0&1&1&0&1\\1&0&0&1&1&0&0\\1&0&0&2&2&0&0\\1&0&0&2&2&0&0\\3&0&0&1&1&0&1\\0&0&0&0&0&0&0\\0&0&0&0&0&0&0\end{pmatrix}\\
A(\1\0\1\1\1\1)
=&\begin{pmatrix}1&0&3&3&1&0&0\\0&0&1&2&1&0&0\\1&0&1&3&2&0&0\\1&0&1&3&2&0&0\\1&0&3&3&1&0&0\\0&0&0&0&0&0&0\\0&0&0&0&0&0&0\end{pmatrix}\quad
A(\1\0\1\1\1\2)
=\begin{pmatrix}3&0&0&1&4&0&3\\2&0&0&0&1&0&2\\3&0&0&1&2&0&3\\3&0&0&1&2&0&3\\3&0&0&1&4&0&3\\0&0&0&0&0&0&0\\0&0&0&0&0&0&0\end{pmatrix}
\end{align*}
\begin{align*}
A(\1\0\1\1\2)
=&\begin{pmatrix}1&0&0&2&3&0&1\\0&0&0&1&2&0&0\\1&0&0&1&3&0&1\\1&0&0&1&3&0&1\\1&0&0&2&3&0&1\\0&0&0&0&0&0&0\\0&0&0&0&0&0&0\end{pmatrix}\quad
A(\1\0\1\2\0)
=\begin{pmatrix}3&2&0&0&0&0&1\\1&1&0&0&0&0&0\\2&1&0&0&0&0&1\\2&1&0&0&0&0&1\\3&2&0&0&0&0&1\\0&0&0&0&0&0&0\\0&0&0&0&0&0&0\end{pmatrix}\\
A(\1\0\1\2\1)
=&\begin{pmatrix}0&0&3&2&0&0&0\\0&0&1&1&0&0&0\\1&0&2&1&0&0&0\\1&0&2&1&0&0&0\\0&0&3&2&0&0&0\\0&0&0&0&0&0&0\\0&0&0&0&0&0&0\end{pmatrix}\quad
A(\1\0\1\2\2)
=\begin{pmatrix}2&0&0&0&3&0&2\\1&0&0&0&1&0&1\\1&0&0&1&3&0&1\\1&0&0&1&3&0&1\\2&0&0&0&3&0&2\\0&0&0&0&0&0&0\\0&0&0&0&0&0&0\end{pmatrix}
\end{align*}

$\bullet$ $\2\0\1$ is a prefix of each word in $\2\0\1\mathcal{L}''$ and 
$$
A(\2\0\1)
=\begin{pmatrix}0&0&2&2&0&1&0\\0&0&0&0&0&0&0\\0&0&1&1&0&1&0\\0&0&1&1&0&0&0\\0&0&1&1&0&0&0\\0&0&0&0&0&0&0\\0&0&0&0&0&0&0\end{pmatrix}
$$
\subsection{Second case}\label{secondcaseofAppendix}To verify that $A(W)$ satisfies (i-ii-iii), for any $W\in\mathcal{L}'\2\1\mathcal{L}'''$, notice that such a $W$ has a suffix in $\2\1\mathcal{L}'''$ and moreover,
\begin{align*}
A(\2\1\0)
=&\begin{pmatrix}0&0&0&2&2&0&0\\0&0&0&0&0&0&0\\0&0&0&1&1&0&0\\1&0&0&1&1&0&0\\0&0&0&1&1&0&0\\0&0&0&0&0&0&0\\0&0&0&0&0&0&0\end{pmatrix}
\quad 
A(\2\1\1\0)
=\begin{pmatrix}3&0&0&0&0&0&2\\0&0&0&0&0&0&0\\2&0&0&0&0&0&1\\1&0&0&1&1&0&1\\1&0&0&0&0&0&1\\0&0&0&0&0&0&0\\0&0&0&0&0&0&0\end{pmatrix}\\
A(\2\1\1\1\0)
=&\begin{pmatrix}2&0&0&3&3&0&0\\0&0&0&0&0&0&0\\1&0&0&2&2&0&0\\3&0&0&1&1&0&1\\1&0&0&1&1&0&0\\0&0&0&0&0&0&0\\0&0&0&0&0&0&0\end{pmatrix}
\quad 
\!\!\!A(\2\1\1\1\1)
=\begin{pmatrix}1&0&2&5&3&0&0\\0&0&0&0&0&0&0\\1&0&1&3&2&0&0\\1&0&3&3&1&0&0\\0&0&1&2&1&0&0\\0&0&0&0&0&0&0\\0&0&0&0&0&0&0\end{pmatrix}\\
A(\2\1\1\1\2)
=&\begin{pmatrix}5&0&0&1&3&0&5\\0&0&0&0&0&0&0\\3&0&0&1&2&0&3\\3&0&0&1&4&0&3\\2&0&0&0&1&0&2\\0&0&0&0&0&0&0\\0&0&0&0&0&0&0\end{pmatrix}\quad
A(\2\1\1\2)
=\begin{pmatrix}1&0&0&2&5&0&1\\0&0&0&0&0&0&0\\1&0&0&1&3&0&1\\1&0&0&2&3&0&1\\0&0&0&1&2&0&0\\0&0&0&0&0&0&0\\0&0&0&0&0&0&0\end{pmatrix}
\end{align*}
\begin{align*}
A(\2\1\2\0)
=&\begin{pmatrix}3&2&0&0&0&0&1\\0&0&0&0&0&0&0\\2&1&0&0&0&0&1\\3&2&0&0&0&0&1\\1&1&0&0&0&0&0\\0&0&0&0&0&0&0\\0&0&0&0&0&0&0\end{pmatrix}\quad
A(\2\1\2\1)
=\begin{pmatrix}1&0&3&2&0&0&0\\0&0&0&0&0&0&0\\1&0&2&1&0&0&0\\0&0&3&2&0&0&0\\0&0&1&1&0&0&0\\0&0&0&0&0&0&0\\0&0&0&0&0&0&0\end{pmatrix}\\
A(\2\1\2\2)
=&\begin{pmatrix}2&0&0&1&4&0&2\\0&0&0&0&0&0&0\\1&0&0&1&3&0&1\\2&0&0&0&3&0&2\\1&0&0&0&1&0&1\\0&0&0&0&0&0&0\\0&0&0&0&0&0&0\end{pmatrix}
\end{align*}

\subsection{Third case}\label{thirdcaseofAppendix} The matrix $A(W)$ satisfies (i-ii-iii) for $W\in\{\2\0\2\0,\2\0\0\2\0,\0\0\0\2\0,\2\2\0\}$, since:

\begin{align*}
\quad
A(\2\0\2\0)
=&\begin{pmatrix}4&2&0&0&0&0&2\\0&0&0&0&0&0&0\\2&1&0&0&0&0&1\\2&1&0&0&0&0&1\\2&1&0&0&0&0&1\\0&0&0&0&0&0&0\\0&0&0&0&0&0&0\end{pmatrix}\quad
\!\!\!A(\2\0\0\2\0)
=\begin{pmatrix}5&3&0&0&0&0&2\\0&0&0&0&0&0&0\\3&2&0&0&0&0&1\\2&1&0&0&0&0&1\\2&1&0&0&0&0&1\\0&0&0&0&0&0&0\\0&0&0&0&0&0&0\end{pmatrix}\\
A(\0\0\0\2\0)
=&\begin{pmatrix}2&1&0&0&0&0&1\\2&1&0&0&0&0&1\\2&1&0&0&0&0&1\\0&0&0&0&0&0&0\\3&2&0&0&0&0&1\\2&1&0&0&0&0&1\\2&0&0&0&0&0&2\end{pmatrix}\;\;\quad 
A(\2\2\0)
=\begin{pmatrix}3&1&0&0&0&0&2\\0&0&0&0&0&0&0\\2&1&0&0&0&0&1\\2&0&0&0&0&0&2\\1&0&0&0&0&0&1\\0&0&0&0&0&0&0\\0&0&0&0&0&0&0\end{pmatrix}.
\end{align*}


\end{document}